\documentclass[12pt]{book}

\usepackage{amsfonts}
\usepackage{amsmath}
\usepackage{amsthm}
\usepackage{amssymb}
\usepackage[english]{babel}
\usepackage{graphics}
\usepackage{latexsym}
\usepackage{longtable} \usepackage{mathrsfs}
\usepackage{graphicx}
\usepackage{wrapfig} 
\usepackage[pdftex,colorlinks=true]{hyperref}

\newtheorem{theorem}{Theorem}
\newtheorem{corollary}[theorem]{Corollary}
\newtheorem{lemma}[theorem]{Lemma}

\newtheorem{claim}[theorem]{Claim}
\newtheorem{example}[theorem]{Example}

\theoremstyle{definition}
\newtheorem{definition}[theorem]{Definition}

\newcommand{\mU}{\mathcal{U}}

\newcommand{\mO}{\mathcal{O}}
\renewcommand{\S}{\mathcal{S}}
\newcommand{\B}{\mathcal{B}}
\newcommand{\J}{\mathcal{J}}

\newcommand{\R}{\mathbb{R}}
\newcommand{\N}{\mathbb{N}}
\newcommand{\mS}{\mathbb{S}}
\newcommand{\mB}{\mathbb{B}}

\newcommand{\noi}{\noindent}
\newcommand{\ms}{\medskip}

\newcommand{\al}{\alpha}
\newcommand{\be}{\beta}
\newcommand{\ga}{\gamma}
\newcommand{\de}{\delta}
\newcommand{\De}{\Delta}
\newcommand{\e}{\varepsilon}
\newcommand{\si}{\sigma}
\newcommand{\la}{\lambda}

\newcommand{\Om}{\Omega}
\newcommand{\om}{\omega}

\newcommand{\lharpoonup}{-\!\!\!\!\rightharpoonup}

\newcommand{\larrow}{\longrightarrow}
\newcommand{\ot}{\otimes}

\newcommand{\ri}{\rightarrow}

\newcommand{\p}{\partial}
\newcommand{\sub}{\subseteq}
\newcommand{\set}{\setminus}
\newcommand{\by}{\times}

\newcommand{\Lip}{\textrm{Lip}}
\newcommand{\tr}{\textrm{tr}}

\newcommand{\diam}{\textrm{diam}}
\newcommand{\dist}{\textrm{dist}}

\newcommand{\Div}{\textrm{Div}}

\newcommand{\bt}{\begin{theorem}}\newcommand{\et}{\end{theorem}}
\newcommand{\bd}{\begin{definition}}\newcommand{\ed}{\end{definition}}
\newcommand{\bl}{\begin{lemma}}\newcommand{\el}{\end{lemma}}
\newcommand{\beq}{\begin{equation}}\newcommand{\eeq}{\end{equation}}
\newcommand{\bc}{\begin{claim}}\newcommand{\ec}{\end{claim}}
\newcommand{\bex}{\begin{example}}\newcommand{\eex}{\end{example}}
\newcommand{\bcor}{\begin{corollary}}\newcommand{\ecor}{\end{corollary}}
\newcommand{\bp}{\begin{proof}}\newcommand{\ep}{\end{proof}}

\numberwithin{equation}{chapter}

\begin{document}

\title{ AN\\
INTRODUCTION TO\\
 VISCOSITY SOLUTIONS FOR FULLY NONLINEAR PDE\\ WITH\\ APPLICATIONS TO CALCULUS OF VARIATIONS IN $L^\infty$\\
 \ms\ms\ms
  \tiny{\textbf{This is a preprint. The final publication is available at Springer via http://dx.doi.org/10.1007/978-3-319-12829-0}}}

\author{\textsl{Nikos Katzourakis}\ms\\
Department of Mathematics and Statistics,\\
University of Reading, Reading, UK.
\ms\\
\textit{n.katzourakis@reading.ac.uk}}

\maketitle




\date{ }


\[
\phantom{L}
\]
\[
\phantom{L}
\]
\[
\phantom{L}
\]
\[
\phantom{L}
\]

\[
\text{\hspace{280pt} To my beloved wife}
\]
\maketitle

\newpage

\tableofcontents

\chapter*{Preface}

This set of notes corresponds to the lectures of a post-graduate short course given by the author at the BCAM - Basque Centre for Applied Mathematics at Bilbao of Spain in early May 2013.

The purpose of these notes is to give a quick and elementary, yet rigorous, presentation of the rudiments of the so-called theory of Viscosity Solutions which applies to fully nonlinear 1st and 2nd order Partial Differential Equations (PDE). For such equations, particularly for 2nd order ones, solutions generally are nonsmooth and standard approaches in order to define a ``weak solution" do not apply: classical, strong almost everywhere (a.e.), weak, measure-valued and distributional solutions either do not exist or may not even be defined. The main reason for the latter failure is that, the standard idea of  using ``integration-by-parts" in order to pass derivatives to smooth test functions by duality, is not available for non-divergence structure PDE. 

The name of this theory originates from the ``vanishing viscosity method" developed first for 1st order fully nonlinear PDE (Hamilton-Jacobi PDE). Today, though, it comprises an independent theory of ``weak" solutions which applies to fully nonlinear elliptic and parabolic PDE and in most cases has no or little relation to the idea of adding a viscosity term. The formal notions have been introduced by P.L. Lions and M.G. Crandall in the early 1980s for 1st order PDE, following preceding contributions of L.C. Evans. The extension to the case of 2nd order PDEs came around the 1990s by H. Ishii and P.L. Lions.

Interesting PDE to which the theory applies arise in Geometry and Geometric Evolution (Monge-Amper\'e PDE, Equations of Motion by Mean Curvature), Optimal Control and Game Theory (Hamilton-Jacobi-Bellman PDE, Isaacs PDE, Differential Games) and Calculus of Variations in $L^p$ and $L^\infty$ (Euler-Lagrange PDE, $p$-Laplacian, Aronsson PDE, $\infty$-Laplacian). 

Due to the vastness of the subject, a drastic choice of material is required for a brief and elementary introduction to the subject of Viscosity Solutions. In the case at hand, our main criterion has been nothing more but personal taste. Hence, herein, we shall restrict ourselves to the 2nd order degenerate elliptic case, focusing in particular on applications in the modern field of Calculus of Variations in $L^\infty$. 

An inspection of the Table of Contents gives an idea about the organisation of the material. A rather immediate observation of the expert is that, unlike most standard texts on Viscosity Solutions where Uniqueness and Comparison form the centre of gravity of the exposition, herein they appear rather late in the presentation and are not over-emphasised. This shift of viewpoint owes to that the author is directed mostly towards extensions of the theory to the vector case of systems. In this realm, the primary focus switches to existence methods, while comparison and uniqueness are not true in general, not even in the most ideal cases.

Throughout these notes, no previous knowledge is assumed on the reader's behalf. Basic graduate-level mathematical maturity suffices for the reading of the first six chapters which is the general theory. For the next two chapters which concern applications to Calculus of Variations, some familiarity with weak derivatives and functionals is assumed, which does not go much deeper than the definitions. The last chapter collects, mostly without proofs, extensions and perhaps more advanced related topics which we were not able to cover in detail in this introduction to the subject. 

There exist several excellent expository texts in the literature on Viscosity Solutions and Calculus of Variations in $L^\infty$. In particular, we would like to point out \cite{CIL, C2, L1, Ko, C1, B, G, MR, Dr}. However, all the sources we are aware of are either of more advanced level, or have a different viewpoint. We believe that the main contribution of the present text is that is \emph{elementary and is mostly addressed to students and non-experts}. This ``textbook style" is reflected also to the fact that references do not appear in the main text. Most of the results are not optimal, but instead there has been a huge, and hopefully successful effort the main ideas to be illustrated as clearly as possible. We hope that these notes serve as a suitable \emph{first reading} on the theory of Viscosity Solutions.

\ms 
\ms

\ms 
\ms

\hspace{280pt} N. K.,

\hspace{280pt}  Reading, UK

\hspace{280pt}  May, 2014

\newpage

$$ $$

$$ $$

\noi \textbf{Acknowledgement.} I am indebted to Enrique Zuazua for inviting me to write this Brief for Springer. Special thanks are due to Juan Manfredi for his share of expertise on the subject of the present notes during several personal communications. I would like to thank Tristan Pryer for the careful reading of an earlier version of this manuscript. His suggestions certainly made the material more readable. I also thank Federica Dragoni warmly for her remarks and suggestions. Last but not least, I am grateful to the referees of this monograph whose numerous suggestions improved the quality of the presentation.

\chapter*{General Theory}

\chapter[History, Examples, Motivation, ...]{History, Examples, Motivation and First Definitions}

What is it all about?

{
\center{\fbox{\parbox[pos]{320pt}{
\ms
\emph{Viscosity Solutions} form a general theory of ``weak" (i.e.\ non-differentiable) solutions which applies to certain fully nonlinear Partial Differential Equations (PDE) of 1st and 2nd order.

\ms
}}

\ms
\ms
}}

\noi Let $u:\Om \sub \R^n \larrow \R$ be a function in $C^2(\Om)$, $n\geq 1$. In the standard way, $C^k(\Om)$ denotes the space of $k$-times continuously differentiable functions over the domain $\Om$. 

Consider the PDE

\beq \label{eq1}
F(\cdot , u,Du,D^2u)\, =\, 0,
\eeq
where
\[
Du\, =\, (D_1u,...,D_n u),\ \ D_i \equiv \frac{\p}{\p x_i},\ i=1,...,n \text{ (Gradient vector)},
\]
\[
D^2u\,=\, \big(D^2_{ij}u \big)_{j=1,...,n}^{i=1,...,n}\  \text{ (Hessian matrix)},
\]
and 
\[
F\ :\ \Om \by \R \by \R^n \by \mS(n) \larrow \R 
\]
is the function defining the equation\footnote{It is customary to call $F$ the ``coefficients" and we will occasionally follow this convention. This terminology is inherited from the case of linear equations, where $F$ consists of functions multiplying the solution and its derivatives. We will also consistently use the notation ``$\cdot $" for the argument of $x$, namely $F(\cdot, u,Du,D^2u)|_{x}=F\big(x,u(x),Du(x),D^2u(x)\big)$. We will \textit{not use} the common clumsy notation $F(x, u,Du,D^2u)$.} whose arguments are denoted by
\[
F(x,r,p,X), \ \ \ (x,r,p,X) \in Om \by \R \by \R^n \by \mS(n).
\]
$F$ is always assumed to be continuous, unless stated otherwise. The set $\mS(n) $ denotes the symmetric $n \by n$ matrices:
\[
\mS(n)\, :=\, \Big\{ X=(X_{ij})_{j=1,...,n}^{i=1,...,n} \in \R^{n \by n} \ \Big|\ X_{ij}=X_{ji} \Big\}. 
\]
The 2nd order PDE \eqref{eq1} trivially includes the case of fully nonlinear 1st order PDE (Hamilton-Jacobi PDE):
\beq \label{eq2}
F(\cdot , u,Du)\, =\, 0.
\eeq
The PDE \eqref{eq1} (and hence \eqref{eq2}) is our \emph{basic object of study}. As usually, here and subsequently,

{\center{

\fbox{\parbox[pos]{340pt}{
\ms

the modifier \textbf{``fully nonlinear"} stands for \textbf{``$F$ \emph{may} not linear in any of its arguments"}, but of course linear equations are not excluded. We will say that an equation is \textbf{genuinely} fully nonlinear, if is nonlinear with respect to the highest order derivatives. 

\ms
}}

\ms\ms }}

\noi If the partial function
\[
X\, \mapsto \, F(x,r,p,X)
\]
is linear, the PDE is typically called \emph{quasilinear} and hence $F$ can be written in the form
\beq \label{eq3}
F(x,r,p,X)\, =\, A_{ij}(x,r,p) X_{ij}\, +\, B(x,r,p)
\eeq
for a matrix function
\[
A\ :\ \Om \by \R \by \R^n \larrow \R^{n \by n}
\]
and a scalar function $B : \Om \by \R \by \R^n \larrow \R$. The summation in $i,j$ in \eqref{eq3} is tacitly assumed without explicitly writing the sum (``Einstein's convention"). By denoting the Euclidean inner product of the matrix space $\R^{n \by n}$ by ``$:$", that is
\[
A:B\, :=\, A_{ij}B_{ij}\, = \, \tr( A^\top\! B) ,
\]
the general quasilinear PDE can be written as
\beq \label{eq4}
A(\cdot, u,Du):D^2u \, +\, B(\cdot,u,Du)\, =\, 0.
\eeq

\noi \textbf{Historical comment:} 

{\center{

\fbox{\parbox[pos]{320pt}{
\ms
Viscosity Solutions were first introduced in the 1980s by Crandall and Lions for $F(\cdot , u,Du)=0$ as a
\ms

\centerline{\emph{uniqueness criterion}, }
\ms

in order to select one of the infinitely-many strong a.e.\ Lipschitz solutions of the Dirichlet problem for $F(\cdot , u,Du)=0$. 
\ms
}}

\ms\ms }}

\ms
\noi \textbf{Example 1 (Non-uniqueness of strong solutions).} The Dirichlet problem for the $1$-dimensional Eikonal PDE:
\[
\left\{
\begin{array}{l}
|u'|^2-1=0, \text{ on }(-1,+1), \ms\\
u(\pm 1)=0,
\end{array}
\right.
\]
admits infinitely many Lipschitz continuous solutions $u:(-1,+1)\ri \R$, with $|u'|^2=1$ a.e.\ on $(-1,+1)$. Indeed, one of these solutions is the distance from the endpoints of the interval:
\[
u_1(x):= 
\left\{
\begin{array}{l}
x-1, \text{ on }(-1,0), \ms\\
1-x, \text{ on }(0,1),
\end{array}
\right.
\]
which solves $|u'|^2=1$ on $(-1,0)\cup(0,1)$. Then, by reflection of the graph of $u_1$ with respect to the horizontal line $u=1/2$, the function
\[
u_2(x):= 
\left\{
\begin{array}{l}
x-1, \text{ on }(-1,-1/2), \ms\\
\ \ \ -x, \text{ on }(-1/2,0), \ms\\
\ \ \ \ \ x, \text{ on }(0,1/2), \ms\\
1-x, \text{ on }(0,1),
\end{array}
\right.
\]
is also a solution of $|u'|^2=1$ on $(-1,-1/2)\cup(-1/2,0)\cup(0,1/2)\cup(1/2,1)$. By reflection of the graph of $u_2$ with respect to the line $u=1/4$, we construct a new solution $u_3$, and so on.
\[
\includegraphics[scale=0.22]{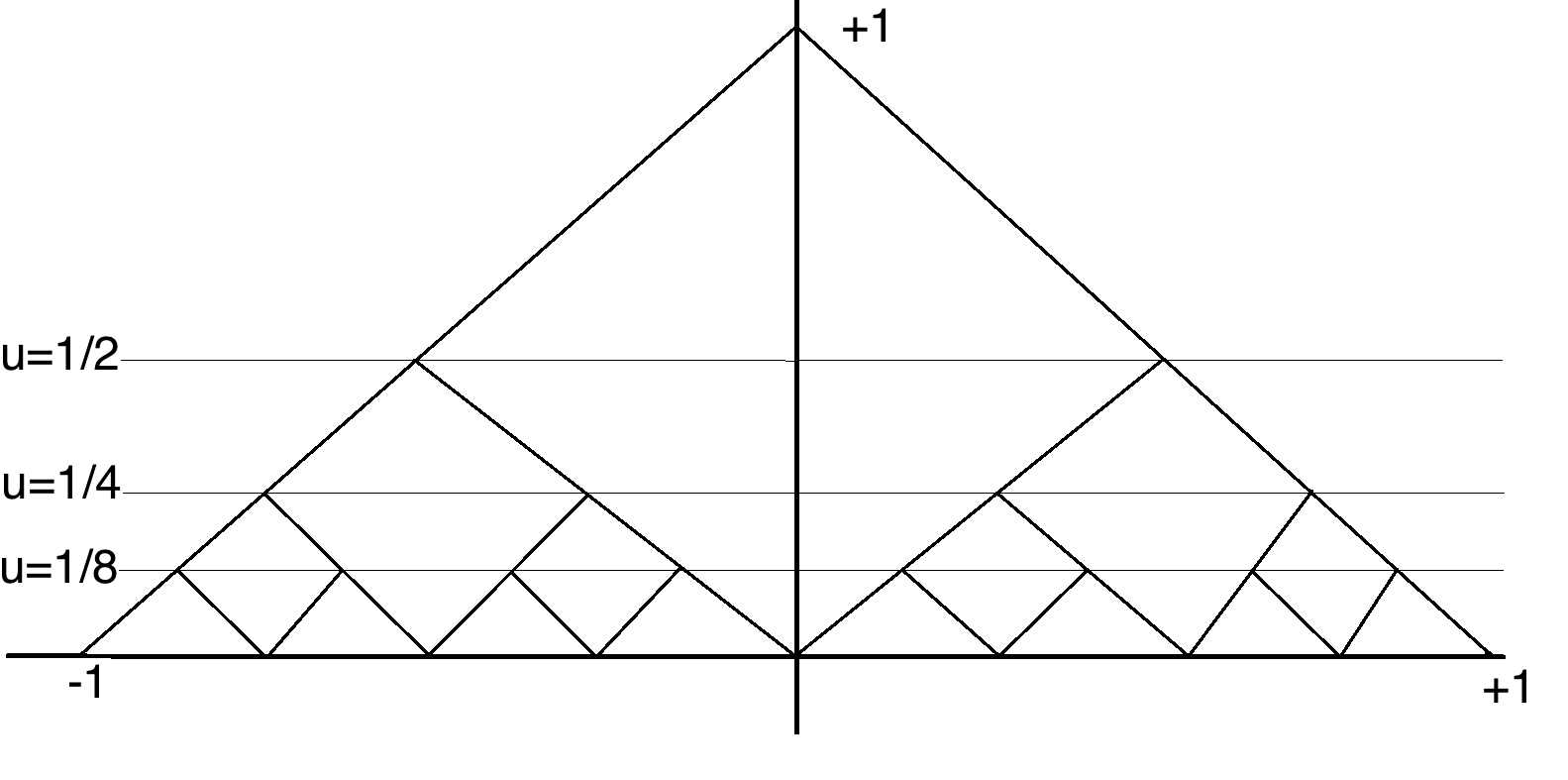}
\]

\noi In many problems, strong a.e.\ solutions of Hamilton-Jacobi PDE were constructed by utilising the \emph{Vanishing Viscosity method}.

\ms

\noi  \textit{\underline{The idea of the Vanishing Viscosity method}}: In order to solve the PDE
\[
F(\cdot , u,Du)\, =\, 0,
\]
we approximate by solutions of the elliptic PDE
\[
\e \De u_\e \, +\, F(\cdot , u_\e,Du_\e)\, =\, 0,
\]
and investigate the limit of $(u_\e)_{\e>0}$ as $\e \ri 0$. The hope then is that we can obtain some sort of estimates which do not collapse as $\e \ri 0$ that will allow to pass to a limit $u_\e \ri u$ which ``solves" the original PDE (``$\De$" of course stands for the Laplacian).

{\center{

\fbox{\parbox[pos]{320pt}{
\ms
Viscosity Solutions originate from the ``Vanishing Viscosity method", but should not be confused with the latter!!!

\ms
In most cases, especially for 2nd order PDE, Viscosity Solutions have nothing to do with adding an ``$\e \De u$" term.

\ms
}}

\ms\ms }}

\noi \textbf{Historical comment:} 

{\center{

\fbox{\parbox[pos]{320pt}{
\ms
For 2nd order PDE, problems more fundamental than uniqueness arise: ``singular solutions" appear, which can not be rigorously justified, because 

\ms
\centerline{$F(\cdot , u ,Du, D^2 u)=0$}
\ms

does not support integration by parts!!!
\ms
}}

\ms\ms }}

\ms

\noi The meaning of the above statement is that we can not ``pass derivatives to test functions" by integration by parts in order to have, for example, a function which is only (weakly) once differentiable as solution of a 2nd order equation.
\ms

\noi \textbf{Remark 2 (Digression into divergence structure PDE).} For the 2nd order PDE
\beq \label{eq5}
D_i \big(A_i(\cdot, u,Du) \big)\, +\, B(\cdot,u,Du)\, =\, 0,
\eeq
where $u:\Om \sub \R^n \ri \R$, we have the option to define weak solutions by duality, that is by taking ``integration by parts" as definition when testing against smooth functions with compact support:
\[
\int_\Om -A_i(\cdot, u,Du) \big)D_i\psi  + B(\cdot,u,Du)\psi \, =\, 0,
\]
for all $\psi \in C^\infty_c(\Om)$.

In particular, divergence structure PDE are always quasilinear and \emph{not genuinely} fully nonlinear. Indeed, by assuming that $A,B$ are $C^1$ and distributing derivatives in \eqref{eq5}, we have
\begin{align}  \label{eq6}
{A_i}_{p_j}(\cdot, u,Du) D^2_{ij}u & \, +\, \Big\{{A_i}_r(\cdot, u,Du)D_i u \nonumber\\
&+\, {A_i}_{x_i}(\cdot, u,Du) + B(\cdot, u,Du)\Big\}\, =\,0.
\end{align}
Hence, the dependence in $D^2u$ is linear. In \eqref{eq6}, the subscripts $p_j, r, x_i$ denote differentiation of the coefficients $A,B$ with respect to the respective arguments.

\ms
\noi \textbf{Remark 3.}

{\center{
\fbox{\parbox[pos]{240pt}{
\ms

\hspace{15pt }Viscosity Solutions form a duality-free 

\ms
``NONLINEAR DISTRIBUTION THEORY", 

\ms
whose ideas are needed even for linear PDE!!!
\ms
}}

}}

\ms\ms
\noi \textbf{Example 4:} Let $A=( A_{ij})^{j=1,...,n}_{i=1,...,n} : \Om \sub \R^n \larrow \mS(n)$  be a map with values non-negative matrices, such that either $A$ has a vanishing eigenvalue, or $A$ is not locally Lipschitz continuous, then the linear PDE
\[
A_{ij}(x)D^2_{ij}u(x)\, =\, 0,
\]
may not have any non-trivial solutions (the affine ones) in all standard senses! Explicit $A$'s will be given later.

On the other hand, standard elliptic estimates (Schauder) imply that if $A$ is both H\"older continuous $C^\al$ and strictly positive, then the PDE has only smooth $C^{2,\al}$ solutions. 

This problem is the same as the problem to interpret (and solve) rigorously the equation
\[
H'\,=\, \de,
\] 
where $\de$ is the Dirac ``function" and $H=\chi_{(0,\infty)}$ is a discontinuous solution, without having Distributions (Generalised Functions) and Measures at hand.
\ms

\noi However, there are \emph{limitations}: Viscosity Solutions apply to \emph{degenerate elliptic} and \emph{degenerate parabolic} PDE. The appropriate ellipticity notion will be given after we will many classes of equations to which our theory applies.

\ms
\noi \textbf{Example 5 (PDE to which Viscosity Solutions apply).} 

\ms

\noi (1) \textbf{An all-important equation: the $\infty$-Laplacian}. For $u\in C^2(\Om)$ and $\Om\sub \R^n$, the infinity-Laplacian is 
\beq \label{eq7}
\De_\infty u \, :=\, Du \ot Du : D^2u\,=\, 0,
\eeq
where the expression in \eqref{eq7} is understood as the double sum $D_i u\,  D_j u \, D^2_{ij}u$. The geometric meaning of \eqref{eq7} is that the rank-one tensor product matrix $Du \ot Du$ is pointwise normal to $D^2u$ in $\R^{n\by n}$. \eqref{eq7} is the fundamental PDE of \emph{Calculus of Variations in $L^\infty$} (the analogue of the Euler-Lagrange PDE), when considering variational problems for the so-called Supremal functionals, the simplest of which is
\beq  \label{eq8}
E_\infty(u,\Om')\, :=\, \|Du\|_{L^\infty(\Om')},\ \ \Om' \Subset \Om.
\eeq
Here $\Om'$ is an open set compactly contained in $\Om$ and $u$ is locally Lipschitz continuous, that is, in the Sobolev space 
\[
W^{1,\infty}_{\text{loc}}(\Om)\, :=\, \Big\{ u \in L^\infty_{\text{loc}}(\Om) \ \Big|\ \exists\,  Du \text{ a.e.\ on }\Om \text{ and } Du \in L^\infty_{\text{loc}}(\Om)^n\Big\}.
\]
The study of such problems has been initiated by Aronsson in the 1960s, but their rigorous study started in the 1990s by utilising the theory of Viscosity Solutions. Both \eqref{eq7} and \eqref{eq8} arise also in the geometric problem of Optimisation of Lipschitz Extensions of functions, Differential games, (implicitly) in Motion by Mean Curvature ...

The $\infty$-Laplacian is the formal limit of the $p$-Laplacian as $p \ri \infty$. The $p$-Laplacian is
\[
\De_p u \, :=\, D_i\big( |Du|^{p-2}D_i u\big)\,=\, 0,
\]
and $\De_p u=0$ is the Euler-Lagrange PDE of the $p$-Dirichlet functional
\[
E_p(u,\Om)\, :=\, \int_{\Om}|Du|^p.
\]
The $\infty$-Laplacian is a quasilinear, degenerate (elliptic) PDE in non-divergence form, and will be our primary example in these notes.

\ms

\noi \emph{Hard-to-interpret ``singular solutions" of the  $\infty$-Laplacian:} 

(a) by rewriting \eqref{eq7} as
\beq \label{eq9}
\De_\infty u \, =\, D_iu \, D_i \Big(\frac{1}{2}|Du|^2 \Big)\,=\, 0,
\eeq
we see that every classical (i.e.\ $C^1$) solution of the Eikonal PDE
\[
|Du|^2 - 1\,=\, 0,
\]
formally solves \eqref{eq9}, but not  \eqref{eq7} (When we write \eqref{eq7} expanded in the form of \eqref{eq8}, we have the product of distribution with a $C^0$ function, which is not well defined) !!!

\ms

(b) In two dimensions, the $\infty$-Laplacian takes the form 
\[
\De_\infty u \,=\, (u_x)^2u_{xx}\, +\, 2u_x u_y u_{xy}\, +\, (u_y)^2u_{yy} \,=\, 0
\]
(in $x$-$y$ notation), when applied to a function $u=u(x,y)$, $u:\R^2 \ri \R$. By separating variables and looking for solutions of the form $u(x,y)=f(x)+g(y)$, a simple calculation leads to the saddle ``solution"
\[
u(x,y)\, =\, x^{4/3}- y^{4/3}
\]
which is merely $C^{1,1/3}(\R^2)$ H\"older continuous and lacks 2nd derivatives along the axes $x=0$ and $y=0$ !!!

\ms

In both $(a)$ and $(b)$, it might be that the singular (not twice differentiable) ``solution" is the unique $C^1$ solution of (an appropriately formulated) minimisation problem for the supremal functional \eqref{eq8}. The latter makes perfect sense if only 1st derivatives exist, at least a.e.\ on the domain, with \emph{no need} for reference to 2nd derivatives.

\ms

\noi (2) \textbf{More examples}:

(a) Linear 2nd order PDE (static, or 1st order in time) with $C^0$ coefficients:
\begin{align}
A_{ij}(x) D_{ij}^2u(x)\, +\, B_k(x)D_ku(x)\, +\, c(x)u(x)\, &=\, 0, \nonumber\\
-u_t\, +\, A_{ij}(x) D_{ij}^2u(x)\, +\, B_k(x)D_ku(x)\, +\, c(x)u(x)\, &=\, 0. \nonumber
\end{align}

(b) Quasilinear 2nd order divergence PDE (static, or 1st order in time) with $A$ in $C^1$ and $B$ in $C^0$:
\begin{align}
D_i\big(A_i(\cdot,u,Du)\big)\, +\,  B(\cdot,u,Du)\, &=\, 0, \nonumber\\
-u_t\, +\, D_i\big(A_i(\cdot,u,Du)\big)\, +\, B(\cdot,u,Du)\, &=\, 0. \nonumber
\end{align}

(c) Quasilinear 2nd order non-divergence PDE (static, or 1st order in time) with $A$, $B$ in $C^0$:
\begin{align}
A_{ij}(\cdot,u,Du)D_{ij}^2u\, +\,  B(\cdot,u,Du)\, &=\, 0, \nonumber\\
-u_t\, +\, A_{ij}(\cdot,u,Du)D_{ij}^2u\, +\, B(\cdot,u,Du)\, &=\, 0. \nonumber
\end{align}
This case includes $(a)$ and $(b)$ as special instances.

(d) Hamilton-Jacobi-Bellman type and Isaacs type equations (static, or 1st order in time):
\begin{align}
\sup_{a \in A} \big\{A^a_{ij}(\cdot)D_{ij}^2u \big\}\, +\,  B(\cdot,u,Du)\, &=\, 0, \nonumber\\
-u_t\, +\, \sup_{a \in A} \big\{A^a_{ij}(\cdot)D_{ij}^2u \big\}\, +\, B(\cdot,u,Du)\, &=\, 0, \nonumber
\end{align}
and
\begin{align}
\inf_{b\in B}\sup_{a \in A} \big\{A^{ab}_{ij}(\cdot)D_{ij}^2u\big\} \, +\,  B(\cdot,u,Du)\, &=\, 0, \nonumber\\
-u_t\, +\, \inf_{b\in B}\sup_{a \in A} \big\{A^{ab}_{ij}(\cdot)D_{ij}^2u \big\}\, +\, B(\cdot,u,Du)\, &=\, 0, \nonumber
\end{align}
where $\{A^{a}_{ij}\}_{a\in A}$ and $\{A^{ab}_{ij}\}_{a\in A}^{b\in B}$ are families of linear coefficients. Both classes of PDE are fully nonlinear.

(e) Functions of the eigenvalues of the Hessian (static, or 1st order in time): the Hessian $D^2u$ is (pointwise) a real symmetric $n \by n$ matrix. By the Spectral Theorem, the spectrum $\si(D^2u)$ consists of $n$ real eigenvalues, which we place in increasing order:
\[
\si(D^2u)\, =\, \big\{\la_1(D^2u),...,\la_n(D^2u) \big\}, \ \ \la_i\leq \la_{i+1}.
\]
Then, consider the PDE of the form
\begin{align}
G\big(\cdot,u,Du;\la_1(D^2u),..., \la_n(D^2u) \big)\, &=\, 0, \nonumber\\
-u_t\, +\, G\big(\cdot,u,Du;\la_1(D^2u),..., \la_n(D^2u)\big)\, &=\, 0. \nonumber
\end{align}
This class consists of fully nonlinear PDE.

Interesting special cases: 

(i) take $G(l_1,...,l_n,p,r,x):=\sum_{i=1}^{n}l_i-f(r)$. Then we obtain the nonlinear \emph{Poisson} equation:
\[
\De u \ =\, f(u).
\]

(ii) take $G(l_1,...,l_n,p,r,x):=\Pi_{i=1}^{n}l_i-f(x)$. Then we obtain the \emph{Monge-Amper\'e} equation
\[
\det(D^2u)\, =\,f. 
\]

(f) Obstacle problems and Gradient Constraint Problems: this includes fully nonlinear PDE of the general form
\begin{align}
&\max\big\{ F(D^2 u), u-f\big\}\,=\, 0,  \nonumber\\
&\min\big\{ F(D^2 u), |Du|-f\big\}\,=\, 0.  \nonumber
\end{align}

\ms

\noi Below we give the appropriate notion of ellipticity, which is necessary for the development of the theory.

\ms

\noi \textbf{Definition 6 (Degenerate Ellipticity)\footnote{In other texts, ellipticity is defined with the \textbf{opposite inequality}. This choice of convention that we make is more convenient, since we consider non-divergence operators like $\De_\infty$. We do not use integration by parts, hence there is no reason to consider ``minus the operator". According to our convention, $\De$ and $\De_\infty$ are elliptic, while for text using the opposite convention,  $-\De$ and $-\De_\infty$ are elliptic.}} \emph{The PDE
\[
F\big(\cdot,u,Du,D^2u\big)\, =\, 0
\]
is called degenerate elliptic, when the coefficient
\[
F\ :\ \Om \by \R \by \R^n \by \mS(n) \larrow \R 
\]
satisfies the weak monotonicity
\beq \label{eq10}
X\leq Y \text{ in }\mS(n)\ \ \Longrightarrow \ \ F(x,r,p,X)\leq F(x,r,p,Y),
\eeq
for all $(x,r,p) \in \Om \by \R \by \R^n$.}
\ms

\noi \textbf{Remark 7.} 

\noi \textit{$\al$) We recall that, by definition, the inequality $A\geq 0$ in the matrix space $\mS(n) \sub \R^{n\by n}$ means that $\si(A) \sub [0,\infty)$, that is all the eigenvalues of the symmetric matrix $A$ are non-negative:
\[
A:q \ot q\, =\, A_{ij}\, q_i\, q_j\, \geq \, 0, \ \ q \in \R^n.
\]
Of course, $A\leq B$ means $B-A\geq 0$. This definition is weak enough to include the case of Hamilton-Jacobi PDE $F(\cdot,u,Du)=0$ which has \emph{no dependence} on 2nd derivatives. Obviously, ``degenerate parabolic PDE" 
\[
-u_t\, +\, F\big(\cdot,u,Du,D^2u\big)\, =\, 0
\]
are defined to be degenerate elliptic in the previous sense, since derivative in ``time" is of first order. }

\textit{In these notes we will focus on the elliptic case. All the definitions, the techniques and the arguments in the proofs extend to the ``parabolic" case with minor modifications.}

\ms

\noi $\be$) \emph{Degenerate ellipticity in the examples:}

(a) \textit{the $\infty$-Laplacian is degenerate elliptic. Indeed, by defining
\[
F_\infty(p,X)\, :=\, X:p \ot p,
\]
we have $F_\infty(p,X)-F_\infty(p,Y)=(X-Y):p \ot p \leq 0$, when $X-Y\leq 0$. }

\textit{In the examples (b)-(d), the PDE is degenerate elliptic when the coefficient matrix of $D^2u$ is non-negative. }

\textit{For (e), the PDE is degenerate elliptic when for each $i=1,...,n$, the function $l_i \mapsto G(x,r,p;l_1,...,l_n)$ is non-decreasing. In particular, the Monge-Amper\'e PDE is degenerate elliptic when we restrict ourselves to the cone of \emph{convex functions $u$} and we also assume $f\geq0$. }

\textit{In (f), the PDE is degenerate elliptic when $F$ is.}

\ms

Now we have enough material in order to proceed to the motivation and the definition of the central notion of ``weak solutions" we will use in these notes. 

\ms

\noi \textbf{The pedagogical idea behind Viscosity Solutions} is to:

{\center{

\fbox{\parbox[pos]{320pt}{
\ms

Use the Maximum Principle in order to ``pass derivatives to smooth test functions" in a nonlinear fashion, without duality!

\ms
}}

\ms\ms }}

\noi \textbf{Motivation of the Definition.} Suppose that $u\in C^2(\Om)$ is a classical solution of the PDE
\beq \label{eq11}
F\big(x,u(x),Du(x),D^2u(x)\big)\, =\, 0,\ \ x\in \Om, 
\eeq
and that \eqref{eq11} is degenerate elliptic, that is, $F$ satisfies \eqref{eq10}. 

Assume further that at \emph{some} $x_0 \in \Om$, $u$ can be ``touched from above" by \emph{some} smooth function $\psi \in C^2(\R^n)$ at $x_0$. By the latter, we mean that the difference has a vanishing (local) maximum in a neighbourhood of $x_0$:
\beq  \label{eq12}
u-\psi \, \leq 0 =(u-\psi)(x_0), \ \ \text{ on a ball }\mB_r(x_0) \sub \Om.
\eeq
\[
\includegraphics[scale=0.22]{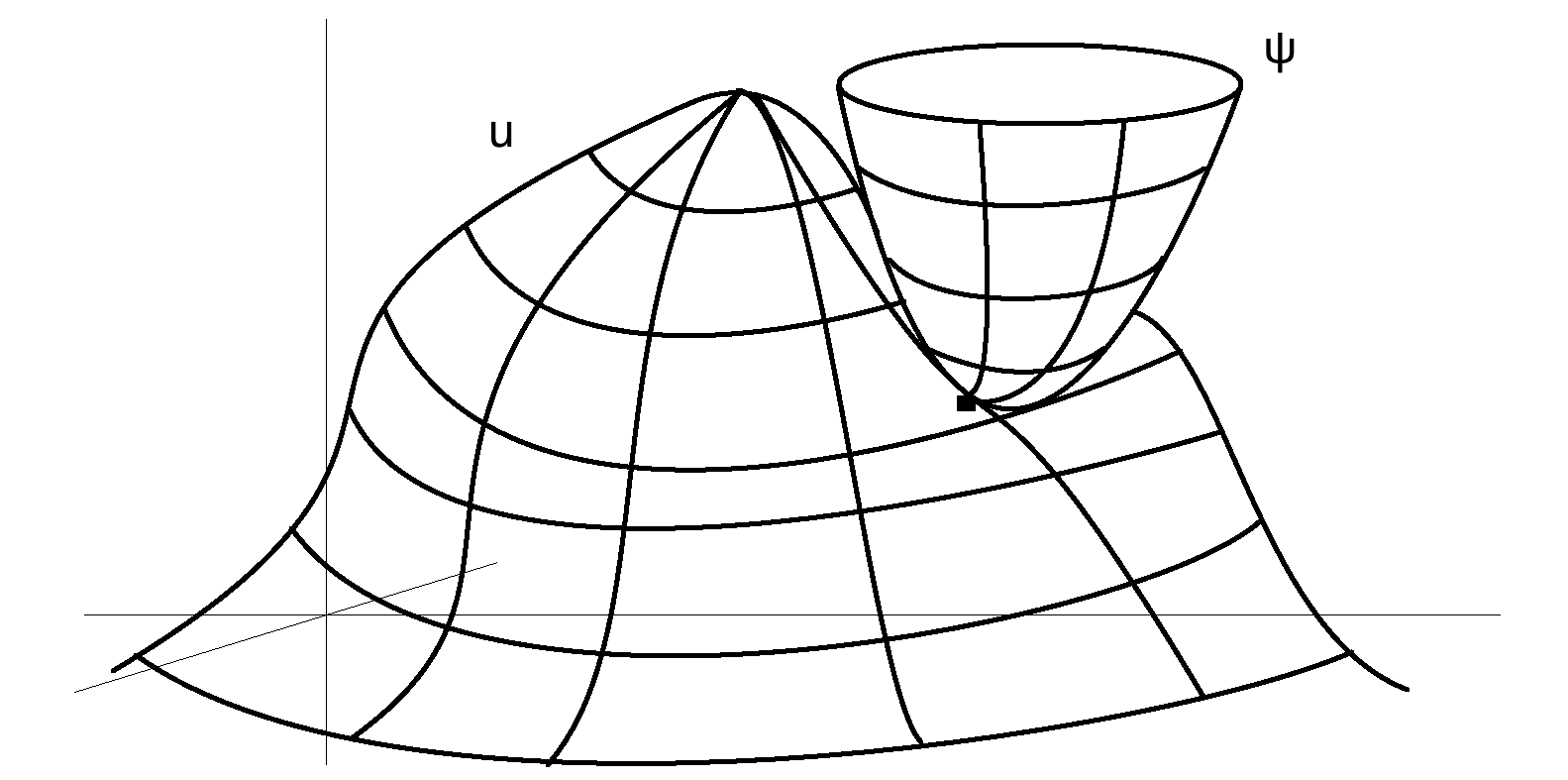}
\]
\[
\includegraphics[scale=0.18]{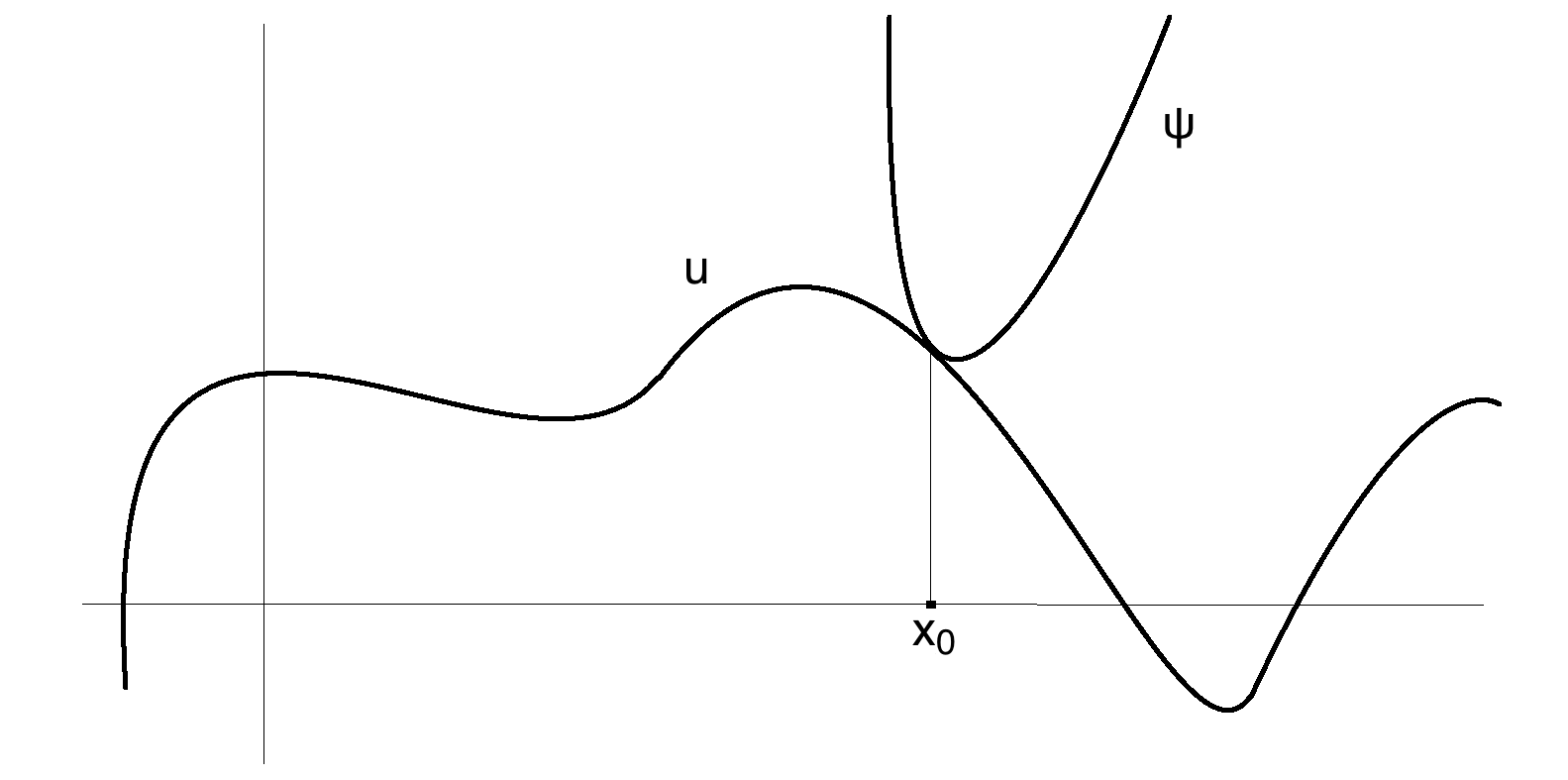}
\]
By \eqref{eq12} and Calculus, at the maximum point $x_0$ we have vanishing gradient and non-positive hessian:
\begin{align}
D(u-\psi)(x_0) &= 0, \ \text{ in } \R^n  \label{eq13}\\
D^2(u-\psi)(x_0) &\leq 0,  \ \text{ in } \mS(n).  \label{eq14}
\end{align}
By using that $u$ is a solution, we have
\begin{align}
 \hspace{40pt} 0 \ &= \ F\big(x_0,u(x_0),Du(x_0),D^2u(x_0)\big) \nonumber\\
&=\ F\big(x_0,\psi (x_0),D\psi(x_0),D^2u(x_0)\big) \hspace{30pt}\big(\text{by }\ \eqref{eq12}, \eqref{eq13}\big) \nonumber\\
& \leq\  F\big(x_0,\psi (x_0),D\psi(x_0),D^2\psi(x_0)\big)  \hspace{30pt}\big(\text{by }\ \eqref{eq11}, \eqref{eq14}\big).\nonumber
\end{align}
Hence, we conclude that if $u$ is a solution of \eqref{eq11}, then
\beq  \label{eq15}
\left.
\begin{array}{c}
u-\psi \leq 0 =(u-\psi)(x_0)\\
 \text{on a ball }\mB_r(x_0) \sub \Om,\\
x_0 \in \Om,\ \psi \in C^2(\R^n)
\end{array}
\right\}
 \ \ \Rightarrow \ \ F\big(x_0,\psi (x_0),D\psi(x_0),D^2\psi(x_0)\big) \geq 0.
\eeq

\noi Similarly, we conclude that if $u$ is a solution of \eqref{eq11}, then
\beq  \label{eq16}
\left.
\begin{array}{c}
u-\phi \leq 0 =(u-\phi)(y_0)\\
 \text{on a ball }\mB_r(y_0) \sub \Om,\\
y_0 \in \Om,\ \phi \in C^2(\R^n)
\end{array}
\right\}
 \ \ \Rightarrow \ \ F\big(y_0,\phi (y_0),D\phi(y_0),D^2\phi(y_0)\big) \leq 0.
\eeq

\noi The crucial observation now is that

{\center{

\fbox{\parbox[pos]{320pt}{
\ms

The implications \eqref{eq15} and \eqref{eq16} have been derived by using that $u \in C^2(\Om)$ is a classical solution of the PDE \eqref{eq11}, but DO NOT DEPEND ON THE EXISTENCE OF $(Du,D^2u)$ !!!

\ms
}}

\ms\ms }}

\noi In particular,  \eqref{eq15} and \eqref{eq16} make sense when $u$ is merely continuous. Hence, this leads to the idea that we can extract $(Du,D^2u)$ from the PDE, by testing against smooth ``touching" functions, and using the Maximum Principle together with the ellipticity of the PDE. Consequently, we have the next definition:

\ms

\noi \textbf{Definition 8 (Viscosity Solutions).\footnote{In other texts Viscosity sub/super solutions are defined with the \textbf{opposite inequalities}. The direction of the inequalities \emph{corresponds to the choice of convention in the ellipticity notion}. However, \textbf{no confusion should arise for the readers because the definition is essentially the same in both cases}: when it comes down to writing the inequalities for, say, $\De_\infty$, in either choice of conventions we have
\[
D\psi(x)\ot D\psi(x):D^2\psi(x)\geq 0 
\]
when $\psi$ touches $u$ from above at $x \in \Om$ and
\[
D\phi(y)\ot D\phi(y):D^2\phi(y)\leq 0 
\]
when $\phi$ touches $u$ from below at $y \in \Om$.}} \emph{Let $u\in C^0(\Om)$, $\Om\sub \R^n$, and consider the degenerate elliptic PDE (see \eqref{eq10})
\beq
F\big(\cdot,u,Du,D^2u\big)\, =\, 0.
\eeq
(a) We say that $u$ is a Viscosity Subsolution of the PDE (or a Viscosity Solution of $F\big(\cdot,u,Du,D^2u\big)\geq 0$) when
\beq  \label{eq17}
\left.
\begin{array}{c}
u-\psi \leq 0 =(u-\psi)(x_0)\\
 \text{on a ball }\mB_r(x_0) \sub \Om,\\
x_0 \in \Om,\ \psi \in C^2(\R^n)
\end{array}
\right\}
 \ \ \Rightarrow \ \ F\big(x_0,\psi (x_0),D\psi(x_0),D^2\psi(x_0)\big) \geq 0.
\eeq
(b) We say that $u$ is a Viscosity Supersolution of the PDE (or a Viscosity Solution of $F\big(\cdot,u,Du,D^2u\big)\leq 0$) when
\beq  \label{eq18}
\left.
\begin{array}{c}
u-\phi \leq 0 =(u-\phi)(y_0)\\
 \text{on a ball }\mB_r(y_0) \sub \Om,\\
y_0 \in \Om,\ \phi \in C^2(\R^n)
\end{array}
\right\}
 \ \ \Rightarrow \ \ F\big(y_0,\phi (y_0),D\phi(y_0),D^2\phi(y_0)\big) \leq 0.
\eeq
(c) We say that $u$ is a Viscosity Solution, when it is both a Viscosity Subsolution and a Viscosity Supersolution.
}
\ms

\noi In words, the definition means that 

{\center{

\fbox{\parbox[pos]{320pt}{
\ms

\emph{when we touch from above (respectively, below) by a smooth test function at a point, then, the test function is a subsolution of} 

\ \ \, \, \emph{the equation (respectively, supersolution) at that point.}

\ms
}}

\ms\ms }}

\noi \textbf{Remark 9.} \ms

\textit{(a) It might look strange to the reader that the definition is actually a \emph{criterion}: \underline{whenever} we can touch from either above or below, then we have an appropriate inequality. At points where no upper or lower touching functions exist, the candidate for solution automatically solves the PDE at these points.}\ms

\textit{(b) It is not a priori clear that the definition is \emph{not void}, that is that we can indeed find such touching functions at ``sufficiently many" points on $\Om$. If this were not the case for a continuous candidate solution $u$, this $u$ would \emph{automatically} be a solution since we would have nothing to check.}
\ms

\textit{(c) Splitting the notion of solution to 1-sided halves is essential. If we had assumed that we can \emph{simultaneously} touch from above and below at any point, then as we will see in the next Chapter, we would not really define a new truly ``weak" solution, since for this hypothetical notion 1st derivatives would have to be Lipschitz continuous at the touching point).}\ms

\textit{(d) It is not clear yet that the notion is compatible with classical solutions. That is, we do not yet know whether the family of $C^2$ Viscosity Solutions and the family of classical solutions of a degenerate elliptic PDE coincide.}

\ms

\noi All the issues raised in the remark, will be answered in the next Chapter, where we also investigate another equivalent definition and exploit its analytic properties. This alternative point of view is based on the possibility to define ``pointwise weak derivatives $(Du,D^2u)$", without reference to ``test functions".

Note that{

\center{\fbox{\parbox[pos]{220pt}{
\ms

merely $C^0$ functions can be interpreted as solutions of fully nonlinear 2nd order PDE!

\ms
}}}

}\ms \ms

\noi The fact that we can drop all derivatives from the PDE, has the following outstanding implication that we will see in following chapters:

{

\center{\fbox{\parbox[pos]{280pt}{
\ms

VISCOSITY SOLUTIONS PASS TO LIMITS UNDER MERELY $C^0$ CONVERGENCE, NO CONTROL ON 

\ \ \ \ \ \ \ \ \ ANY DERIVATIVES IS REQUIRED !!!

\ms
}}}

}\ms 
\noi In the above, ``$C^0$ convergence" stands for locally uniform convergence, namely, the convergence arising from the natural topology of the space $C^0(\Om)$. Also,

{
\center{\fbox{\parbox[pos]{250pt}{
\ms

The Dirichlet problem
\[
\left\{
\begin{array}{l}
F\big(\cdot,u,Du,D^2u\big)\,=\,0, \text{ in }\Om,\\
\hspace{77pt}u\, =\, b, \text{ on }\p \Om,
\end{array}
\right.
\]
has a unique continuous Viscosity Solution, under natural assumptions on $F,\, \Om,\, b$!
\ms
}}}

}\ms  \ms

\noi The latter solution might be nowhere differentiable on $\Om$, and hence $(Du,D^2u)$ might have no classical sense! Such strong convergence and existence results hold only for the linear theory of Distributions (Generalised Functions) for PDE with constant coefficients.

\ms

\ms

\ms

\noi \textbf{Remarks on Chapter 1.} Viscosity solutions were introduced first by Crandall and Lions in \cite{CL} as a uniqueness criterion for 1st order PDE. The essential idea regarding passage to limits was observed earlier by Evans in \cite{E1, E2}. As usually, the original point of view of the notions was slightly different. Crandall-Evans-Lions \cite{CEL} have definitions similar to the ones we use today. The viscosity notions were later extended to the 2nd order case by Ishii in \cite{I1}. Aronsson was the first to consider $L^\infty$ variational problems in \cite{A1}. The $\infty$-Laplacian made its first public appearance in his subsequent papers \cite{A3, A4}. Aronsson himself was also the first to observe  in \cite{A6, A7} ``singular solutions" which at the time were hard-to-justify (visosity theory had not been develop yet). Battacharya-DiBenedetto-Manfredi were the first to consider the $\infty$-Laplacian in the Viscosity sense in \cite{BDM}. Background material on divergence structure PDE (although remotely needed for this chapter) can be found in standard references like the textbook of Evans \cite{E4}. For PDE theory of classical solutions to nonlinear elliptic equations, the standard reference is Gilbarg-Trudinger \cite{GT}. For regularity estimates of viscosity solutions to fully nonlinear equations, the interested reader may consult (the advanced book of) Cabr\'e-Caffarelli \cite{CC}.

\chapter[Second Definitions, Analytic Properties]{Second Definitions and Basic Analytic Properties of the notions}

In the previous chapter we defined a notion of non-differentiable ``weak" solutions which applies to degenerate elliptic PDE of the general form
\beq \label{2.1}
F\big(\cdot,u,Du,D^2u\big)\, =\, 0
\eeq
that makes sense when $u$ is merely in $C^0(\Om)$, $\Om\sub \R^n$. As always, the ``coefficients" $F=F(x,r,p,X)$ is a continuous and possibly nonlinear function defined on $\Om \by \R \by \R^n \by \mS(n)$. We recall that the appropriate degenerate ellipticity condition we need to impose is
\beq \label{2.2}
X\leq Y \text{ in }\mS(n)\ \ \Longrightarrow \ \  F(x,r,p,X)\leq F(x,r,p,Y) .
\eeq
In this chapter we will handle the issues raised at the end of the previous chapter, the most important of which were

\ms

(a) \emph{Is the definition of Viscosity Solutions non-void? Namely, is it true that any continuous function can be touched from above and below by a smooth function at ``sufficiently many" points of the domain, in order to uniquely determine a solution on open subsets?}

\ms

(b) \emph{Are Viscosity Solutions compatible with classical solutions? Namely, a twice differentiable Viscosity Solutions of \eqref{2.1} is the same as a classical solution?}

\ms

(c) \emph{Why is it necessary to split the notion to two halves?}

\ms

Were the answer ``no" to (a), then the definition would be trivial since, if such a hypothetical continuous function existed, would be a solution to ALL equations in the Viscosity sense!

\ms

Were the answer ``no" to (b), the definition would obviously be naive. If one was to prove existence of a a ``weak" solution which happens to be twice differentiable, then this solution should definitely solve the PDE pointwise.

\ms

Concerning (c), as we shall see a little later, if we were to require that we can touch from above and below by smooth functions \emph{simultaneously at the same} points, then we would not go in a class much weaker than that of strong solutions.
\ms

We begin by introducing an alternative but equivalent definition.

\ms 
\noi \textbf{Motivation for the Definition.} The idea is to somehow

{

\center{\fbox{\parbox[pos]{310pt}{
\ms

DEFINE POINTWISE GENERALISED DERIVATIVES $(Du,D^2u)$ FOR NON-SMOOTH FUNCTIONS BY MEANS

 \hspace{60pt} OF THE MAXIMUM PRINCIPLE !

\ms
}}}

}
\ms

\noi We begin by slightly relaxing the classical Calculus definition of derivatives. The ideas is simple -  we just take as definition the Taylor expansion at a point. This makes no difference for 1st derivatives, but is essential for 2nd (and higher order) derivatives.

\ms

\noi \textbf{Definition 1 (Pointwise Derivatives as Taylor Expansions).} \emph{Let $u\in C^0(\Om)$ and $x\in \Om\sub \R^n$. We say that $u$ is twice differentiable at $x$ if there exist $p\in \R^n$ and $X\in \mS(n)$ such that
\beq \label{2.3}
u(z+x)\, =\, u(x)\, +\, p\cdot z\, +\, \frac{1}{2}X:z\ot z \, +\, o(|z|^2),
\eeq
as $z\ri 0$. In such case, we call $p$ the gradient of $u$ at $x$ and denote it by $Du(x)$ and we call $X$ the hessian of $u$ at $x$ and denote it by $D^2u(x)$. }

\emph{In the standard way, we shall say that $u$ is once differentiable at $x$ when
\beq \label{2.4}
u(z+x)\, =\, u(x)\, +\, p\cdot z\, +\, o(|z|),
\eeq
as $z\ri 0$, and then $p$ will be denoted by $Du(x)$.}

\ms

\noi \textbf{Remark 2.} \textit{We use ``$\cdot$" to denote the Euclidean inner product in $\R^n$. The meaning of \eqref{2.3} is that there exists a continuous increasing function $\si \in C^0[0,\infty)$ with $\si(0)=0$ such that
\beq \label{2.5}
\left|u(y)\, -\, u(x)\, -\, p\cdot (y-x)\, -\, \frac{1}{2}X:(y-x)\ot (y-x) \right| \, \leq \, \si(|y-x|)|y-x|^2,
\eeq
for $y \in \Om$. If such a pair $(Du(x),D^2u(x))$ exists, then it is unique. \eqref{2.4} is nothing but the classical notion of 1st derivative. }

\textsl{The new ingredient of \eqref{2.3} is that we defined $D^2u(x)$ at $x$, although the gradient $Du$ may not exist as a function near $x$. In classical calculus, $D^2u(x)$ is defined by differentiating the gradient map $Du:\Om\sub \R^n \larrow \R^n$ at $x$.
}

\ms

\noi \textbf{Example 3.} Let $W\in C^0(\R)$ be any continuous nowhere differentiable function, for instance the classical Weierstrass function. Then, the function
\[
w\, \in C^0(\R)\ , \ \ w(x)\, :=\,  W(x)|x|^{3},
\]
has two vanishing derivatives in the sense of \eqref{2.3} (i.e.\ $w'(0)=w''(0)=0$), since 
\[
|w(x)| \, \leq\, \Big(\max_{[-1,+1]}|W|\Big) |x|^3\, =\, o(|x^2|), 
\]
as $x\ri 0$. However, $w''$ can not be defined in the classical way since $w'$ does not exist pointwise anywhere near $x=0$ ($w'$ is a genuine distribution in $\mathcal{D}'(\R)$).

\ms

\noi \textbf{The idea.} \eqref{2.3} defines a pointwise notion of derivative, but the problem is that if merely $u\in C^0(\Om)$, then $(Du(x),D^2u(x))$ may not exist anywhere in $\Om$ !!! We rectify this problem by ``splitting the notion to two one-sided halves":

\ms

\noi \textbf{Definition 4 (Generalised Pointwise Derivatives).} \emph{Let $u\in C^0(\Om)$, $\Om\sub \R^n$ and $x\in \Om$. we define the sets
\begin{align} \label{2.6a}
\J^{2,+}u(x)\, :=\, \Big\{ (p,X) \in \R^n \by \mS(n)& \ \Big|  \ \text{ as }z\ri0, \text{ we have}: \nonumber\\
u(z+x)\, \leq\, u(x) \, &+\, p\cdot z\, +\, \frac{1}{2}X:z\ot z \, +\, o(|z|^2) \Big\},
\end{align}
and
\begin{align} \label{2.7a}
\J^{2,-}u(x)\, :=\, \Big\{ (p,X) \in \R^n \by \mS(n)& \ \Big|  \ \text{ as }z\ri0, \text{ we have}: \nonumber\\
u(z+x)\, \geq\, u(x) \, &+\, p\cdot z\, +\, \frac{1}{2}X:z\ot z \, +\, o(|z|^2) \Big\}.
\end{align}
We call $\J^{2,+}u(x)$ the 2nd order Super-Jet of $u$ at $x$ and $\J^{2,-}u(x)$ the 2nd order Sub-Jet of $u$ at $x$.}

\ms
\noi \textbf{Remark 5.} \textit{(a) A pair $(p,X) \in \J^{2,+}u(x)$ at some $x\in \Om$ plays the role of an 1-sided (non-unique, 1st and 2nd) superderivative(s) of $u$ at $x$. Similarly, a pair $(p,X) \in \J^{2,-}u(x)$ at some $x\in \Om$ plays the role of an 1-sided (non-unique, 1st and 2nd) subderivative(s) of $u$ at $x$.}

\textit{(b) It may may well happen that $\J^{2,+}u(x) \cap \J^{2,-}u(x)=\emptyset$ for all $x\in \Om$, but as we will see next there exist ``many" points $x$ for which  $\J^{2,+}u(x) \neq \emptyset$ and $\J^{2,-}u(x) \neq \emptyset$.}

\[
\underset{\underset{\text{\footnotesize and ``X" the quadratic 1-sided approximations. \hspace{150pt}}}{\text{\footnotesize Illustration of sub- and super- Jets. ``p" denotes the linear 1-sided approximations}}}{\includegraphics[scale=0.22]{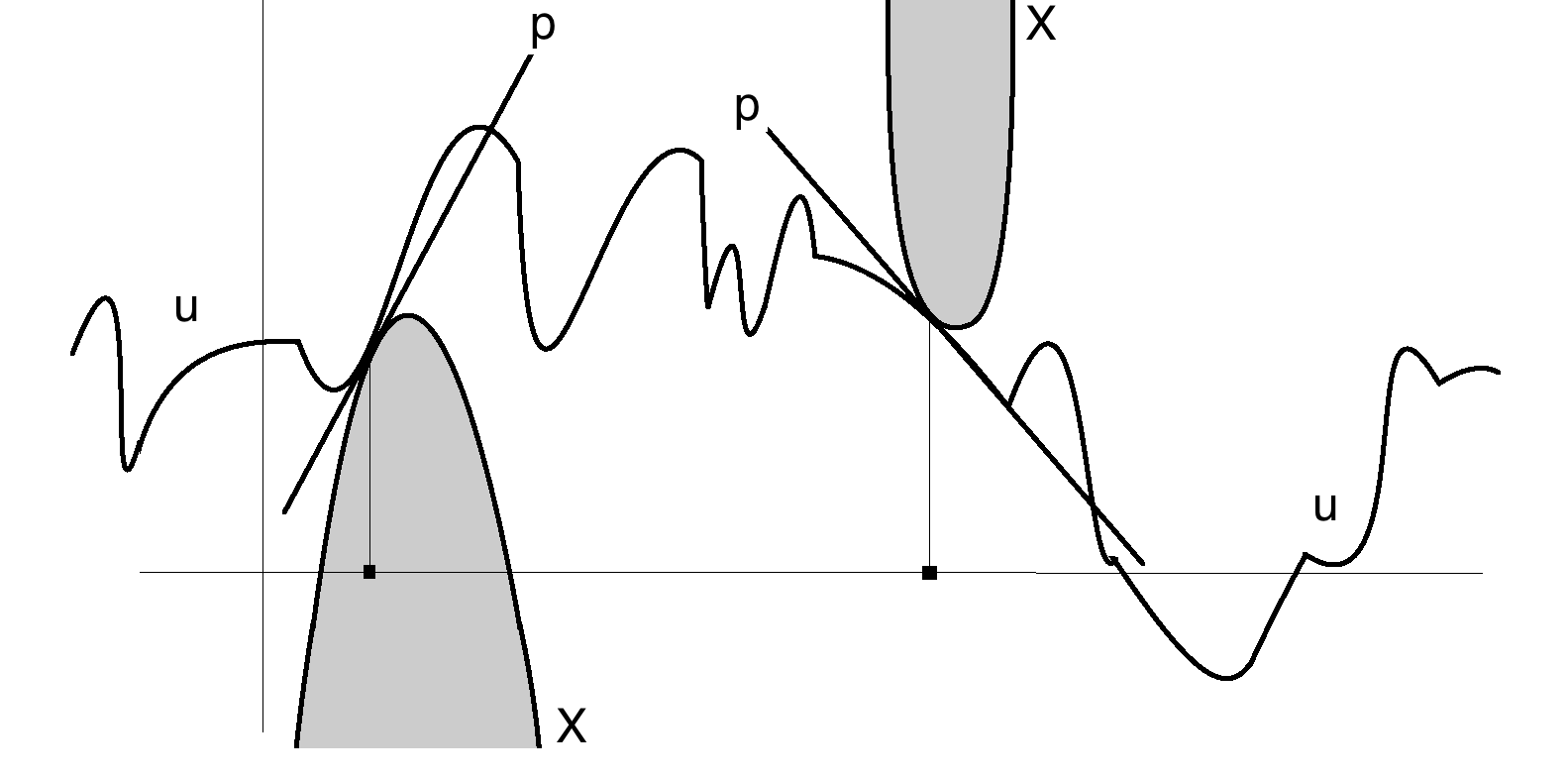}}
\]
Before proceeding to study the properties of Jets, we re-define Viscosity Solutions. The idea is to
{
\center{\fbox{\parbox[pos]{250pt}{
\ms

split the notion of derivatives $(Du,D^2u)$ to pointwise generalised 1-sided (upper/lower) 

\ms
\centerline{and} 
\ms

split the notion of solution of $F=0$ to 1-sided (sub-/super-) for the respective derivatives.

\ms
}}}

}
\ms

\noi \textbf{Definition 6 (Viscosity Solutions)\footnote{See the footnote of Definition 8 of Chapter 1, regarding the alternative popular convention in the definition of the notions regarding the directions of the inequalities.}.} \emph{Consider the degenerate elliptic PDE
\[
F\big(\cdot,u,Du,D^2u\big)\, =\, 0,
\]
where $F\in C^0\big(\Om\by\R \by\R^n \by\mS(n) \big)$ (see \eqref{eq10}). }

(a) \emph{ The function $u\in C^0(\Om)$ is a Viscosity Subsolution of the PDE (or Viscosity Solution of $F\big(\cdot,u,Du,D^2u\big)\geq 0$) on $\Om$, when 
\beq \label{2.8a}
(p,X)\in \J^{2,+}u(x)\ \ \Longrightarrow \ \ F(x,u(x),p,X)\geq 0,
\eeq
for all $x\in \Om$.}

(b) \emph{ The function $u\in C^0(\Om)$ is a Viscosity Supersolution of the PDE (or Viscosity Solution of $F\big(\cdot,u,Du,D^2u\big)\leq 0$) on $\Om$, when 
\beq  \label{2.9a}
(p,X)\in \J^{2,-}u(x)\ \ \Longrightarrow \ \ F(x,u(x),p,X)\leq 0,
\eeq
for all $x\in \Om$.}

(c) \emph{ The function $u\in C^0(\Om)$ is a Viscosity Solution of the PDE when it is both a Viscosity Subsolution and a Viscosity Supersolution.}

\ms

\noi \textbf{Remark 7.} \textit{(a) The Viscosity notions given by \eqref{2.8a} and \eqref{2.9a}, can be restated as 
\beq   \label{2.10}
\inf_{(p,X)\in \J^{2,+}u(x)} F(x,u(x),p,X) \geq 0
\eeq
for all $x\in \Om$, and
\beq   \label{2.11}
\sup_{(p,X)\in \J^{2,-}u(x)} F(x,u(x),p,X) \leq 0
\eeq
for all $x\in \Om$. Obviously, if $\J^{2,+}u(x)=\emptyset$ for some $x\in \Om$, then 
\[
\inf_{(p,X)\in \J^{2,+}u(x)}F(x,u(x),p,X) =+\infty
\] 
and $u$ trivially is a subsolution of the PDE at this point. A similar observation holds for $\J^{2,-}u(x)$ and supersolutions as well. }

\textit{(b) First order Sub-Jets and Super-Jets can be defined in the obvious way:
\[ 
\J^{1,+}u(x)\, :=\, \Big\{ p \in \R^n \ \Big|  \ u(z+x)\, \leq\, u(x) \, +\, p\cdot z\, +\, o(|z|) ,\ \text{ as }z\ri0 \Big\},
\]
\[ 
\J^{1,-}u(x)\, :=\, \Big\{ p \in \R^n \ \Big|  \ u(z+x)\, \geq\, u(x) \, +\, p\cdot z\, +\, o(|z|) ,\ \text{ as }z\ri0 \Big\}.
\]
Since in these notes we are mostly concerned with the 2nd order theory, first order Jets will be of minor importance.
}

\ms

The equivalence between the two seemingly different definitions of Viscosity Solutions will be established, once we first establish the following properties of Jets.

\ms

\noi \textbf{Lemma 7 (Properties of Jets).}  \emph{Let $u\in C^0(\Om)$ and $x\in \Om\sub \R^n$. Then,}

\noi  \emph{(a) $\J^{2,-}u(x)=-\J^{2,+}(-u)(x)$. }

\ms

\noi  \emph{(b) $\J^{2,+}u(x)$ is a convex subset of $\R^n \by \R^{n \by n}$. }
\ms

\noi \emph{(c) If $\J^{1,+}u(x)\neq \emptyset$ and $\J^{1,-}u(x)\neq \emptyset$, then both sets are singletons, $u$ is differentiable at $x$ and both sets coincide with the gradient:
\[
\J^{1,\pm}u(x)\neq \emptyset\ \ \Longrightarrow\ \ \J^{1,+}u(x)\, =\, \J^{1,-}u(x)\,=\, \{Du(x)\}.
\]
}
\noi \emph{(d) If $\J^{2,+}u(x)\neq \emptyset$ and $\J^{2,-}u(x)\neq \emptyset$, then $u$ is $C^{1,1}$ at $x$, that is
\[
\big|u(z+x)-u(x)-Du(x)\cdot z \big|\, \leq \, O(|z|^2),
\]
as $z\ri 0$. Moreover,
\[
(p,X^{\pm}) \in \J^{2,\pm}u(x)\ \ \Longrightarrow\ \ X^- \leq X^+.
\]
}
\noi  \emph{ (e) (compatibility with classical derivatives) If $u$ is twice differentiable at $x$ (see \eqref{2.3}), then
\beq
\J^{2,+}u(x)\, =\, \Big\{\big(Du(x),D^2u(x)+A\big)\ \big| \ A\geq 0 \text{ in }\mS(n) \Big\}
\eeq}

\noi  \emph{(f) If $\J^{2,+}u(x) \neq \emptyset$, then it has infinite diameter: $\diam\big(\J^{2,+}u(x)\big) =+\infty$.}

\ms

\noi  \emph{(g) (Jets non-empty on dense sets) The set
\[
\big\{ x\in \Om \ \big| \ \J^{2,+}u(x) \neq \emptyset \big\}
\]
is dense in $\Om$.}

\noi  \emph{(h) For each $p\in \R^n$, the ``slice" set
\[
\big\{ X\in  \mS(n)\ \big| \ (p,X)\in \J^{2,+}u(x) \big\}
\]
is closed. }

\ms

\noi \textbf{Proof.} (a) is a direct consequence of the definitions \eqref{2.6a}, \eqref{2.7a}.

\ms

(b) If $(p',X')$, $(p'',X'') \in \J^{2,+}u(x)$, then for any $\la \in [0,1]$, \eqref{2.6a} implies that
\begin{align}
\la u(z+x)\, &\leq \, \la u(x)\, +\, \la p' \cdot z\, +\, \frac{1}{2}\la X' \, +\, o(|z|^2), \nonumber\\
(1-\la) u(z+x)\, &\leq \, (1-\la)  u(x)\, +\, (1-\la)  p'' \cdot z\, +\, \frac{1}{2}(1-\la)  X'' \, +\, o(|z|^2),\nonumber
\end{align}
as $z\ri 0$. Hence, by summing we obtain
\[
\big(\la p' +(1-\la) p'', \la X' +(1-\la)X''\big) \in \J^{2,+}u(x).
\]

(c) Let $p^\pm \in \J^{1,\pm}u(x)$. Then, 
\begin{align}
u(z+x)\, -\, u(x)\,  &\leq\,  p^+\cdot z \, +\, o(|z|), \label{2.14a}\\
-u(z+x)\, +\, u(x)\,  &\leq\,  - p^- \cdot z \, +\, o(|z|), \label{2.15a}
\end{align}
as $z\ri 0$. We set $z:=\e w$, where $\e>0$ and $|w|=1$ and and the inequalities to obtain
\[
(p^- -p^+)\cdot w\, \leq\, \frac{o(\e)}{\e},
\]
as $\e \ri 0$. Hence, by passing to the limit, we get
\[
|p^- -p^+|\, =\, \max_{|w=1|}\{(p^- -p^+)\cdot w\} \, \leq \, 0, 
\]
which gives $p^+=p^-$. Let $p$ denote the common value of $p^\pm$. Then, by \eqref{2.14a}, \eqref{2.15a}  we obtain
\[
u(z+x)\, =\, u(x)\,  +\,  p\cdot z \, +\, o(|z|)
\]
as $z\ri 0$, which says that $p=Du(x)$.

\ms

(d) Let $(p^\pm , X^\pm) \in \J^{2,\pm}u(x)$. Then, since $p^\pm \in \J^{1,\pm}u(x)$, (c) above implies that $p^+=p^-=Du(x)$ and hence $u$ is differentiable at $x$.  Moreover, we have
\begin{align}
u(z+x)\, -\, u(x)\,  -\,  Du(x)\cdot z \, \leq\, \frac{1}{2}X^+:z\ot z\, +\, o(|z|^2), \label{2.14b}\\
-u(z+x)\, +\, u(x)\,  +\,  Du(x) \cdot z \, \leq \, -\frac{1}{2}X^-:z\ot z\, +\, o(|z|^2), \label{2.15b}
\end{align}
and by setting $z:=\e w$, $\e>0$, $|w|=1$ and adding the inequalities, we obtain
\[
(X^- -X^+ ):w \ot w \, \leq\, \frac{o(\e^2)}{\e^2},
\]
as $\e \ri 0$. Hence, $X^- \leq X^+$. Again by \eqref{2.14b} and \eqref{2.15b}, we have  
\[
\big| u(z+x)\, -\, u(x)\,  -\,  Du(x)\cdot z  \big|\, \leq \, \Big(\frac{1}{2}\max \big|\si (X^\pm)\big|+1\Big)|z|^2,
\]
as $z\ri 0$, says that $u$ is $C^{1,1}$ at $x$. Here ``$\si$" denotes the set of eigenvalues of a matrix.

\ms

(e) If $u$ is twice differentiable at $x$, then we have
\begin{align}
u(z+x)\, &=\, u(x) \, +\, Du(x)\cdot z\, +\, \frac{1}{2}D^2u(x):z\ot z\, +\, o(|z|^2) \nonumber\\
  &\leq\, u(x) \, +\, Du(x)\cdot z\, +\, \frac{1}{2}(D^2u(x)+A):z\ot z\, +\, o(|z|^2), \nonumber
\end{align}
as $z\ri 0$. Hence,
\[
\big(Du(x),D^2u(x)+A\big) \in \J^{2,+}u(x).
\]
Conversely, if $(p,X) \in \J^{2,+}u(x)$ and $(Du(x),D^2u(x))$ exist, we have
\begin{align}
u(z+x)\, -\, u(x) \, &=\, Du(x)\cdot z\, +\, \frac{1}{2}D^2u(x):z\ot z\, +\, o(|z|^2) \nonumber\\
u(z+x)\, -\, u(x) \,   &\leq\,  p\cdot z\, +\, \frac{1}{2}X:z\ot z\, +\, o(|z|^2), \nonumber
\end{align}
which by subtraction give
\beq \label{2.13}
(Du(x)-p)\cdot z \, \geq \, \frac{1}{2}\big(X-D^2 u(x)\big):z\ot z\, +\, o(|z|^2)
\eeq
as $z\ri 0$. By setting $z:=\e w$ for $|w|=1$ and $\e>0$, \eqref{2.13} implies
\[
(p-Du(x))\cdot w \, \leq \, O(\e),
\]
as $\e\ri 0$. By passing to the limit, we get 
\[
|p-Du(x)| \, =\,  \max_{|w|=1} \big\{ (p-Du(x))\cdot w \big\} \, \leq  \, 0,
\]
and hence $Du(x)=p$. Again by \eqref{2.13}, since the left hand side vanishes, for $z=\e w$ we have
\[
\big(X-D^2 u(x)\big):w\ot w\, \geq \, \frac{o(\e^2)}{\e^2}.
\]
By letting $\e \ri 0$, we see that all the eigenvalues of $X-D^2 u(x)$ are non-negative. Hence, $X-D^2 u(x)\geq 0$ in $\mS(n)$. Consequently, (c) follows by setting $A:=X-D^2 u(x)$.

\ms

(f) follows from (e) be choosing $A:=t I$, $t>0$, where $I$ the identity of $\R^{n\by n}$. Indeed, we have that the ray
\[
\big\{ \big(Du(x),D^2u(x)+tI \big)\ : \ t\geq 0\big\}
\]
is contained in $\J^{2,+}u(x)$.

\ms 
(g) We simplify things by assuming in addition that $u\in C^0(\overline{\Om})$ and $\Om \Subset \R^n$. The general case follows by covering $\Om$ with an increasing sequence of bound open sets. Fix $x_0 \in \Om$ and $\e>0$. For $a\in \R$, consider the paraboloid centred at $x_0$ with curvature $1/\e$:
\[
Q^a(z)\, :=\, a\ +\ \frac{|z-x|^2}{2\e}.
\]
For $a\geq \sup_{\Om}u$, we have $Q^a \geq u$ on $\Om$. We start decreasing the parameter $a$ and slide the paraboloid downwards, until the graph of $Q^a$ touches the graph of $u$ at some point $x_\e \in \overline{\Om}$ for the first time (that is, without violating $Q^a \geq u$). 
\[
\includegraphics[scale=0.22]{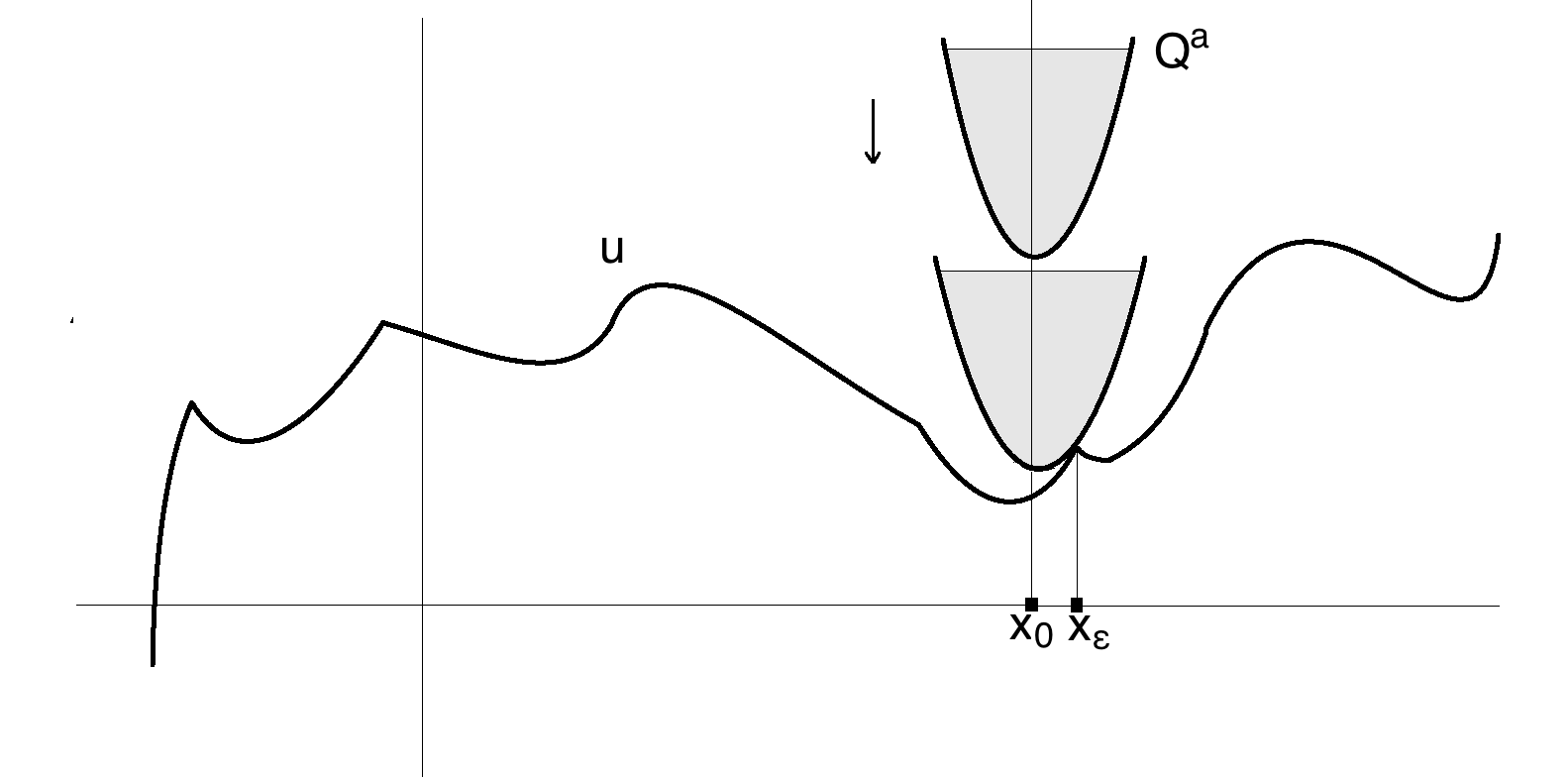}
\]
Since $Q^a \geq u$ on $\overline{\Om}$ and $Q^a(x_\e)=u(x_\e)$, we have constructed a maximum of $Q^a-u$ at $x_\e$:
\beq \label{2.14}
(Q^a-u)(z) \, \leq \, 0\, =\, (Q^a-u)(x_\e), \ \ \text{ for }z\in   \overline{\Om}.
\eeq
For $z:=x_0$ in \eqref{2.14}, we get
\[
a- u(x_0) \, \leq \, a -\frac{|x_0-x_\e|^2}{2\e} -u(x_\e)
\]
which gives
\[
|x_\e -x_0|\, \leq\, 2 \sqrt{\|u\|_{C^0(\Om)}\e}.
\]
Hence, for $\e >0$ small, $x_\e$ is an interior maximum in $\Om$, and also $x_\e \ri x_0$ as $\e \ri 0$. By Taylor expansion of the quadratic function $Q^a$ at $x_\e$, \eqref{2.14} implies
\begin{align}
u(z+x_\e)\, & \leq \, u(x_\e)\, +\, DQ^a (x_\e) \cdot z +\frac{1}{2}D^2  Q^a (x_\e) :z \ot z \nonumber\\
& \leq \, u(x_\e)\, +\, \left(\frac{x_\e - x_0}{\e}\right) \cdot z +\frac{1}{2}\frac{I}{\e}:z \ot z,\nonumber
\end{align}
as $z\ri 0$. Hence,
\[
\left( \frac{x_\e - x_0}{\e} , \frac{I}{\e} \right) \in \J^{2,+}u(x_\e).
\]
Since $x_0 \in \Om$ is arbitrary and $x_\e \ri x_0$ as $\e \ri 0$, it follows that the set of $x$'s for which  $\J^{2,+}u(x)\neq \emptyset$ is dense in $\Om$.

\ms

(h) Fix $p\in \R^n$ and suppose $\{(p,X_m)\}_1^\infty \sub \J^{2,+}u(x)$, where $X_m \ri X$ as
$m\rightarrow \infty$. Fix $\e >0$ and $m(\e)\in \N$ such that $|X_{m(\e)} -X|\leq \e$. For $\de>0$ fixed, we have
 \begin{align} \label{eq4.42}
\max_{|z|\leq \de}&
\left(\frac{u(z+x)-u(x)-p\cdot z-\frac{1}{2}X
:z\ot z}{|z|^2}\right) \nonumber\\
& \leq \, \e \ + \ \max_{|z|\leq \de}
\left(\frac{u(z+x)-u(x)-p\cdot z - \frac{1}{2}X_{m(\e)}:
z\ot z}{|z|^2}\right)\\
 &  \leq \, \e \ + \ o(1), \nonumber
 \end{align}
as $\de \rightarrow 0$. By passing to the limit $\de \rightarrow
0$ in \eqref{eq4.42} and then letting $\e \rightarrow 0$, we
obtain 
\[
\underset{z\ri 0}{\lim\sup}\left(\frac{u(z+x)-u(x)-p\cdot z-\frac{1}{2}X
:z\ot z}{|z|^2}\right) \, \leq \, 0,
\]
which says that $(p,X) \in \J^{2,+}u(x)$.  Consequently, (f) follows.      \qed

\ms

Having at hand these properties of Jets, we many now establish the equivalence between the two different definitions  of Viscosity Solutions, the one involving Jets and and the one involving touching by smooth test functions.

\ms

\noi \textbf{Theorem 8 (Equivalence of Definitions of Viscosity Solutions).} \emph{Let $u \in C^0(\Om)$ and $x\in \Om \sub \R^n$. Then, we have
\beq
\J^{2,+}u(x)\, =\, \Big\{ (D\psi(x),D^2\psi(x))\ \Big| \ \psi \in C^2(\R^n),\ u-\psi \leq 0 = (u-\psi)(x) \text{ near }x \Big\}.
\eeq
In particular, the two definitions of Viscosity (sub-, super-) Solutions (see \eqref{eq17}, \eqref{eq18} and \eqref{2.8}, \eqref{2.9}) coincide.}
\ms

\noi \textbf{Proof.} If $\psi \in C^2(\R^n)$ and $u-\psi \leq 0 = (u-\psi)(x)$ near $x$, then by Taylor expansion of $\psi$ at $x$ we have
\begin{align}
u(z+x)\, &\leq \, \psi(z+x) \, +\, u(x)\, -\, \psi(x) \nonumber\\
      & = \, u(x)\, +\, D\psi(x)\cdot z\, +\, \frac{1}{2}D^2\psi(x):z\ot z \, +\, o(|z|^2), \nonumber
\end{align}
as $z\ri 0$. Hence,
\[
\big(D\psi(x),D^2\psi(x)\big)\in \J^{2,+}u(x).
\]
Conversely, let $(p,X)\in \J^{2,+}u(x)$. This means that there is a $\si \in C^0[0,\infty)$ increasing with $\si(0)=0$ such that 
\[
u(z+x)\, \leq \, u(x)\, +\, p\cdot z\, +\, \frac{1}{2}X:z\ot z \, +\, \si(|z|)|z|^2, 
\]
as $z\ri 0$. The right hand side defines a function $\phi \in C^0(\R^n)$ such that $\phi(x)=u(x)$, $D\phi(x)=p$ and $D^2\phi(x)=X$, while 
\[
u-\phi \leq 0=(u-\phi)(x) ,\ \ \text{near }x.
\]
The problem is that \textit{$\phi$ is merely continuous and not in $C^2(\R^n)$}. The next lemma completes the proof of the Theorem, by showing that we can find $\tau \geq \si$ which in addition is in $C^2(0,\infty)$ and has the same properties as $\si$. We conclude by setting
\[
\psi(y)\, :=\, u(x)\, +\, p\cdot (y-x)\, +\, \frac{1}{2}X:(y-x)\ot (y-x) \, +\, \tau(|y-x|)|y-x|^2
\]
which is in $C^2(\R^n)$ and has the desired properties.     \qed

\ms

\noi \textbf{Lemma 9 (Mollifying the remainder).} \emph{Let $\si$ be an increasing function in $C^0(0,\infty)$ with $\si(0^+)=0\footnote{We use the popular notation $\si(0^+):=\lim_{t \ri 0^+}\si(t)$.}$. Then, there exists $\tau \in C^2(0,\infty)$ such that  
\[
\si(t)\, \leq \,\tau(t)\, \leq \, 8 \si(4t),\ \ t>0.
\]
}

\noi \textbf{Proof.} We set
 \[
\tau(t) \ := \ \frac{1}{2t^2}\int_0^{4t}\int_0^r\si(s)\, ds\, dr.
 \]
Since both $\si$ and $r\mapsto \int_0^r\si(s)ds$ are increasing, we have
\begin{align}
\tau(t)\, &\leq \, \frac{4t}{2t^2}\left(\sup_{0<r<4t}\int_0^r \si(s)ds\right) \nonumber\\
&=\, \frac{2}{t}\int^{4t}_0 \si(s)ds \nonumber\\
&\leq \, 8\, \si (4t). \nonumber
\end{align}
On the other hand, we have
\begin{align}
\int_0^{4t}\int_0^r \si(s)ds\, dr \, &\geq \, \int_{2t}^{4t}\int_0^r \si(s)ds\, dr\nonumber\\
&\geq \, 2t\int_0^{2t} \si(s)ds\nonumber\\
&\geq \, 2t\int_t^{2t} \si(s)ds\nonumber\\
&\geq \, 2t^2\si(t),  \nonumber
\end{align}
and hence $\tau \geq \si$. The lemma follows.     \qed
\ms

We now establish that Viscosity Solutions are consistent with classical solutions for degenerate elliptic PDE (and hence degenerate parabolic PDE) at points of twice differentiability. 

\ms

\noi \textbf{Theorem 10 (Compatibility of Viscosity with classical notions).} \emph{Let $\Om\sub \R^n$ be open and consider the degenerate elliptic PDE
\[
F\big(\cdot,u,Du,D^2u\big)\, =\, 0
\]
(that is $F$ is continuous on $\Om \by \R \by \R^n \by \mS(n)$ and satisfies \eqref{2.2}). If $u \in C^0(\Om)$ is a Viscosity Solution on $\Om$, then $u$ solves the PDE classically at points of twice differentiability. Conversely, if $u$ is twice differentiable and satisfies the PDE pointwise on $\Om$, then $u$ is a Viscosity Solution on $\Om$.
}
\ms

\noi \textbf{Proof.}  Let $u \in C^0(\Om)$ be a Viscosity Solution and suppose that $u$ is twice differentiable at some $x\in \Om$. Then, by (c) of Lemma 7, we have that 
\[
(Du(x),D^2u(x))\in \J^{2,\pm}u(x).
\]
Hence, by \eqref{2.8}, \eqref{2.9}, we have
\[
0\, \leq \, F \big(x,u(x),Du(x),D^2u(x)\big)\, \leq\, 0
\]
and hence $u$ solve the PDE pointwise at $x$.

Conversely, assume that $u$ is a twice differentiable and solves the PDE classically. Then the motivation of the definition of Viscosity Solutions in Chapter 1 implies that $u$ is a Viscosity Solution on $\Om$. 

For the reader's convenience, we provide an alternative proof of the latter fact which uses the language of Jets. Fix $x\in \Om$ such that $u$ is twice differentiable at $x.$ Then, by (c) of Lemma 7 we have that $\J^{2,+}u(x) \neq \emptyset$ and for any $(p,X) \in \J^{2,+}u(x)$, we have
\[
p\, =\, Du(x) \text{ in }\R^n,\ \ X\,\geq\, D^2u(x) \text{ in }\mS(n). 
\]
Hence, in view of degenerate ellipticity (see \eqref{2.2}), we have
\begin{align}
0\, &\leq \, F \big(x,u(x),Du(x),D^2u(x)\big) \nonumber\\
&=\, F \big(x,u(x),p,D^2u(x)\big) \nonumber\\
&\leq \, F \big(x,u(x),p,X\big). \nonumber
\end{align}
Hence, $u$ is a Viscosity Solution.  The Theorem ensues.        \qed

\ms

\noi \textbf{Remark 11.} 

(i) \textit{The only point that degenerate ellipticity is actually needed in what we have established so far is in the proof of Theorem 10, that is, for the consistency of viscosity and classical solutions. This is the main reason that viscosity solutions apply primarily to degenerate elliptic and degenerate parabolic PDE, namely, to PDE which support the Maximum Principle.}

\ms

(ii) \textit{Lemma 7 and Theorem 10 answer the fundamental questions (a)-(c) posed in the beginning of this Chapter. In particular, super-Jets and sub-Jets of any continuous function are non-empty on dense subsets of the domain, although these sets might be disjoint for an arbitrary continuous (and perhaps nowhere differentiable) function. Splitting the notions to 1-sided halves is essential: in general, we can not achieve having both upper and lower second derivatives (equivalently, touching from above and below by smooth functions in the same time), unless the function is $C^{1,1}$.}

\ms

When dealing with Viscosity Solutions, it is rarely the case that we have to actually calculate Jets of functions. This is in accordance with the fact that when dealing with classical solutions (e.g.\ Schauder theory) we almost never have to calculate derivatives, and when dealing with weak solutions of divergence PDE, we almost never have to calculate weak derivatives. For the sake of clarity, we do however give an explicit calculation of Jets.

\ms

\noi \textbf{Example 12 (Calculations of Jets).} Consider the cone function 
\[
C(z)\, :=\, -|z|, \ \ z\in \R^n.
\]
Then, we have 
\begin{align}
\J^{2,+}C(0)\, &=\, \Big\{\mB_1(0)\by \mS(n)\Big\} \bigcup \big\{(p,X)\ \big| \ |p|=1,\  X:p\ot p\geq 0 \big\},
\label{2.21}\\
\J^{2,-}C(0)\, &=\, \emptyset, \label{2.22}\\
\J^{2,\pm}C(x)\, &=\, \left\{-\left(\frac{x}{|x|},\frac{1}{|x|}\Big(I-\frac{x\ot x}{|x|^2}\Big)\pm A \right)\ :\ A\geq 0\right\}, \ \ x\neq 0.  \label{2.23}
\end{align}
These formulas can be calculated analytically, by arguing as in Lemma 7. For the reader's convenience, we provide an alternative geometric more intuitive proof. The figures below depict our claims in the case of 1 and 2 dimensions. 
\[
\includegraphics[scale=0.22]{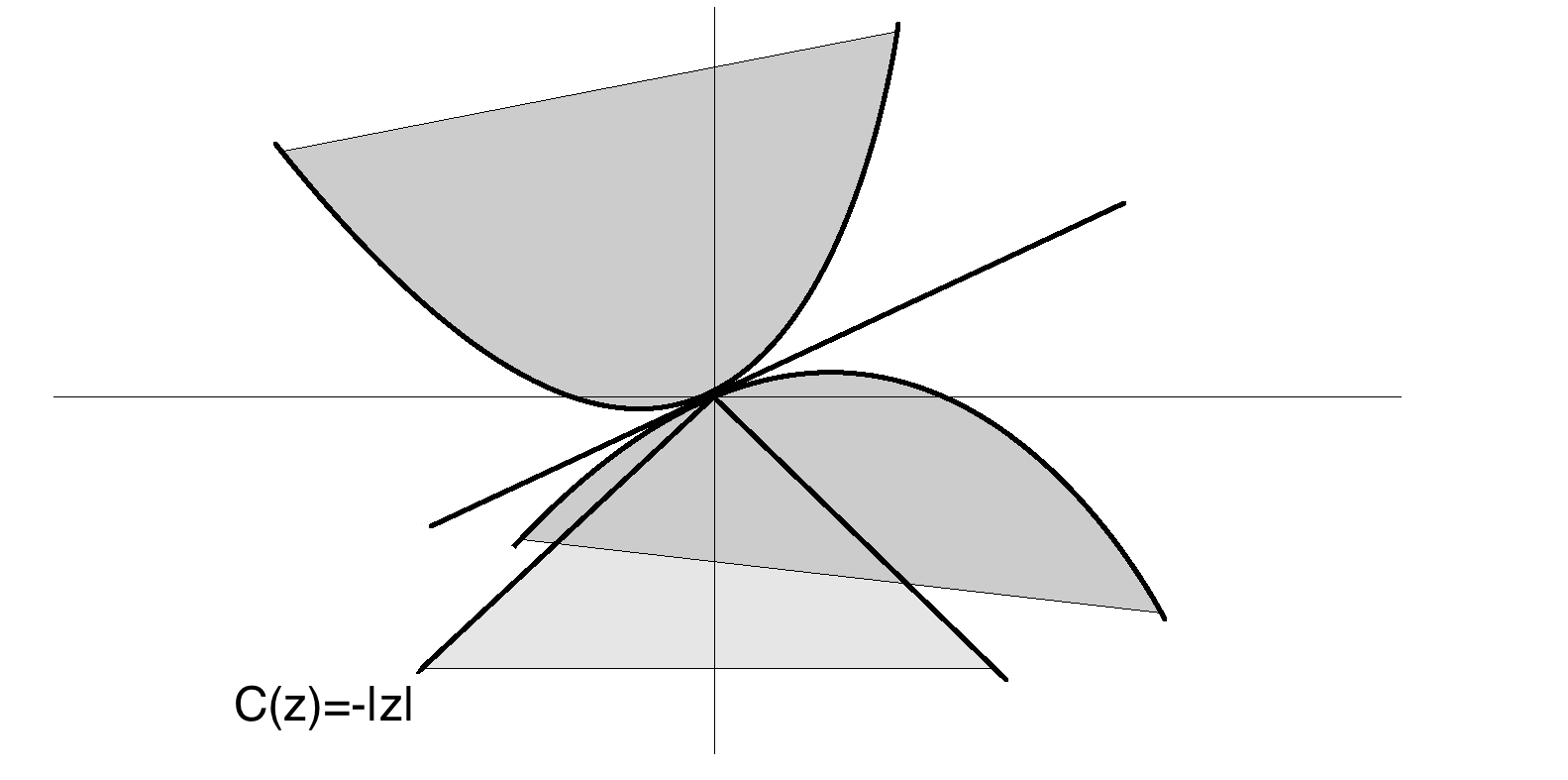}
\]
\[
\includegraphics[scale=0.22]{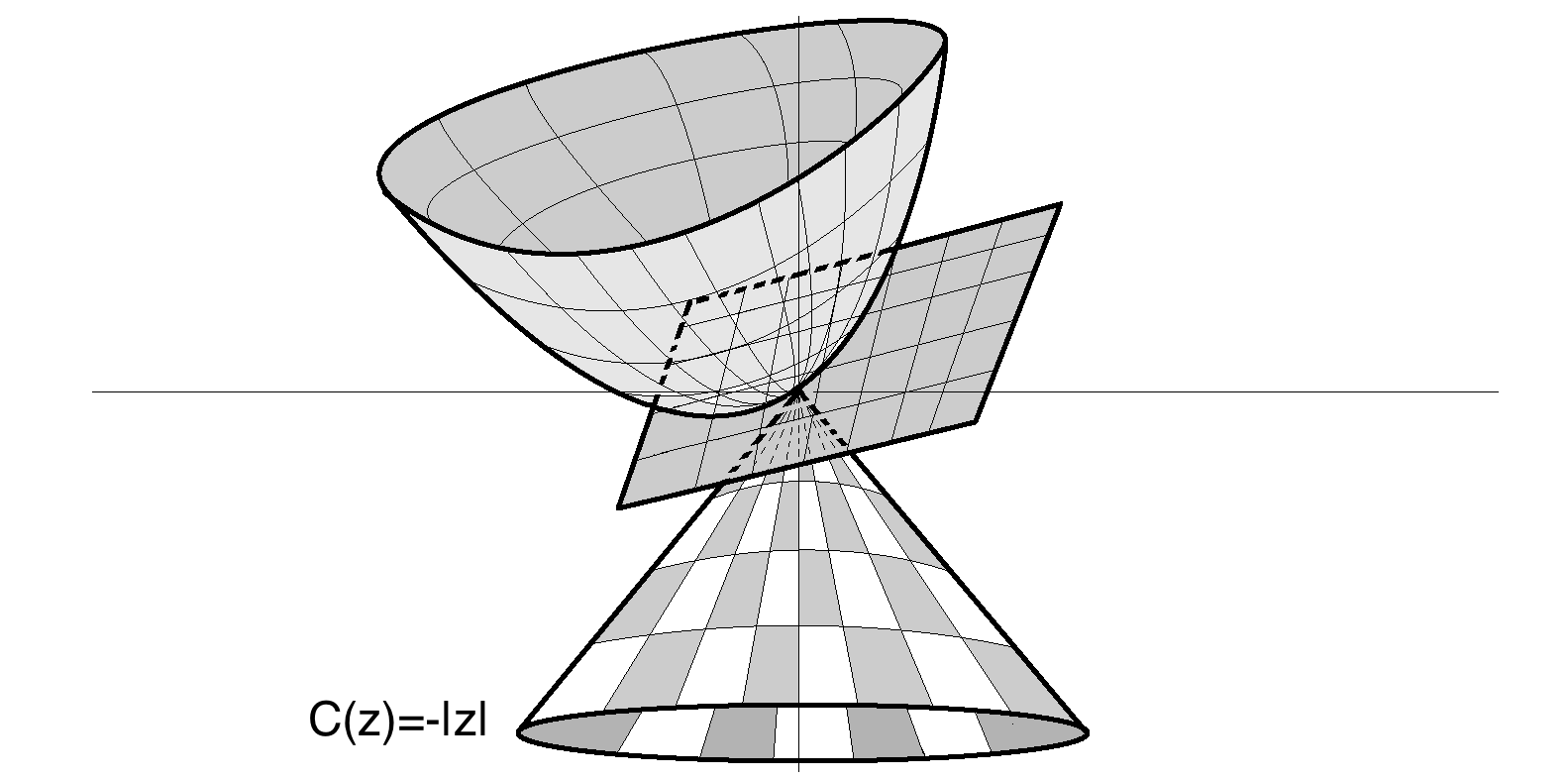}
\]
Formula \eqref{2.23} follows from Lemma 7 (c), since $C$ is smooth off the origin. Formula \eqref{2.22} follows immediately from Theorem 8 and the observation that there is no smooth function $\psi$ touching from below the cone at zero. Every hyperplane that passes through the origin with slope at most 1, lies above the cone and hence $\J^{1,+}C(0)$ equals the closed unit ball $\overline{\mB_1(0)}$. If $p\in \J^{1,+}C(0)$ and $|p|<1$, then any such hyperplane determined by $p$ touches the cone only at the origin, and hence any smooth function $\psi$ having this hyperplane as its tangent as zero lies above the cone near the origin. As a result, $\mB_1(0)\by \mS(n) \sub \J^{2,+}C(0)$. If on the other hand  $p\in \J^{1,+}C(0)$ and $|p|=1$, then the hyperplane touches the cone along some generator. In such an event, any function $\psi$ touching from above the cone at zero with this hyperplane as its tangent must be convex (near zero) in the direction of $p$. In view of Theorem 8, we conclude that \eqref{2.23}  ensues.

\ms

\ms

\ms

\noi \textbf{Remarks on Chapter 2.} Our exposition has been more thorough than standard texts on Viscosity Solutions, which are more advanced and omit the details we have analysed. However, the handbook of Crandall-Ishii-Lions \cite{CIL} remains the best reference on the subject, as well as the article of Crandall \cite{C2} and should be consulted as a further reading. The author has also greatly benefited from the lecture notes of Koike \cite{Ko} and Dragoni \cite{Dr}, as well as from those of Mallikarjuna Rao \cite{MR}, which have similar but different viewpoints on the subject and could be read parallel to these notes.

\chapter[Stability, Approximation, Existence]{Stability Properties of the notions and Existence via Approximation} 

In the previous two chapters we introduced and studied some basic properties of Viscosity Solutions $u\in C^0(\Om)$ of fully nonlinear degenerate elliptic PDE of the general form
\beq \label{3.1}
F \big(\cdot,u,Du,D^2u \big)\, =\, 0,
\eeq
where $\Om\sub \R^n$ and the nonlinearity
\[
F\in C^0\big(\Om \by \R \by \R^n \by \S(n) \big) 
\]
is monotone in the sense of \eqref{2.2}, that is that the function $X \mapsto F(x,r,p,X)$ is non-decreasing. An equivalent way to say that is the more elegant inequality
\beq \label{3.1a}
(X-Y)\big( F(x,r,p,X)-F(x,r,p,Y)\big)\, \geq\, 0.
\eeq
In this chapter we are concerned with limiting operations of Viscosity Solutions. It is a remarkable fact that  

{
\center{\fbox{\parbox[pos]{300pt}{
\ms

VISCOSITY SOLUTIONS ARE EXTREMELY STABLE UNDER PASSAGE TO SEVERAL TYPES OF LIMITS.

\ms
}}}

\ms\ms
}

\noi In particular, one of the most outstanding properties of this ``nonlinear pointwise theory of distributions" is that 
{
\center{\fbox{\parbox[pos]{330pt}{
\ms

\emph{in order to prove existence of solution for a PDE, it suffices to construct a sequence of approximate PDE whose solutions are locally bounded and equicontinuous, without control on derivatives!} 

\ms
}}}

\ms\ms
}

\noi Roughly speaking, the reason for this kind of strong stability is that ``we have passed all derivatives to test functions". The heuristic is that 

{

\ms

\center{\textit{``the less a priori regularity on the solutions' behalf we require in a ``weak" theory, the more flexible under limit operations the solutions are."}}

\ms

}
\noi Before proceeding, we record a couple of handy observations regarding the ``duality" interplay between Jets and maxima, that will be needed in the sequel.

\ms

\noi \textbf{Remark 1.} (a)
\textit{Every Jet $(p,X)\in \J^{2,+}u(x)$ of $u\in C^0(\Om)$ at $x\in \Om \sub \R^n$ can be realised by more than one $\psi \in C^2(\R^n)$ for which $u-\psi \leq 0=(u-\psi)(x)$ near $x$. It is only the values of the derivatives $(D\psi(x),D^2\psi(x))$ that matter.}

\ms

(b) \textit{Given any Jet $(p,X)\in \J^{2,+}u(x)$ of $u\in C^0(\Om)$ at $x\in \Om \sub \R^n$ (or equivalently any $\psi \in C^2(\R^n)$ for which $u-\psi \leq 0=(u-\psi)(x)$ near $x$), we may always assume that the Jet $(p,X)$ is realised by a function $\bar{\psi} \in C^2(\R^n)$ for which 
\[
u-\bar{\psi}  \, <\,  0\, =\, (u-\psi)(x)\  \text{ \underline{strictly} on a punctured ball } \  \mB_r(x)\set \{x\}.
\]
This can be easily achieved by setting $\bar{\psi} (z):= \psi(z)+|z-x|^4$ and restricting the neighbourhood of $x$ appropriately. Indeed, we have $\bar{\psi}\geq \psi$ on $\R^n$ and also $(D\bar{\psi}(x),D^2\bar{\psi}(x)) = (D\psi(x),D^2\psi(x))$.
\[
\includegraphics[scale=0.22]{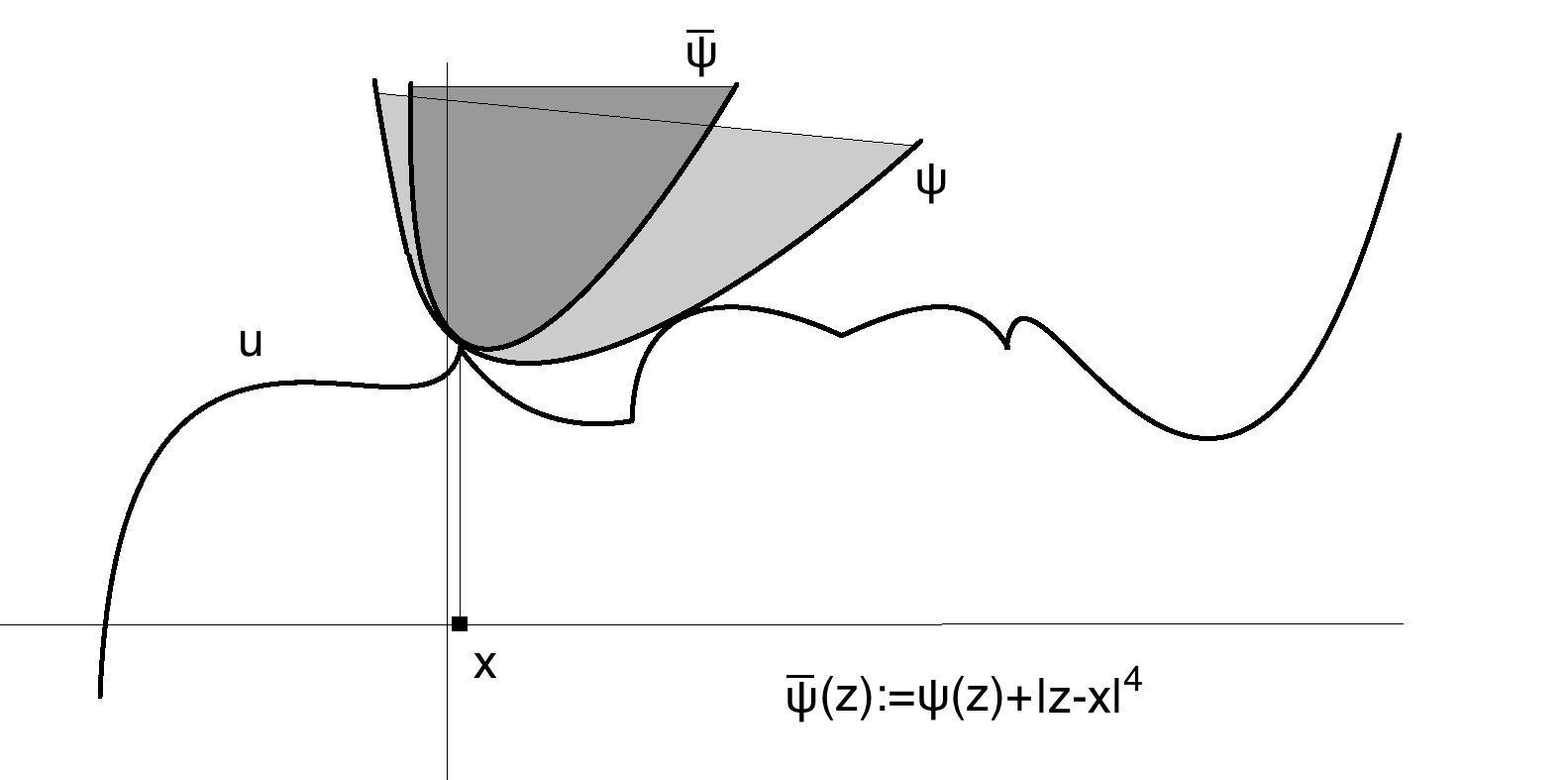}
\]
} 
Viscosity solutions are based on the coupling of two 1-sided notions, that of subsolutions and that of supersolutions. Passage to limits in the case of solutions is possible under uniform (local) convergence, but passage to limits is also possible in the case of subsolutions/supesolutions separately, \textit{and actually  under weaker 1-sided convergence}. Since the 2-sided case of solutions is somewhat simpler in both statements and proofs, we decouple the 1-sided and 2-sided stability results. 

The reason for this repetition is merely pedagogical: 1-sided stability for both sub- and super- solutions, implies stability of solutions.

\ms 

\noi {\underline{\textbf{Case 1:}}

{

\center{\textbf{2-SIDED STABILITY UNDER LOCALLY UNIFORM CONVERGENCE.} }

}

\ms

\noi We begin with the next fundamental result:

\ms

\noi \textbf{Theorem 2 (Stability under $C^0$ perturbations).} \emph{Let $\Om\sub \R^N$ and let also $\{F_m\}_1^\infty$, $F$ be in $C^0\big(\Om \by \R \by \R^n \by \S(n) \big) $ such that
\[
F_m \larrow F
\]
locally uniformly as $m\ri \infty$. Suppose that $\{u_m\}_1^\infty \sub C^0(\Om)$ is a sequence of Viscosity (Sub-/Super-) Solutions of 
\[
F_m\big(\cdot,u_m,Du_m,D^2u_m\big)\, =\, 0,
\]
on $\Om$. If 
\[
u_m \larrow u
\]
locally uniformly to some $u\in  C^0(\Om)$ as $m\ri \infty$, then the limit function $u$ is a Viscosity (Sub-/Super-) Solution of the limit PDE
\[
F\big(\cdot,u,Du,D^2u\big)\, =\, 0,
\]
on $\Om$. }

\ms

In practical situations of PDE theory, the standard way to prove convergence $u_m\ri u$ of the approximate solutions is by deriving \emph{appropriate estimates}. When treating a PDE in the Viscosity sense, it follows that

{
\center{\fbox{\parbox[pos]{320pt}{
\ms

very weak estimates (e.g.\ H\"older $C^\al$) suffice to pass to limits !

\ms
}}}

\ms\ms
}

By $C^\al$ bounds, we mean bounds on the functions \textit{and not $C^{1,\al}$, that is not on any kind of derivatives!} Indeed, directly from Theorem 1 and application of the classical Ascoli-Arz\'ela compactness theorem, we have the next consequence:
\ms

\noi \textbf{Corollary 3 (Bounds for Stability of Viscosity Solutions).} \emph{Let $\Om\sub \R^N$ and let also $\{F_m\}_1^\infty$, $F$ be in $C^0\big(\Om \by \R \by \R^n \by \S(n) \big) $ such that
\[
F_m \larrow F
\]
locally uniformly as $m\ri \infty$. Suppose that $\{u_m\}_1^\infty \sub C^0(\Om)$ is a sequence of Viscosity (Sub-/Super-) Solutions of 
\[
F_m\big(\cdot,u_m,Du_m,D^2u_m\big)\, =\, 0,
\]
on $\Om$. Then:} 

(a) \emph{If $\{u_m\}_1^\infty \sub C^0(\Om)$ is locally bounded and equicontinuous on $\Om$, then there exists $u\in  C^0(\Om)$ such that $u_m \ri u$ as $m\ri \infty$ along a subsequence and $u$ is a Viscosity (Sub-/Super-) Solution of the limit PDE
\[
F\big(\cdot,u,Du,D^2u\big)\, =\, 0,
\]
on $\Om$. }

(b) \emph{If there exists an $\al>0$ such that for any $\Om'\sub \Om$ the sequence $\{u_m\}_1^\infty$ satisfies the local H\"older estimate
\[
\sup_{m\in N}\|u_m\|_{C^\al(\Om')}\, \leq \, C(\Om'),\ \ \ \Om' \Subset \Om,
\]
then the same conclusion as in (a) above holds.}

\ms

\noi \textbf{Remark 4.} (a) \textit{We recall that the H\"older $C^\al$ semi-norm on $C^0(\Om)$ is defined by
\[
\|u\|_{C^\al(\Om)}\, :=\, \sup_{\Om}|u|\, +\, \sup_{x,y \in \Om,x\neq y}\frac{|u(x)-u(y)|}{|x-y|^\al}.
\]
It is evident that such a H\"older bound (as in (b) of the Corollary above) implies equicontinuity and local boundedness, since
\[
|u_m(x)|\, \leq \, C ,\ \ \ |u_m(x)-u_m(y)| \leq \, C|x-y|^\al, \text{ for }x,y \in \Om',
\]
for $C:=\sup_{m\in \N}\|u\|_{C^\al(\Om')}$.}

(b) \textit{In no other context such a strong stability result holds (for non-classical solutions), except for the linear theory of Distributions (Generalised Functions) of linear equations with \emph{smooth} coefficients!}

\ms

\noi \textbf{The idea.} The proof of Theorem 1 is based on the next lemma, whose essence is the Calculus fact that

{
\center{\fbox{\parbox[pos]{244pt}{
\ms

MAXIMA OF $C^0$ FUNCTIONS PERTURB TO 

\ \ MAXIMA UNDER $C^0$ PERTURBATIONS.

\ms
}}}

\ms\ms
}

\noi \textbf{Lemma 5 (Uniform perturbations of maxima).} \emph{Let $\Om\sub \R^n$ and suppose $u,\{u_m\}_1^\infty \sub C^0(\Om)$ such that $u_m \ri u$ locally uniformly on $\Om$. If $x\in \Om$ and $\psi \in C^2(\R^n)$ are such that
\[
u-\psi\, <\, 0\,=\, (u-\psi)(x)
\]
near $x$, then there exist $\{x_m\}_1^\infty \sub \Om$ and $\psi_m \in C^2(\R^n)$  such that 
\begin{align}
&x_m \larrow x,\ \text{ in }\Om, \nonumber\\
&\psi_m \larrow \psi,\ \text{ in }C^2(\R^n),\nonumber
\end{align}
as $m\ri \infty$, with the property
\[
u_m-\psi_m\, \leq\, 0\,=\, (u-\psi_m)(x),
\]
near $x$. Moreover, $\psi_m$ can be chosen to be a shift of $\psi$ by a constant $a_m \ri 0$: $\psi_m=\psi+a_m$.
}

\ms

In view of the preceding remark, the strictness of the maximum can always be achieved given a (perhaps non-strict) maximum. Before giving the proof of Lemma 5, we reformulate it the language of our generalised derivatives.

\ms

\noi \textbf{Lemma 5* (Stability of Jets).} \emph{Let $\Om\sub \R^n$ and suppose $u,\{u_m\}_1^\infty \sub C^0(\Om)$ such that $u_m \ri u$ locally uniformly on $\Om$. If $x\in \Om$ and $(p,X)\in \J^{2,+}u(x)$, then there exist $\{x_m\}_1^\infty \sub \Om$ and 
\[
(p_m,X_m)\in \J^{2,+}u_m(x_m)
\]
such that 
\[
(x_m,p_m,X_m) \larrow (x,p,X)
\]
as $m\ri \infty$.
}

\ms

\noi \textbf{Proof of Lemmas 5, 5*.} Fix $R>0$ smaller than $\frac{1}{2}\dist(x,\p \Om)$, in order to achieve $\overline{\mB_R(x)} \sub \Om$. By assumption we have that $u_m \ri u$ in $C^0(\overline{\mB_R(x)})$ as $m\ri \infty$, and also
\[
u-\psi\, <\, 0\,=\, (u-\psi)(x),
\]
on $\mB_R(x) \set \{x\}$, by decreasing $R$ further if necessary. By compactness and continuity, we can choose a point $x_m \in \overline{\mB_R(x)}$ such that 
\[
\max_{\overline{\mB_R(x)}}\, \{u_m-\psi\}\, =\, u_m(x_m)-\psi(x_m).
\]
Hence, we have
\beq \label{3.2}
u_m-\psi\, \leq\,  (u_m-\psi)(x_m),
\eeq
on $\overline{\mB_R(x)}$. We set
\[
a_m \, :=\, u_m(x_m)-\psi(x_m).
\]
By considering the shift 
\[
\psi_m \, :=\, \psi \, +\, a_m
\]
of $\psi$, we have
\[
u_m-\psi_m\, \leq\, 0\,=\, (u-\psi_m)(x),
\]
on $\overline{\mB_R(x)}$. Since the maximum of $u-\psi$ is strict on $\mB_R(x)$, we have
\[
\max_{\p \mB_R(x)}\, \{u-\psi\}\, <\, \sup_{\mB_R(x)}\, \{u-\psi\}.
\]
Since $u_m \ri u$ uniformly on $\mB_R(x)$, we have that for $m$ large enough the point $x_m$ of maximum of $u_m-\psi_m$ lies strictly inside the ball $\mB_R(x)$.
\[
\includegraphics[scale=0.22]{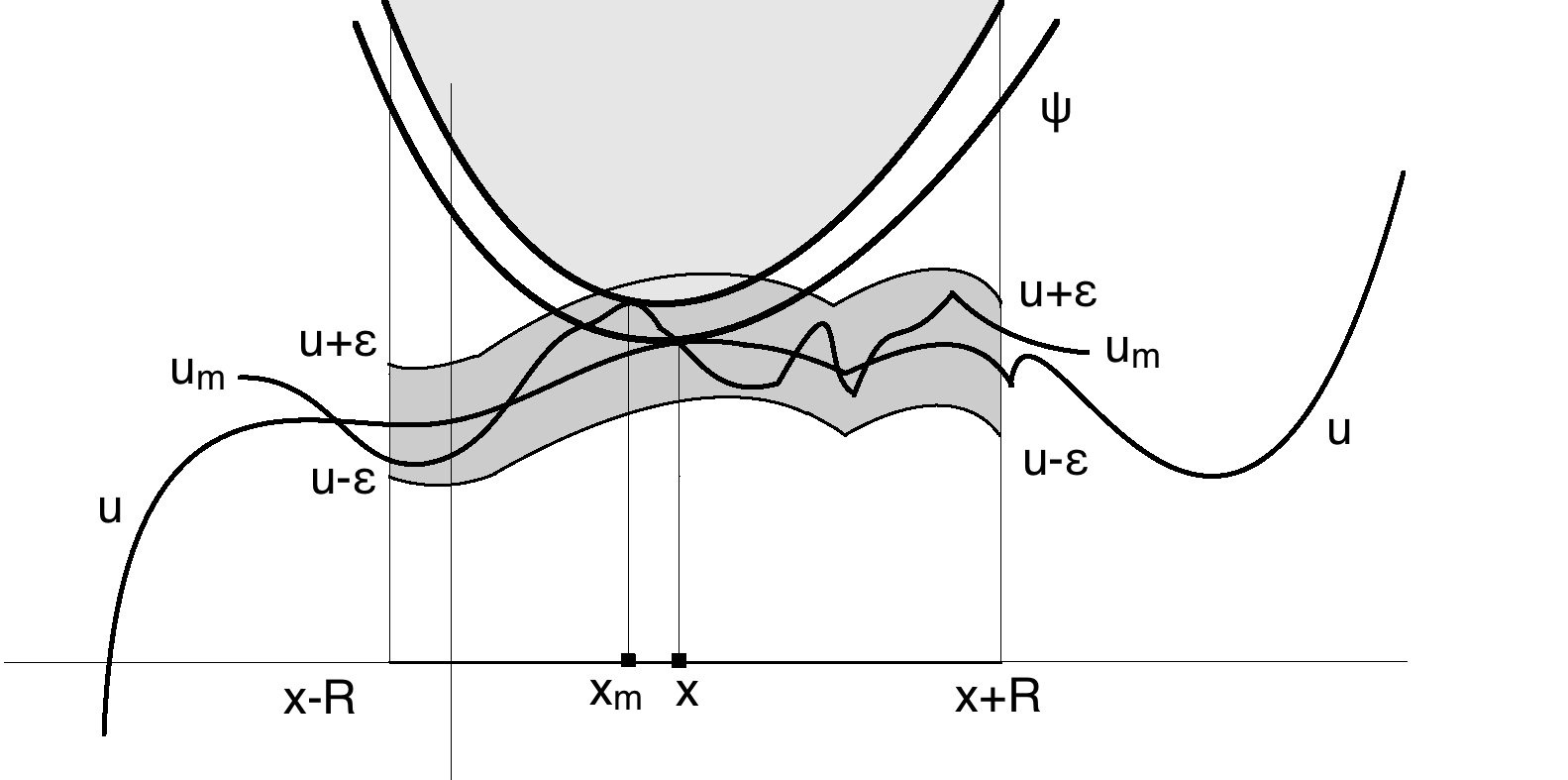}
\]
By evaluating  \eqref{3.2} at $x$ we have
\beq \label{3.3}
(u_m-\psi)(x)\, \leq\,  (u_m-\psi)(x_m).
\eeq
By compactness, there exists $x^* \in \mB_R(x)$ such that $x_m \ri x^*$, along a subsequence as $m\ri \infty$, which we denote again by $x_m$. By passing to the limit in \eqref{3.3} we have
\[
(u-\psi)(x)\, \leq\,  (u-\psi)(x^*).
\]
Since $x$ is the unique maximum of $u-\psi$ in $\mB_R(x)$, we obtain that $x^*=x$ and hence $x_m \ri x$ as $m\ri \infty$. 
\[
\includegraphics[scale=0.22]{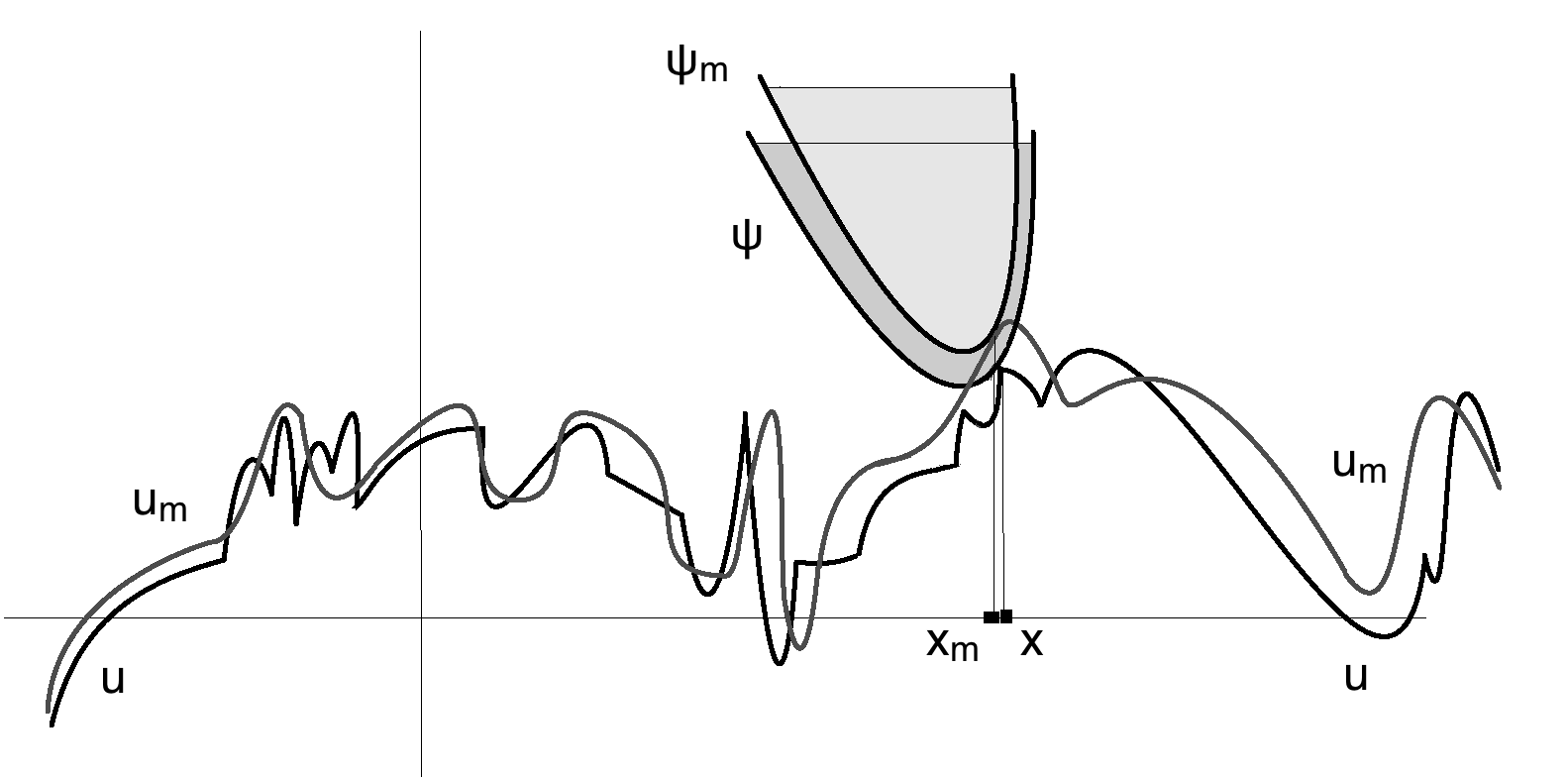}
\]
Hence, Lemma 5 follows. Lemma 5* follows by invoking Theorem 8 of Chapter 2. \qed   

\ms

\noi \textbf{Proof of Theorem 2.} Let $(p,X) \in \J^{2,+}u(x)$ for some $x\in \Om$. Since $u_m \ri u$ in $C^0(\Om)$, by Lemma 5* there exists $x_m \ri x$ and
\[
(p_m,X_m) \in \J^{2,+}u_m(x_m)
\]
such that $(p_m,X_m)\ri (p,X)$, as $m\ri \infty$. Since $u_m$ is a Viscosity Solution of
\[
F_m\big(\cdot,u_m,Du_m,D^2u_m\big)\, \geq\, 0
\] 
on $\Om$, we have
\[
F_m\big(x_m,u_m(x_m),p_m,X_m\big)\, \geq\, 0.
\] 
Since $F_m\ri F$ locally uniformly on $\Om \by \R \by \R^n \by\mS(n)$ and $u_m \ri u$ locally uniformly on $\Om$, by passing to the limit as $m\ri \infty$ we get
\[
F\big(x,u(x),p,X\big)\, \geq\, 0.
\] 
Hence, $u$ is a Viscosity Solution of
\[
F\big(\cdot,u,Du,D^2u\big)\, \geq\, 0
\] 
on $\Om$. The supersolution property follows by arguing symmetrically. The proof of the theorem in complete.      \qed

\ms

\noi \textbf{An application: Existence for the Dirichlet Problem.} The outstanding stability properties of Viscosity Solutions under very weak bounds imply that we can always solve a  PDE, once we construct an approximate PDE whose solutions are precompact in $C^0(\Om)$. Directly from Theorem 2 and Corollary 3, we have the next

\ms

\noi \textbf{Corollary 6 (Solvability of the Dirichlet Problem via Stability)}. \emph{Consider a bounded open set  $\Om \Subset \R^n$ and $b\in C^0(\overline{\Om})$. Let $F$ be continuous and degenerate elliptic on $\Om \by \R \by \R^n \by \mS(n)$ in the sense of \eqref{3.1a}.
Consider the Dirichlet problem
\beq \label{3.4}
\left\{
\begin{array}{l}
F\big(\cdot,u,Du,D^2u\big)\,=\,0, \, \text{ in }\Om,\\
\hspace{77pt}u\, =\, b, \text{ on }\p \Om.
\end{array}
\right.
\eeq
Assume that there exist approximate problems 
\beq \label{3.5}
\left\{
\begin{array}{l}
F_m(\cdot,u_m,Du_m,D^2u_m)\,=\,0, \  \ \text{ in }\Om,\\
\hspace{99pt}u_m\, =\, b_m, \text{ on }\p \Om,
\end{array}
\right.
\eeq
for coefficients $\{F_m\}_1^\infty \sub C^0\big( \Om \by \R \by \R^n \by \mS(n) \big)$ and boundary conditions  $\{b_m\}_1^\infty \sub C^0(\overline{\Om})$, such that
\begin{align}
F_m &\larrow F  \ \text{ in }C^0\big( \Om \by \R \by \R^n \by \mS(n) \big), \nonumber\\
b_m & \larrow b \ \ \text{ in }C^0(\overline{\Om}), \nonumber
\end{align}
as $m\ri \infty$. If the problems \eqref{3.5} have Viscosity Solutions $\{u_m\}_1^\infty \sub C^0(\overline{\Om})$ which are bounded and equicontinuous on $\Om$, then the problem \eqref{3.4} has a Viscosity Solution $u \in C^0(\overline{\Om})$, which is addition can be approximated by $u_m$.
}
\ms

In Chapter 5, we will utilise Corollary 6 in order to prove existence of $\infty$-Harmonic functions with given boundary values, that is of solution to the Dirichlet Problem for the $\infty$-Laplacian
\[
\De_\infty u \, :=\, Du \ot Du :D^2 u\, =\, 0.
\]
The apt approximating sequence in the case at hand is the sequence of $p$-Harmonic functions, that is of solutions to 
\[
\De_p u \, :=\, \Div \big( |Du|^{p-2}Du \big)\, =\, 0.
\]
as $p\ri \infty$. Now we continue with the 

\ms 

\noi {\underline{\textbf{Case 2:}}

{

\center{\textbf{1-SIDED STABILITY UNDER $\Gamma$-CONVERGENCE.} }

}

\ms

$\Gamma$-Convergence is a type of 1-sided convergence introduced by Ennio De Giorgi in the infinite-dimensional context of Calculus of Variations, as a notion of convergence whose job is to pass minimisers of (a sequence of) functionals to minimisers of the limit functional.

In the 1-sided results that follow, we will only consider the result for subsolutions. Symmetric results holds for supersolutions as well, in the obvious way. We leave the statements and proofs as a simple exercise for the reader.

In order to handle the 1-sided case efficiently, we recall from Calculus the next
\ms

\noi \textbf{Definition 7 (Semicontinuity).} \emph{Let $\Om \sub \R^n$ and $u : \Om \ri \R$. We say that $u$ is upper semicontinuous on $\Om$ and write $u\in USC(\Om)$, when
\[
 \underset{y\ri x}{\lim\sup}\, u(y) \, \leq\, u(x),
\]
for all $x\in \Om$. Symmetrically, we say that $u$ is lower semicontinuous on $\Om$ and write $u\in LSC(\Om)$, when
\[
\underset{y\ri x}{\lim\inf}\, u(y) \, \geq\, u(x),
\]
for all $x\in \Om$. }
\ms

The formal definition of $\Gamma$-Convergence, adapted to our case, is

\ms

\noi \textbf{Definition 8 ($\pm\Gamma$-Convergence).} \emph{Let $u,\{u_m\}_1^\infty : \Om \sub \R^n \larrow \R$. }

\emph{(i) We say that \underline{$u_m$ plus-Gamma converges to $u$ on $\Om$} and write
\[
u_m \, \overset{+\Gamma\ }{\larrow}\, u, \ \text{ as }m\ri \infty,
\]
when for all $x\in \Om$, the following hold:}

\noi \emph{(a) for \underline{any} sequence $x_m \ri \infty$ as $m\ri \infty$, we have
\[
\underset{m\ri \infty}{\lim\inf}\, u_m(x_m) \, \geq\,  u(x).
\]
(b) there \underline{exists} a sequence $\hat{x}_m \ri \infty$ as $m\ri \infty$, we have
\[
\underset{m\ri \infty}{\lim}\, u_m(\hat{x}_m) \, =\,  u(x).
\]
}

\noi \emph{(ii) We say that \underline{$u_m$ minus-Gamma converges to $u$ on $\Om$} and write
\[
u_m \, \overset{-\Gamma\ }{\larrow}\, u, \ \text{ as }m\ri \infty,
\]
when $-u_m$ plus-gamma converges to $-u$ on $\Om$, that is when for all $x\in \Om$, the following hold:}

\noi \emph{(a) for \underline{any} sequence $x_m \ri \infty$ as $m\ri \infty$, we have
\[
\underset{m\ri \infty}{\lim\sup}\, u_m(x_m) \, \leq\,  u(x).
\]
(b) there \underline{exists} a sequence $\hat{x}_m \ri \infty$ as $m\ri \infty$, we have
\[
\underset{m\ri \infty}{\lim}\, u_m(\hat{x}_m) \, =\,  u(x).
\]
}
\noi \emph{(iii) We shall say that $u_m$ plus(or minus)-Gamma converges to $u$ at $x$ if the conditions (a), (b) of (i) (or (b)) hold at $x$ only.}

\ms

\noi It is a trivial matter to check that 

{
\center{\fbox{\parbox[pos]{200pt}{
\ms

$\ \ u_m \larrow u$ in $C^0(\Om)$\ \ $\Big|\!\!\!\Longrightarrow$ \ \ $u_m \overset{\pm \Gamma\ }{\larrow} u$. 
\ms
}}}

\ms\ms
}

Having these convergence notions at hand, we may now proceed to the 1-sided analogues of Theorem 2 and Lemma 5. One final observation before the

\noi \textbf{Remark 9.}

{
\center{\fbox{\parbox[pos]{320pt}{
\ms

THE DEFINITION OF VISCOSITY SOLUTIONS REMAINS THE SAME IF WE SPLIT THE A-PRIORI CONTINUITY REQUIREMENT TO 2 HALVES, THAT IS, THAT SUBSOLUTIONS ARE UPPER-SEMICONTINUOUS AND SUPER-

\ \ \ \ \ \  SOLUTIONS ARE LOWER-SEMICONTINUOUS.

\ms
}}}

\ms\ms
}

This requirement makes no difference for the notion of solutions, but relaxes the a priori regularity requirements of sub-/super- solutions to the ones which suffice for 1-sided considerations. We shall use this slight modification of the definition without further notice.

\ms

\noi \textbf{Lemma 10 (1-sided Stability of Super-Jets).} \emph{Let $\Om\sub \R^n$ and suppose $u,\{u_m\}_1^\infty \sub USC(\Om)$ such that 
\[
u_m \, \overset{-\Gamma\ }{\larrow}\, u, \ \text{ as }m\ri \infty,
\]
on $\Om$. If $x\in \Om$ and $(p,X)\in \J^{2,+}u(x)$, then there exist $\{x_m\}_1^\infty \sub \Om$ and 
\[
(p_m,X_m)\in \J^{2,+}u_m(x_m)
\]
such that 
\[
(x_m,u_m(x_m),p_m,X_m) \larrow (x,u(x),p,X),
\]
as $m\ri \infty$.
}

\ms

\noi \textbf{Proof.} Fix $0<R<\frac{1}{2}\dist(x,\p \Om)$ such that $\overline{\mB_R(x)} \sub \Om$. By Theorem 8 of Chapter 2 and Remark 1, there is a $\psi \in C^2(\R^n)$ such that $(p,X)=(D\psi(x),D^2\psi(x))$ and
\[
u-\psi\, <\, 0\,=\, (u-\psi)(x),
\]
on $\mB_R(x) \set \{x\}$.  By compactness and semi-continuity of $u_m$, we can choose a point $x_m \in \overline{\mB_R(x)}$ such that 
\[
\max_{\overline{\mB_R(x)}}\, \{u_m-\psi\}\, =\, u_m(x_m)-\psi(x_m).
\]
Hence, we have
\beq \label{3.2a}
u_m-\psi\, \leq\,  (u_m-\psi)(x_m),
\eeq
on $\overline{\mB_R(x)}$. By compactness, there exists $x^* \in \mB_R(x)$ such that $x_m \ri x^*$, along a subsequence as $m\ri \infty$, which we denote again by $x_m$.  By assumption, there is a sequence $\hat{x}_m \ri x$ such that $u_m(\hat{x}_m) \ri u(x)$ as $m\ri \infty$. 
By evaluating \eqref{3.2a} at $\hat{x}_m$ and taking limsup, we have
\begin{align} \label{3.3a}
(u-\psi)(x)\, & = \, \underset{m\ri \infty}{\lim\sup} \, (u_m-\psi)(\hat{x}_m)  \nonumber\\
& \leq\,  \underset{m\ri \infty}{\lim\sup} \, (u_m-\psi)(x_m)\\
&\leq \, (u-\psi)(x^*). \nonumber
\end{align}
Since $x$ is the unique maximum of $u-\psi$ in $\mB_R(x)$, by \eqref{3.3a} we obtain that $x^*=x$ and hence $x_m \ri x$ as $m\ri \infty$. Consequently, for $m$ large, at the point $x_m$ we have an interior maximum in $\mB_R(x)$. By Taylor expansion of $\psi$ at $x_m$, we obtain
\[
(p_m,X_m)\, :=\, \big(D\psi(x_m),D^2\psi(x_m)\big) \in \J^{2,+}u_m(x_m),
\]
and by smoothness of $\psi$, we get that $(p_m,X_m) \ri (p,X)$ as $m\ri \infty$. Again by  \eqref{3.2a} evaluated at $\hat{x}_m$, we have 
\begin{align} 
u(x)\, & = \, \underset{m\ri \infty}{\lim\inf} \, \Big(u_m(\hat{x}_m)-\psi_m(\hat{x}_m)+\psi(x_m) \Big) \nonumber\\
& \leq\,  \underset{m\ri \infty}{\lim\inf} \, u_m(x_m)\nonumber\\
&\leq\,  \underset{m\ri \infty}{\lim\sup} \, u_m(x_m)\nonumber\\
&\leq\, u(x).
\end{align}
As a result, we get that $u_m(x_m)\ri u(x)$ as $m\ri \infty$. Putting it all together, Lemma 10 ensues.           \qed

\ms

Now we have the 1-sided analogue of Theorem 2. 
\ms

\noi \textbf{Theorem 11 (1-sided Stability of Subsolutions).} \emph{Let $\Om\sub \R^N$ and let also $F,\{F_m\}_1^\infty \sub C^0\big(\Om \by \R \by \R^n \by \S(n) \big)$. Suppose that $\{u_m\}_1^\infty \sub USC(\Om)$ is a sequence of Viscosity Subsolutions of  $F=0$ on $\Om$, i.e.
\[
F\big(\cdot,u_m,Du_m,D^2u_m\big)\,\geq\, 0.
\]
Suppose that $F_m \ri F$ locally uniformly on $\Om \by \R \by \R^n \by \S(n) $ as $m\ri \infty$, and
\[
u_m \, \overset{-\Gamma\ }{\larrow}\, u, \ \text{ as }m\ri \infty,
\]
on $\Om$, and that $u\in USC(\Om)$. Then, the limit $u$ is  a Viscosity Subsolution on $\Om$ of the limit PDE:
\[
F\big(\cdot,u,Du,D^2u\big)\, \geq\, 0.
\]
}

\noi \textbf{Proof.} Let $(p,X) \in \J^{2,+}u(x)$ for some $x\in \Om$. Since $u_m \overset{-\Gamma\ }{\larrow} u$, by Lemma 9 there exists $x_m \ri x$ and
\[
(p_m,X_m) \in \J^{2,+}u_m(x_m)
\]
such that $(u_m(x_m),p_m,X_m)\ri (u(x),p,X)$, as $m\ri \infty$. Since $u_m$ is an $USC(\Om)$ Viscosity Solution of
\[
F_m\big(\cdot,u_m,Du_m,D^2u_m\big)\, \geq\, 0
\] 
on $\Om$, we have
\[
F_m\big(x_m,u_m(x_m),p_m,X_m\big)\, \geq\, 0.
\] 
Hence,
\begin{align}
0\, &\leq \,  \underset{m\ri \infty}{\lim\sup} \, F_m\big(x_m,u_m(x_m),p_m,X_m\big)  \nonumber\\
&=\, \underset{m\ri \infty}{\lim\sup} \, \Big\{F_m\big(x_m,u_m(x_m),p_m,X_m\big)  - F_m\big(x,u_m(x_m),p,X\big) \Big\} \nonumber\\
&\ \ \ \ +\, \underset{m\ri \infty}{\lim\sup} \,  F_m\big(x,u_m(x_m),p,X\big) \nonumber\\
&\leq\,  \underset{m\ri \infty}{\lim\sup}  \max_{[u(x)-1,u(x)+1 ]}\Big\{ F_m\big(x_m, \cdot ,p_m,X_m\big) - F_m\big(x, \cdot ,p,X\big)\Big\}  \nonumber \\
&\ \ \ \ +\, \underset{m\ri \infty}{\lim\sup} \,  F_m\big(x,u_m(x_m),p,X\big) \nonumber\\
& = \, F\big(x,u(x),p,X\big). \nonumber
\end{align}
Hence, $u$ is a Viscosity Solution of
\[
F\big(\cdot,u,Du,D^2u\big)\, \geq\, 0
\] 
on $\Om$. The proof of the theorem in complete.       \qed

\ms

The following result is an important application of Theorem 11 which shows that Viscosity Subsolution are closed under the operation of pointwise supremum.
\ms

\noi \textbf{Theorem 12 (Suprema of Subsolutions is Subsolution).} \emph{Let $\Om\sub \R^N$ and let also $F \in C^0\big(\Om \by \R \by \R^n \by \S(n) \big)$. Suppose that $\mU \sub USC(\Om)$ is a family of Viscosity Subsolutions of  $F=0$ on $\Om$, i.e.
\[
F\big(\cdot,u,Du,D^2u\big)\,\geq\, 0, \ \ u \in \mU.
\]
We set
\[
U(x) \, :=\, \sup\big\{u(x)\ |\ u \in \mU \big\},\ \ x\in \Om, \ms
\]
and suppose that $U(x)<\infty$ for all $x\in \Om$  and that $U \in USC(\Om)$. Then, $U$ is  a Viscosity Solution of $F\geq 0$ on $\Om$.
}

\ms

\noi \textbf{Proof.} By the definition of $U$, we have $u_m \leq U$ on $\Om$ for any sequence $\{u_m\}_1^\infty \sub \mU$. Since $U \in USC(\Om)$, for each $x\in \Om$ and $x_m \ri x$ as $m\ri \infty$, we have
\[
\underset{m\ri \infty}{\lim\sup} \, u_m(x_m)\, \leq \, \underset{m\ri \infty}{\lim\sup} \, U(x_m) \, \leq\, U(x).
\]
On the other hand, for each $x\in \Om$ there is a sequence $\{u_m\}_1^\infty \sub \mU$ such that 
\[
u_m(x) \larrow U(x),\ \ \text{ as } m\ri \infty.
\]
Hence, for each $x\in \Om$ there is a sequence $\{u_m\}_1^\infty \sub \mU$ such that $u_m$ minus-Gamma converges to $U$ at $x$. Fix an $x\in \Om$ such that $(p,X)\in \J^{2,+}U(x)$ and consider the respective sequence  $\{u_m\}_1^\infty \sub \mU$. By Lemma 9, there exist $\{x_m\}_1^\infty \sub \Om$ and 
\[
(p_m,X_m)\in \J^{2,+}u_m(x_m)
\]
such that 
\[
(x_m,u_m(x_m),p_m,X_m) \larrow (x,u(x),p,X),
\]
as $m\ri \infty$. Since each $u_m$ is a Viscosity Solution of 
\[
F\big(\cdot,u_m,Du_m,D^2u_m\big)\,\geq\, 0,
\]
we have
\[
F\big(x_m,u_m(x_m),p_m,X_m \big)\, \geq\, 0.
\]
By passing to the limit as $m\ri \infty$, we have that $U$ is a Viscosity Solution as well.       \qed 

\ms

\ms

\ms

\noi \textbf{Remarks on Chapter 3.} The property of stability under limits is one of the most important features of the notion of viscosity solutions. For further reading, the reader should consult Crandall-Ishii-Lions \cite{CIL}, Crandall \cite{C2}, Lions \cite{L1, L2}, Barles \cite{Ba}, Bardi-Capuzzo Dolcetta \cite{BCD}, Ishii \cite{I1}, Souganidis \cite{So}. The part of the chapter on the stability under $1$-sided convergence follows closely Crandall-Ishii-Lions \cite{CIL}, but with many more details. An excellent introduction to the notion of $\Gamma$-convergence (in infinite dimensions!) is Braides \cite{Br}. Alternative methods to produce (non-unique!) strong a.e.\ solutions of fully nonlinear ``implicit" equations which are not based on stability (or the vanishing viscosity method)  can be found in Dacorogna \cite{D1, D2}.

\chapter[Mollification and Semiconvexity]{Mollification of Viscosity Solutions and Semiconvexity}

In the previous three chapters we defined and studied the basic analytic properties of Viscosity Solutions $u\in C^0(\Om)$ to fully nonlinear degenerate elliptic PDE of the form
\[
F\big(\cdot,u,Du,D^2u\big)\,=\, 0,
\]
where $\Om\sub \R^n$ and $F \in C^0(\Om \by \R \by \R^n \by\mS(n))$. One of the central difficulties when handling an equation in the viscosity sense is that the notion of solution itself is not sufficiently ``handy". It is a pointwise criterion and, apart from continuity, there is no further a priori regularity of solutions. Hence, it is crucial to invent suitable regularisations of Viscosity Solutions is order to be able to manipulate calculations with them satisfactorily. As usually, the heuristic is to \textit{prove a certain statement first for sufficiently a differentiable mollification and then pass to limits}.

\ms

\noi \textbf{Digression into Divergence-Structure Quasilinear PDE.} Before proceeding to the central theme of this chapter, for the sake of comparison we recall the analogous standard mollification scheme which applies to ``weak solutions", defined by means of integration by parts. 

Given $\e>0$ and $u : \R^n \ri \R$ Lebesgue measurable, the (standard) \emph{mollification by convolution} of $u$ is defined by
\[
u^\e(x)\, :=\, \int_{\R^n}u(y)\frac{1}{\e^n}\eta \left( \frac{|x-y|}{\e}\right)dy
\]
where $\eta \in C^\infty_c(\R)$, $\eta \geq 0$ and $\|\eta\|_{L^1(\R)}=1$. It is well known that $u^\e \in C^\infty(\R^n)$ and $u^\e \ri u$ a.e.\ on $\Om$, as $\e \ri 0$. This type of mollification is defined by means of integral convolution, i.e.
\[
u^\e\ = \ u * \eta^\e,\ \ \ \eta^\e(x)\,:=\, \frac{1}{\e^n}\eta \left( \frac{|x|}{\e}\right)
\]
and roughly this implies that it \emph{respects duality and weak derivatives defined via integration by parts}. The figure below gives the most popular choice of $\eta$ function.
\[
 \includegraphics[scale=0.18]{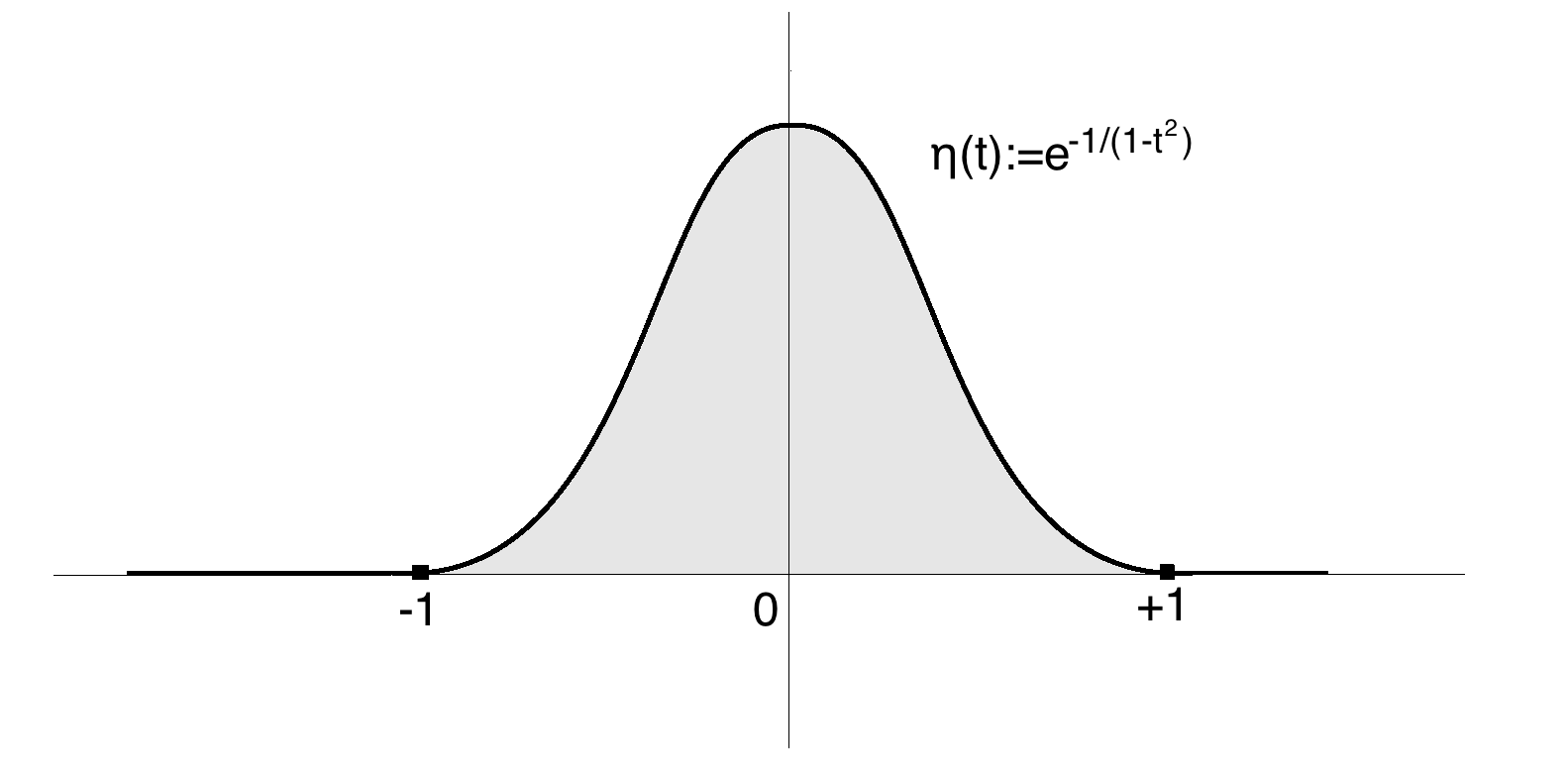} 
\]
It is a well known fact that $u^\e \ri u$ in the respective topology of the $L^p$ space or Sobolev $W^{k,p}$ space that $u$ lies into. For example,  if $u\in W^{1,p}_{\text{loc}}(\R^n)$ then $u^\e \ri u$ in $W^{1,p}_{\text{loc}}(\R^n)$ as $\e \ri 0$, etc.

\ms \ms

The apt mollification scheme for Viscosity Solutions which respects our generalised derivatives, namely the semi-Jets $\J^{2,\pm}$, is 1-sided. To a certain extent, this reflects the 1-sided nature of the notion of sub-/super- solutions. 

\ms

\noi \textbf{Definition 1 (Sup-/Inf- Convolution Approximations).} \emph{
Let $\Om\Subset \R^n$ and $\e>0$ be given. For $u\in C^0(\overline{\Om})$, we define
\begin{align}
u^\e(x)\, &:=\, \sup_{y\in \Om}\left\{u(y)\, -\, \frac{|x-y|^2}{2\e}\right\}, \ \ x\in \Om,\label{4.1}\\
u_\e(x)\, &:=\, \inf_{y\in \Om}\left\{u(y)\, +\, \frac{|x-y|^2}{2\e}\right\},  \ \ x\in \Om.\label{4.2}
\end{align}
We call \underline{$u^\e$ the sup-convolution of $u$} and \underline{$u_\e$ the inf-convolution of $u$}.
}

\ms 
\ms 

\noi \textbf{Remark 2.} \textit{Geometrically, the sup-convolution of $u$ at $x$ is defined as follows: we ``bend downwards" the graph of $u$ near $x$ by subtracting the paraboloid ${|\cdot - x|^2}/{2\e}$ which is centred at $x$. Then, $u^\e(x)$ is defined as the maximum of the ``bent" function
\[
y \mapsto u(y) \,- \, \frac{|x-y|^2}{2\e}.
\]}
\[
\includegraphics[scale=0.24]{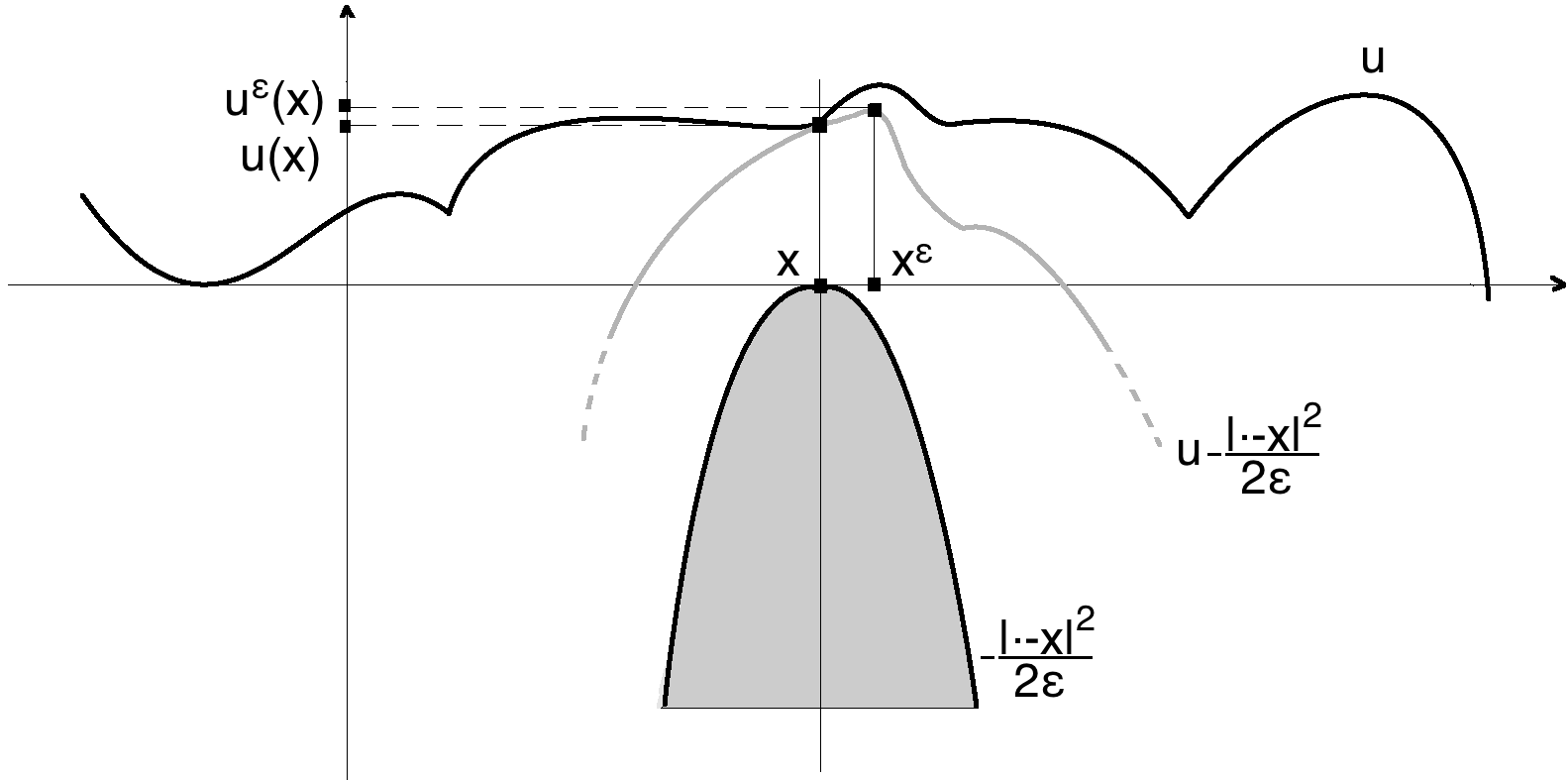}
\]
\textit{The convergence $u^\e \ri u$ that we will establish rigorously later, can be seen geometrically as follows: the factor $1/\e$ of the paraboloid increases its curvature and makes it more and more steep as $\e \ri 0$. This localises the point $x^\e$ (where the maximum is achieved) closer and closer to $x$.}

\ms

Unlike the integral convolution which is in $C^\infty$, the sup-/inf- convolutions are merely \textbf{semiconvex/semiconcave}. Before proceeding to the analytic properties of our approximations, we record the definition of semiconvexity and \textit{Alexandroff's theorem}, which says that semiconvex functions are twice differentiable a.e.\ on its domain.

\ms

\noi \textbf{Definition 3 (Semiconvexity).} \emph{Let $\Om\sub \R^n$ and $f\in C^0(\Om)$. $f$ is called semiconvex when there exists an $\e \in(0,\infty]$ such that the function
\[
x \, \mapsto \, f(x)\, +\, \frac{|x|^2}{2\e}
\]
is convex. If a function is semiconvex for an $\e_0>0$, then it is semiconvex for all $0<\e\leq\e_0$. The inverse of the largest $\e>0$ for which $f$ is semiconvex is called \textbf{semiconvexity constant}: 
\[
\text{Semiconvexity constant}\, :=\, \inf\Big\{\frac{1}{\e} \ \Big| \ x \, \mapsto \, f(x)\, +\, \frac{|x|^2}{2\e} \text{ is convex}\Big\}.
\]
If $1/\e$ is the semiconvexity constant, then $f$ will be called \emph{$\e$-semiconvex}. 
}

\ms

It can be seen that a function $f$ is convex if and only if its semiconvexity constant vanishes. In the obvious way, a function $f$ is called \textbf{semiconcave} if $-f$ is semiconvex.

\ms

\noi \textbf{Remark 4.} \textit{It is well known that a convex function can be characterised by the geometric property that its graph can be touched from below by a hyperplane at every point of its domain of definition and the hyperplane lies below the graph of the function. 
\[
 \includegraphics[scale=0.17]{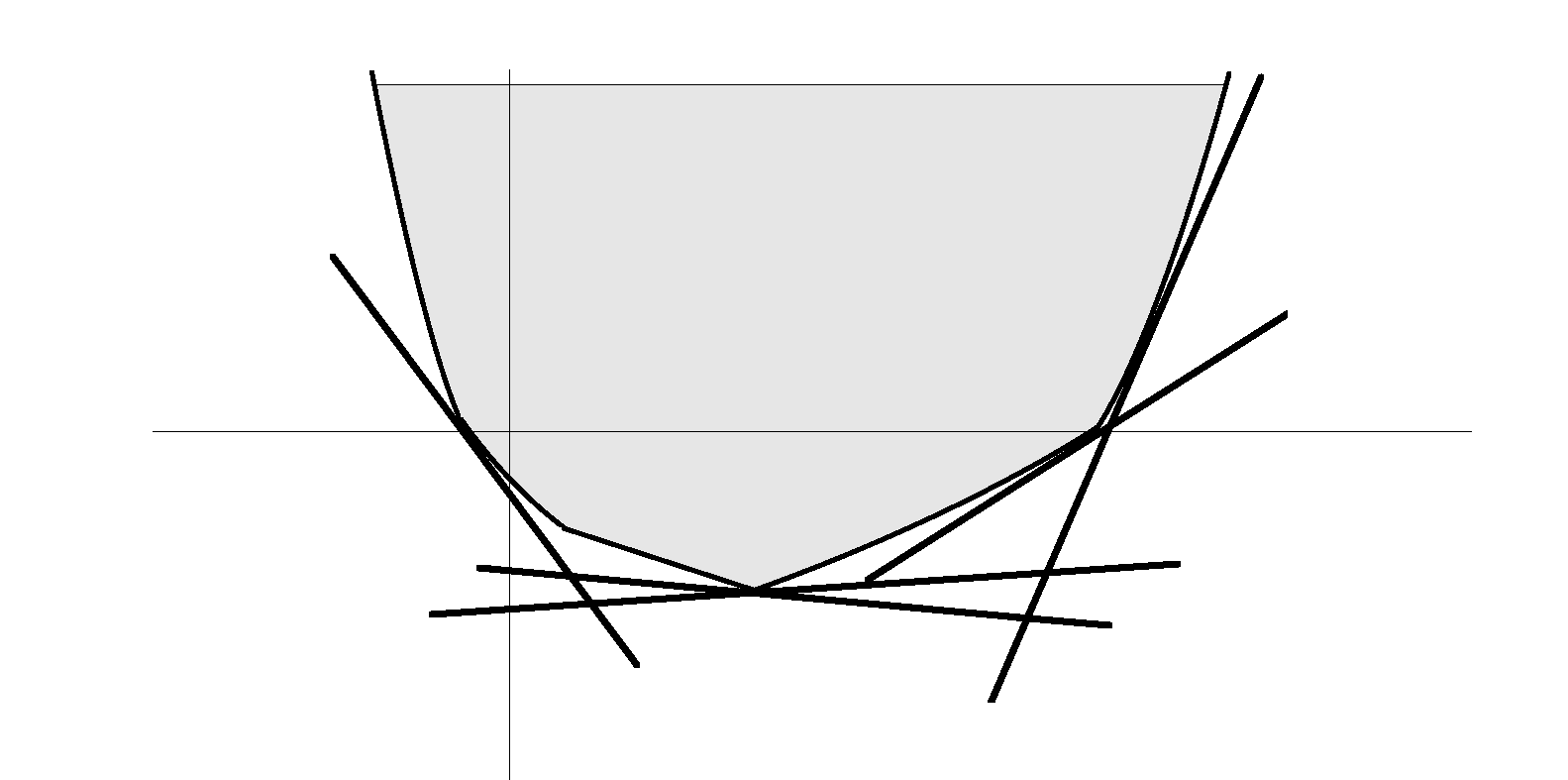} 
\]
It follows that $\e$-semiconvex functions can be geometrically characterised as those whose graph can be touched from below by a paraboloid of curvature $1/\e$  at every point of their domain of definition which lies below the graph of the function.
\[
\includegraphics[scale=0.19]{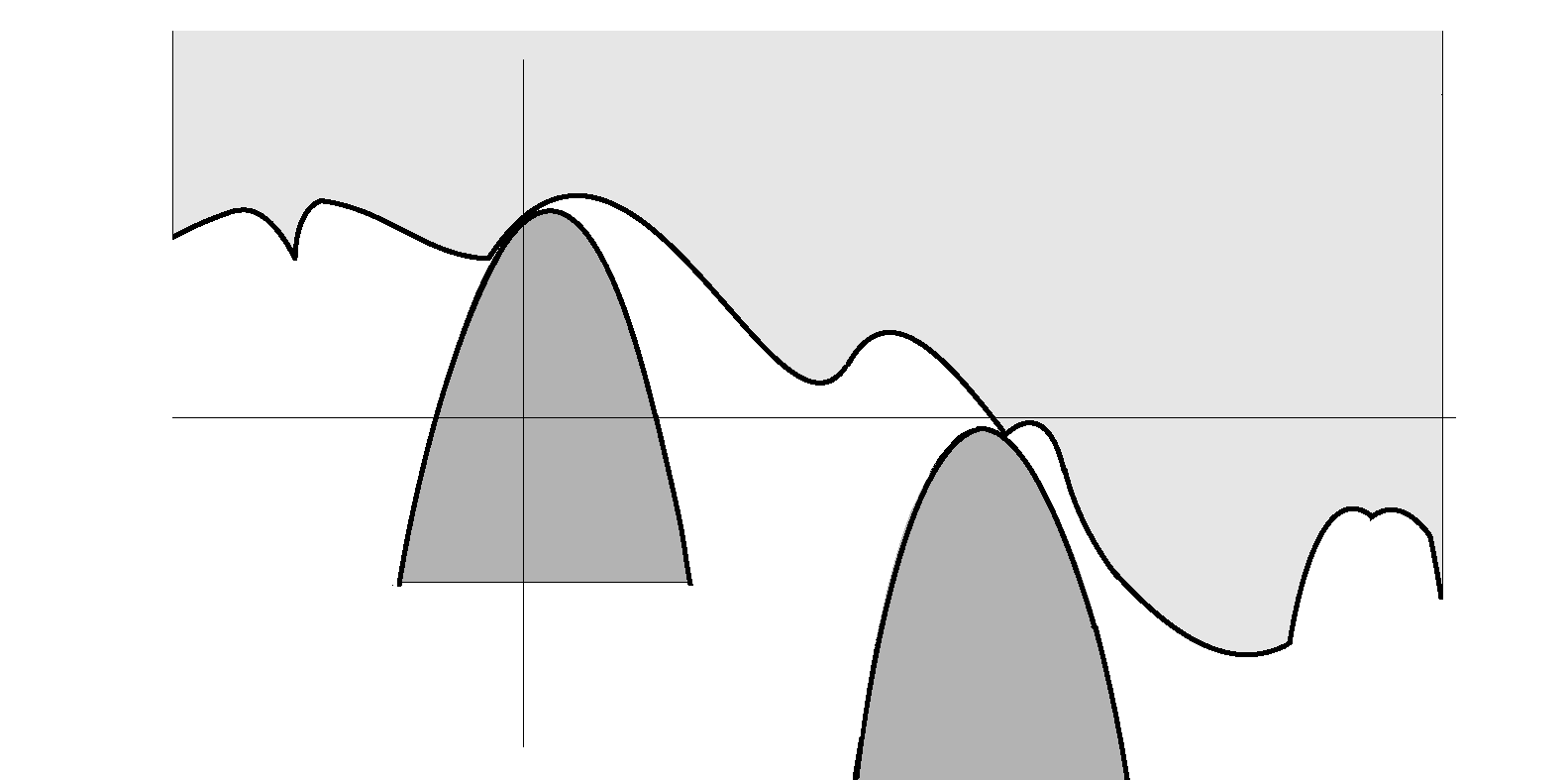}
\]}
\noi \textbf{Remark 5.} 

{\center{
\fbox{\parbox[pos]{260pt}{
\ms

\textsl{Semiconvexity is a sort of ``1-sided $C^{1,1}$ regularity".}

\ms
}}

}}

\ms\ms

\noi \textit{Indeed, a function $u$ belongs to the Lipschitz space $C^{1,1}(\overline{\Om})=W^{2,\infty}(\Om)$ for a Lipschitz domain $\Om$ if and only if there is an $M>0$ such that, for any $x,z\in \Om$, 
\[
\big|u(z+x)\, -\, u(x)\, -\, Du(x)\cdot z \big|\, \leq \, M|z|^2.
\]
This is the same as 
\begin{align}
u(y)\, \leq \, u(x)\, + \, Du(x)\cdot (y-x) \, +\, M|y-x|^2, \nonumber\\
u(y)\, \geq \, u(x)\, + \, Du(x)\cdot  (y-x) \, -\, M|y-x|^2.  \nonumber
\end{align}
This pair of inequalities says that at each $x\in \Om$, $u$ can be touched simultaneously from both above and below by paraboloids with tangents $Du(x)$ at $x$ and curvatures $\pm 2M$. 
\[
\includegraphics[scale=0.19]{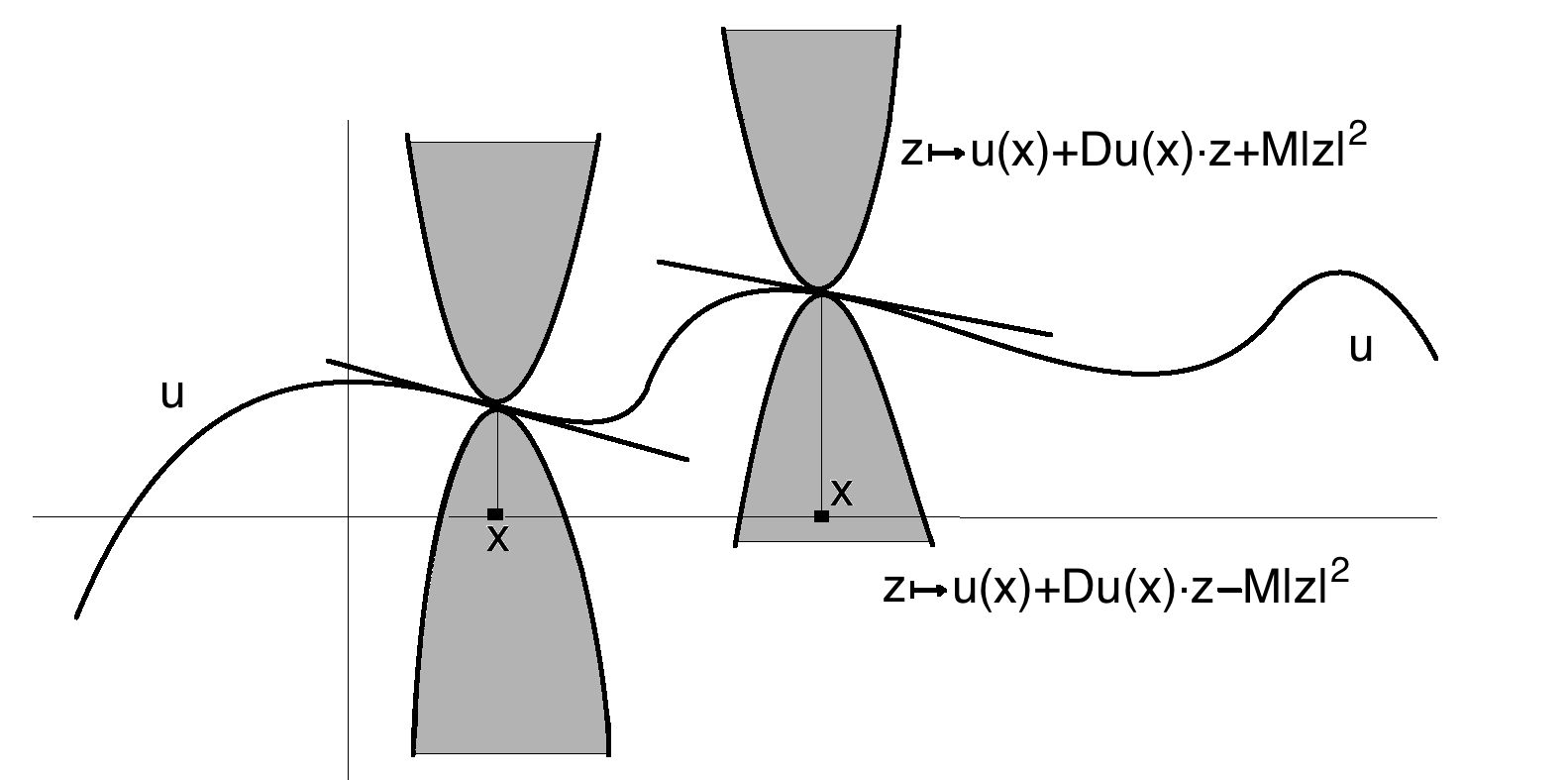}
\]
Equivalently, 
\[
\text{\textsl{$u \in C^{1,1}(\overline{\Om})$ if and only if $u$ is both Semiconvex and Semiconcave.}}
\]}
The basic result concerning the regularity of convex and semiconvex functions is the next

\ms

\noi \textbf{Theorem 6 (Alexandroff).} \emph{Let $f\in C^0(\Om)$ and $\Om\sub \R^n$. If $f$ is semiconvex, then $f$ possesses pointwise second derivatives everywhere on $\Om$ except for a set of Lebesgue measure zero:
\[
\exists \ \big(Df(x),D^2f(x) \big)\ \ \text{for a.e. }x\in \Om.
\]
Moreover, if the semiconvexity constant of $f$ is $1/\e$, then
\[
D^2f(x)  \, \geq \, -\frac{1}{\e}I,\ \ \ \text{for a.e. }x\in \Om.
\]
}
The last matrix inequality is meant in the standard sense, that is that the eigenvalues of  $D^2f(x)$ lie in $[-1/\e,+\infty)$.

\ms

A proof of this important result will be given in the Appendix. The next result collects the main properties of $u^\e$ and $u_\e$.

\ms

\noi \textbf{Theorem 7 (Properties of Sup-Convolution).} Let $u\in C^0(\overline{\Om})$, $\Om\Subset \R^n$.  If $u^\e$ and $u_\e$ are given by \eqref{4.1} and \eqref{4.2} respectively, for $\e>0$ we have:

\ms

\noi (a) $u_\e=-(-u)^\e$.

\ms

\noi (b)  $u^\e\geq u$ on $\Om$.

\ms

\noi (c)  \textbf{(Monotonicity I)}  If $v\in C^0(\overline{\Om})$, then $u\leq v$ implies $u^\e \leq v^\e$.

\ms

\noi (d)  \textbf{(Monotonicity II)} $\e' \leq \e''$ implies $u^{\e'} \leq u^{\e''}$.

\ms

\noi (e) \textbf{(Localisation)}  If we set 
\beq
X(\e)\,  :=\,  \left\{ y \in \overline{\Om}\ \ \Bigg|\ \ u^\e(x)\, =\, u(y)-\frac{|x-y|^2}{2\e}\right\} \label{4.5},
\eeq
then $X(\e)$ is compact and for all $x^\e \in X(\e)$,
\beq \label{4.6}
|x-x^\e|\, \leq \, 2\sqrt{\|u\|_{C^0(\Om)}\e }=:\rho(\e).
\eeq
That is, 
\[
X(\e)\, \sub \, \overline{\mB_{\rho(\e)}(x)}. 
\]
Moreover, 
\beq \label{4.7}
u^\e(x)\, =\, \max_{y\in \overline{\mB_{\rho(\e)}(x)}}\left\{u(y)-\frac{|x-y|^2}{2\e}\right\}\, =\, u(x^\e)-\frac{|x-x^\e|^2}{2\e},
\eeq
for any $x^\e \in X(\e)$. In addition, $\overline{\mB_{\rho(\e)}(x)} \sub \Om$ for $\e$ small.
\ms

\noi (f) \textbf{(Approximation)}   $u^\e \searrow u$ in $C^0(\overline{\Om})$ as $\e \ri 0$.

\ms

\noi (g) \textbf{(Semiconvexity - a.e.\ twice differentiability)}  For each $\e>0$, $u^\e$ is semiconvex and hence twice differentiable a.e. on $\Om$. Moreover, if $(Du^\e, D^2u^\e)$ denote the pointwise derivatives, we have the estimate
\beq \label{4.8}
D^2u^\e \, \geq \, -\frac{1}{\e}  I
\eeq
a.e. on $\Om$.

\ms

\noi (h) \textbf{(``Magic Property")} We have that
\[
(p,X)\in \J^{2,+}u^\e(x) \ \ \Longrightarrow \ \ (p,X)\in \J^{2,+}u(x+\e p),
\]
where
\[
x \,+ \,\e p \,=\, x^\e \ \in \ X(\e).
\]
Moreover, for a.e.\ $x\in \Om$, the set $X(\e)$ is a singleton and 
\[
X(\e) \, =\, x \, +\, \e Du^\e(x).
\]

\ms

\noi \textbf{Proof.} (a) is immediate from the definitions. \ms

(b) By \eqref{4.1}, for all $x,y \in \Om$ we have
\[
u^\e(x)\,  \geq \, u(y) - \frac{|y-x|^2}{2\e}.
\] 
By choosing $y:= x$, we see that $u^\e(x)\geq u(x)$ for all $x\in \Om$.
\ms

(c) If $u\leq v$ and $u,v \in C^0(\overline{\Om})$, for all $x,y \in \Om$ we have
\[
u(y) \, -\, \frac{|y-x|^2}{2\e} \,\leq \, v(y) \, -\, \frac{|y-x|^2}{2\e}.
\]
By taking ``sup" in ${y\in \Om}$, we obtain $u^\e(x) \leq v^\e(x)$.

\ms

(d) If $\e' \leq \e''$, then for all $x,y \in \Om$ we have
\[
u(y) \, -\, \frac{|y-x|^2}{2\e'} \,\leq \, u(y) \, -\, \frac{|y-x|^2}{2\e''}.
\]
By taking ``sup" in ${y\in \Om}$, we obtain $u^{\e'}(x) \leq u^{\e''}(x)$.

\ms

(e) Consider the ``argmax" set $X(\e)$ given by \eqref{4.5}. If $x^\e \in X(\e)$, by \eqref{4.1} we have
\beq  \label{4.19}
u^\e(x) \, =\, u(x^\e)- \frac{|x^\e-x|^2}{2\e}  \, \geq \, u(y) - \frac{|y-x|^2}{2\e} ,
\eeq
for all $y\in \Om$. By choosing $y:=x$ in \eqref{4.19}, we get 
\[
|x^\e-x|^2\, \leq\, 4\|u\|_{C^0(\Om)}\e.
\]
as a result, the set $X(\e)$ is contained in the ball $\mB_{\rho(\e)}(x)$, where 
\[
\rho(\e) \, =\, 2\sqrt{\|u\|_{C^0(\Om)}\e}. 
\]
In addition, since $\Om$ is open, for $\e$ small the closed ball $\overline{\mB_{\rho(\e)}(x)}$ is contained in $\Om$. Hence, \eqref{4.7} follows.

\ms

(f) By assumption, $u$ is uniformly continuous on $\Om$. Hence, by setting
\[
\om (t)\, :=\, \sup \Big\{|u(x)-u(y)| \ : \ |x-y|\leq t,\ x,y \in \Om \Big\},
\]
it follows that there is an increasing $\om \in C^0[0,\infty)$ with $\om(0)=0$ such that 
\[
|u(x)-u(y)| \, \leq \, \om(|x-y|), 
\]
for all $x,y \in \Om$ (modulus of uniform continuity). Hence, for any $x\in \Om$, \eqref{4.7} gives
\begin{align}
u^\e(x) \,&= \, u(x^\e) \, -\frac{|x-x^\e|^2}{2\e}  \nonumber\\
&\leq\, u(x^\e)  \nonumber\\
&\leq \, u(x)\, + \, \om(|x-x^\e|)  \nonumber\\
& \leq \, u(x)\, +\, \om\big(\rho(\e)\big), \nonumber
\end{align}
while by (b) we have $u(x)\leq u^\e(x)$. Hence, $u^\e \ri u $ as $\e \ri 0$ uniformly on $\Om$. By (d), it also follows that $u^\e$ pointwise decreases to $u$ as $\e \ri 0$ .

\ms

(g) By \eqref{4.1}, we have the identity
\[
u^\e(x) \,= \, \sup_{\Om}\left\{ u(y) \, -\frac{1}{2\e}\Big(|x-y|^2-|x|^2 \Big)  \right\}\, -\, \frac{\, |x|^2}{2\e} 
\]
which implies
\[
u^\e(x) \, +\,  -\, \frac{\, |x|^2}{2\e} \, = \, \sup_{\Om}\left\{ \Big(u(y) \, +\frac{\, |y|^2}{2\e}\Big) + \Big(\frac{y}{\e}\Big)\cdot x   \right\}.
\]
It follows that the function
\[
x \, \mapsto \, u^\e(x) \, +\,  \frac{\, |x|^2}{2\e}
\]
is a supremum of the affine functions
\[
x \, \mapsto \,  \Big(u(y) \, +\frac{\, |y|^2}{2\e}\Big) + \Big(\frac{y}{\e}\Big)\cdot x 
\]
and as such, it is convex. By definition, $u^\e$ is semiconvex with semiconvexity constant $1/\e$. The conclusion now follows by invoking Alexandroff's theorem.

\ms

(h) Let $(p,X) \in \J^{2,+}u^\e(x)$. By Theorem 8 of Chapter 2, there is a $\psi \in C^2(\R^n)$ such that 
\[
u^\e-\psi \,\leq \, (u^\e-\psi)(x) 
\]
with $D\psi(x)=p$ and $D^2\psi(x)=X$. Hence, for all $z,y \in \Om$ and $x^\e \in X(\e)$, we have
\begin{align}
u(y)-\frac{|y-z|^2}{2\e}-\psi(z)\, &\leq\, \sup_{y \in \Om}\left\{ u(y)-\frac{|y-z|^2}{2\e}-\psi(z) \right\}  \nonumber\\
&=\,  (u^\e-\psi)(z)  \nonumber\\
 &\leq \, (u^\e-\psi)(x)  \nonumber\\
&=\, u(x^\e)-\frac{|x^\e-x|^2}{2\e}-\psi(x). \nonumber
\end{align}
Hence, we have the inequality
\beq \label{4.25}
u(y)-\frac{|y-z|^2}{2\e}-\psi(z)\, \leq\, u(x^\e)-\frac{|x^\e-x|^2}{2\e}-\psi(x).
\eeq
For $y:=x^\e$ in \eqref{4.25},  we get
\[
|z-x^\e|^2 + 2\e \psi(z)\, \geq\, |x-x^\e|^2 + 2\e \psi(x),
\]
for all $z\in \Om$. This says that the smooth function 
\[
\Psi(z)\, :=\,  |z-x^\e|^2 + 2\e \psi(z)
\]
has a minimum at $x \in \Om$. Consequently, we have that
\[
D\Psi(x)\, = \, 0
\]
which gives that
\[
x^\e-x\,=\, \e p. 
\]
Again by \eqref{4.25} but now for the choice $z:=y-x^\e+x$, we have
\[
u(y)-\psi(y-x^\e+x)\, \leq\, u(x^\e)-\psi(x),
\]
on $\Om$. Hence, by Theorem 8 of Chapter 2 we conclude that 
\[
(p,X) \in \J^{2,+}u\big( x + \e p \big), 
\]
as desired. The last claim follows from Lemma 7 of Chapter 2 and the fact that $u^\e$ is differentiable a.e.\ on $\Om$.                    \qed

\ms

\noi \textbf{Remark 8.} \textit{The ``magic property" of sup-convolutions is a sort of converse to the stability theorems for Jets of Capter 3. The latter statement is meant in the sense that a jet of the approximating sequence $u^\e$ is passed to the limit $u$, while, in the stability theorems, given a jet of $u$ we construct a jets of the approximating sequence. The ``magic property" will be needed latter in the uniqueness theory for the Dirichlet problem.
}

\ms

Now we consider the important ``PDE properties" of sup-/inf- convolutions. We will consider only the case of ``sup-convolution/subsolution". The symmetric case of ``inf-convolution/supersolution" is left as a simple exercise for the reader. For the statement we will the \textit{inner open neighbourhood of $\Om$} of width $\rho(\e)$:
\beq \label{4.a}
\Om^{\rho(\e)} \, :=\, \Big\{x\in \Om\ \big|\ \dist(x,\p \Om)>\rho(\e) \Big\}.
\eeq
\[
\includegraphics[scale=0.21]{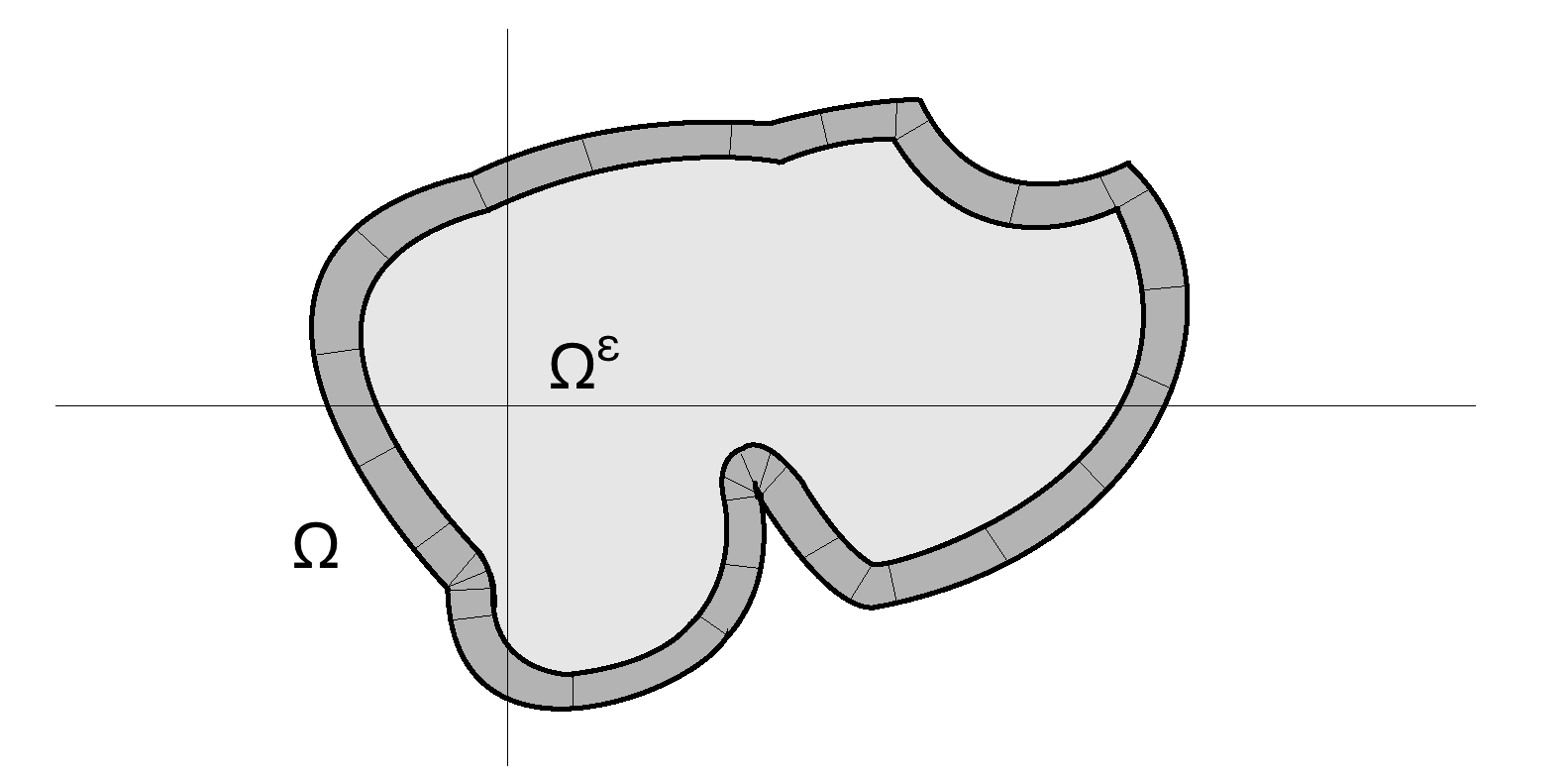}
\]
We will also need to restrict our nonlinearity $F$ defining the PDE by assuming a non-strict monotonicity property in the $u$ argument. This property, will appear again later when we will consider uniquness results for the Dirichlet problem.

\ms

\noi \textbf{Definition 9 (Proper nonlinearities).} \emph{Let $F \in C^0\big( \Om \by \R \by\R^n \by \mS(n) \big)$, $\Om \sub \R^n$. The function $F$ is called \textbf{proper} when the following monotonicity properties hold:
\begin{align}
X\leq Y\ \text{ in }\mS(n)\ \ &\Longrightarrow \ \ F(x,r,p,X)\leq F(x,r,p,Y), \nonumber\\
r\leq s\ \ \text{ in }\R\ \ \ \ \ \, &\Longrightarrow \ \ F(x,r,p,X)\geq F(x,s,p,X).  \nonumber
\end{align}
}
\ms

Hence, $F$ is proper when the respective PDE $F=0$ is degenerate elliptic and the partial function $r\mapsto F(x,r,p,X)$ is monotone non-increasing. One final comment is that we state and prove the next result only for PDE of the \textit{``decoupled in the $x$-dependence'} form
\[
F(u,Du,D^2u)\, \geq\, f.
\]
The result however holds more generally for nonlinearities of the form $F = F(x,r,p,X)$, but under extra assumption and requires a less elementary proof. Now we are ready for the

\ms

\noi \textbf{Theorem 10 (Sup-convolution as Viscosity and a.e.\ Subsolution).} 
\emph{Let $F \in C^0\big(\R \by\R^n \by \mS(n) \big)$ be a proper function,  $f\in C^0(\Om)$, $\Om \sub \R^n$. Suppose that $u$ is a Viscosity Solution of 
\[
F(u,Du,D^2u)\, \geq\, f
\]
on $\Om$. Then for every $\e>0$ small, its sup-convolution approximation $u^\e$ is a Viscosity Solution of 
\[
F(u^\e,Du^\e,D^2u^\e)\, \geq \, f_\e
\]
on $\Om^{\rho(\e)}$, where
\[
f_\e(x) \, :=\, \inf_{\mB_{\rho(\e)}(x)}f
\]
and $f_\e \ri f$ in $C^0(\Om)$, as $\e \ri 0$.
 }

\emph{Moreover, $u^\e$ is a strong solution, that is it is twice differentiable a.e.\ on $\Om^{\rho(\e)}$ and satisfies the PDE classically on $\Om^{\rho(\e)}$ except for a set of Lebesgue measure zero.}

\ms

\noi \textbf{Proof.} Fix $\e>0$ small such that $\Om^{\rho(\e)}$ of \eqref{4.a} is an non-empty open inner neighbourhood of $\Om$. We begin by observing that by \eqref{4.1}, \eqref{4.7}, we may write
\[
u^\e(x)\, =\, \sup_{|z|<\rho(\e)}\left\{  u(x+z)\, -\, \frac{|z|^2}{2\e} \right\}
\]
For $|z|< \rho(\e)$ and $\e>0$ both fixed, we define the function
\beq   \label{4.24}
u^{\e,z}(x)\, :=\, u(x+z)\, -\, \frac{|z|^2}{2\e}.
\eeq
The inequality $|z|< \rho(\e)$ implies that $x+z \in \Om$ and hence each $u^{\e,z}$ is well-defined and continuous on $\Om^{\rho(\e)}$. Now it is a very simple fact to verify (directly by the definition of Jets) that
\[
\J^{2,+}u^{\e,z}(x)\, =\, \J^{2,+}u(x+z). 
\]
Fix $(p,X)\in \J^{2,+}u^{\e,z}(x)$. Then, by using the properness of $F$, we have
\begin{align}
F\big(u^{\e,z}(x),p ,X\big)\, & = \, F\left(u(x+z)-\frac{|z|^2}{2\e},p ,X\right)  \nonumber\\
&\geq\, F\big(u(x+z),p ,X\big).  \nonumber
\end{align}
Since $(p,X)\in \J^{2,+}u(x+z)$ and $u$ is a viscosity solution of 
\[
F(u,Du,D^2u)\, -\, f\, \geq\, 0
\]
on $\Om$, it follows that
\begin{align}
F\big(u^{\e,z}(x),p ,X\big)\, &\geq\, F\big(u(x+z),p ,X\big)  \nonumber\\
 &\geq\, f(x+z)  \nonumber\\
&\geq\,  \inf_{|z|<\rho(\e)}\, f  \nonumber\\
&=\, f_\e(x). \nonumber
\end{align}
Hence, for $|z|< \rho(\e)$ and $\e>0$, each $u^{\e,z}$ is a Viscosity Solution of 
\beq \label{4.14}
F\big(u^{\e,z},Du^{\e,z},D^2u^{\e,z}\big)\, - \,f_\e\, \geq\, 0
\eeq
on $\Om^{\rho(\e)}$. In order to conclude, we invoke \eqref{4.24} and Theorem 12 of Chapter 3: since
\[
u^\e(x)\, =\, \sup_{|z|< \rho(\e)} \, u^{\e,z}(x)
\]
and $u^\e$ is finite and continuous, it follows that $u^\e$ is the supremum of subsolutions which satisfies the assumptions of the result just quoted. Hence, $u^\e$  is a subsolution of \eqref{4.14} itself. By (g) of Theorem 7, it follows that $u^\e$ solves classically the inequality a.e.\ on $\Om^{\rho(\e)}$.         \qed

\ms

\noi \textbf{Remark 11.} \textit{The ``magic" regularisation properties of the sup-/inf- convolution approximations \eqref{4.1} and \eqref{4.2} can be seen as a result of the \textbf{Lax formulas} which represent the (unique viscosity) solution of the initial-value problems for the evolutionary Hamilton-Jacobi PDE
\[
\left\{
\begin{array}{l}
u_t(x,t)\, =\, \dfrac{1}{2}\big|Du(x,t)\big|^2,\ \ x\in \R^n, \ t>0 \ms\\
\, u(x,0)\, =\, u(x), \ \ \ \ \ \ \ \ \ \ \ \ x\in \R^n,
\end{array}
\right.
\]
and
\[
\left\{
\begin{array}{l}
u_t(x,t)\, =\, -\dfrac{1}{2}\big|Du(x,t)\big|^2,\ \ x\in \R^n, \ t>0 \ms\\
\, u(x,0)\, =\, u(x), \ \ \ \ \ \ \ \ \ \ \ \ \ x\in \R^n.
\end{array}
\right.
\]
Indeed, the solutions can be represented by the variational formulas
\[
u(x,t)\, =\, \sup_{y\in \Om}\left\{u(y)\, -\, \frac{|x-y|^2}{2t}\right\}, \ \ x\in \R^n,
\]
and
\[
u(x,t)\, =\, \inf_{y\in \Om}\left\{u(y)\, +\, \frac{|x-y|^2}{2t}\right\}, \ \ x\in \R^n,
\]
respectively. Roughly, we mollify Viscosity Solutions by using the flow map of the Hamilton-Jacobi PDE.
}

\ms

\ms

\ms

\noi \textbf{Remarks on Chapter 4.} An excellent introduction to semiconvexity with applications to Optimal Control (a topic not considered in these notes) is Cannarsa-Sinestrari \cite{CS}. Background material on the simple measure-theoretic notions needed in this chapter as well as on the standard notion of  mollification via integral convolution (for divergence-structure PDE) can be found in the appendix of the textbook of Evans \cite{E4}. Sup/Inf convolutions have been introduced in the context of Viscosity Solutions (in different guises) by Jensen in \cite{J1}, as a regularisation tool for the proof of uniqueness to the Dirichlet problem (which is the subject of Chapter 6 in these notes). A fairly simple proof of Alexandroff's theorem which assumes a certain knowledge of geometric measure theory can by found in Evans-Gariepy \cite{EG}. A shorter proof which uses less measure theory (once the reader is willing to accept the validity of the so-called ``co-area formula") can be found in the appendix of Crandall-Ishii-Lions \cite{CIL}. The proof of Theorem 10 is much more difficult when one assumes dependence on $x$. A proof which requires usage of the comparison principle for viscosity solutions can be found e.g.\ in \cite{I3}. The part of this chapter on the properties of sup-convolutions follows \cite{CIL, C2} closely, but herein we give many more details. The term ``magic property" we used in $(h)$ of Theorem 7 follows the by now standard terminology used in \cite{CIL}.

\chapter[Existence via Perron's Method]{Existence of Solution to the Dirichlet Problem via Perron's Method}

In this chapter we are concerned with the question of existence of a Viscosity Solution $u \in C^0(\overline{\Om})$ to the Dirichlet problem 
\beq \label{5.1}
\left\{
\begin{array}{l}
F\big(\cdot,u,Du,D^2u\big)\,=\,0, \ \text{ in }\Om, \ms\\
\hspace{78pt}u\, =\, b,\, \text{ on }\p \Om,
\end{array}
\right.
\eeq
where $\Om \Subset \R^n$, $b \in C^0(\overline{\Om})$ is the ``boundary condition", and 
\[
F \ :\ \Om \by \R \by \R^n \by \mS(n) \larrow \R
\]
is continuous and degenerate elliptic. As always, the later means that the partial function $X\mapsto F(x,r,p,X)$ is monotone, that is
\[
(X-Y)\Big( F(x,r,p,X)\, -\, F(x,r,p,Y)\Big)\, \geq \, 0,
\]
for all $X,Y \in \mS(n)$ and $(x,r,p)\in \Om \by \R \by \R^n$. It is quite obvious that by \textbf{solution to the Dirichlet problem \eqref{5.1}} above, we mean

{\center{
\fbox{\parbox[pos]{270pt}{
\ms

\textit{a Viscosity Solution to the PDE on $\Om$, which satisfies} 

\, \textit{the boundary condition in the classical sense on $\p \Om$}. 

\ms
}}

}}

\ms

More general nonlinear boundary conditions (that will require weak ``viscosity interpretation") will be considered in later chapters. Herein we are concerned with the adaptation of \emph{Perron's method of sub-/super- solutions} to the Viscosity setting. Roughly, 

{\center{
\fbox{\parbox[pos]{230pt}{
\ms

PERRON'S METHOD GIVES EXISTENCE

\ \ \ \ \ BY ASSUMING UNIQUENESS ! ! !

\ms
}}

}}

\ms

This method has been introduced by Oscar Perron in 1923 for the Laplacian. In the next chapter will be concerned with the \emph{problem of uniqueness of Viscosity Solutions} to \eqref{5.1} and with sufficient conditions that guarantee it. Let us record here the definition of Viscosity Sub-/super- Solutions, under the slight modification imposed by Remark 9 of Chapter 3, in which we split the continuity requirement of the solution to halves as well:

\ms

\noi \textbf{Definition 1\footnote{We re-state the definition only for convenience of the reader. The \emph{only difference} from the previous version we have introduced is that now we require the sub/super solutions to be only $USC$/$LSC$ respectively and not continuous.} (Semicontinuous Viscosity Sub-/Super- Solutions).} 

\noi (a) \emph{The function $u\in USC(\Om)$ is a Viscosity Solution of 
\[
F\big(\cdot,u,Du,D^2u \big)\, \geq \,0
\]
on $\Om$, when for all $x\in \Om$,
\[
 \inf_{(p,X)\in \J^{2,+}u(x)}\, F\big(x,u(x),p,X\big)\geq 0.
\]
\noi (b) The function $u\in LSC(\Om)$ is a Viscosity Solution of 
\[
F\big(\cdot,u,Du,D^2u \big)\, \leq \,0
\]
on $\Om$, when for all $x\in \Om$,
\[ 
 \sup_{(p,X)\in \J^{2,-}u(x)}\, F\big(x,u(x),p,X\big)\leq 0.
\]
\noi (c) The function $u\in C^0(\Om)$ is a Viscosity Solution of 
\[
F\big(\cdot,u,Du,D^2u \big)\, =\,0
\]
when it is both a Viscosity Subsolution and a Viscosity Supersolution.}

\ms

We recall that the definition of semicontinuity is given in Chapter 3 and that the Sub-/Super- Jets $\J^{2,\pm}$ are given by \eqref{2.6a}, \eqref{2.7a}. Let us also record the assumption we will need for our existence result, in the form of definition:

\ms

\noi \textbf{Definition 2 (Comparison for the Dirichlet Problem \eqref{5.1}).} 

\noi \emph{(a) A function $\underline{u} \in USC(\overline{\Om})$ is called a Viscosity Subsolution of \eqref{5.1} when it is a Viscosity Subsolution of
\[
F\big(\cdot,\underline{u},D\underline{u},D^2\underline{u} \big)\, \geq\,0
\]
on $\Om$ and satisfies $u\geq b$ on $\p \Om$.}

\ms

\noi \emph{(b) A function $\overline{u} \in LSC(\overline{\Om})$ is called a Viscosity Supersolution of \eqref{5.1} when it is a Viscosity Supersolution of
\[
F\big(\cdot,\overline{u} ,D\overline{u} ,D^2\overline{u}  \big)\, \leq\,0
\]
on $\Om$ and satisifes $u\leq b$ on $\p \Om$.}

\ms

\noi \emph{(c) We say that \eqref{5.1} satisfies the \textbf{Comparison Principle}, when, for any Subsolution $\underline{u} \in USC(\Om)$  and any Supersolution $\overline{u} \in LSC(\Om)$, we have
\[
\underline{u}\, \leq \, \overline{u},\ \text{ on } \Om.
\]
}
It is a quite evident consequence of the above definition that 

{\center{
\fbox{\parbox[pos]{290pt}{
\ms

Every solution $u$ of the Dirichlet problem is ``sandwiched" 

\ \ \ \ \ between a subsolution $\underline{u}$ and a supersolution $\overline{u}$:
\[
\underline{u}\, \leq \, u \, \leq \, \overline{u},\ \text{ on } \Om.
\]

}}

}}

\ms

\noi The following is the principal result of this chapter.

\ms

\noi \textbf{Theorem 3 (Existence via Perron's Method).} \emph{Let $F \in C^0\big( \Om \by \R \by\R^n \by \mS(n) \big)$ be proper with $\Om \Subset \R^n$. Consider the Dirichlet problem
\[
\left\{
\begin{array}{l}
F\big(\cdot,u,Du,D^2u\big)\,=\,0, \ \text{ in }\Om, \ms\\
\hspace{78pt}u\, =\, b,\, \text{ on }\p \Om,
\end{array}
\right.
\]
for $b\in C^0(\overline{\Om})$, and assume that it satisfies the Comparison Principle (Definition 2). If there exist a Subsolution $\underline{u}$ and a Supersolution $\overline{u}$ of  the Dirichlet problem, then
\beq \label{5.2}
V(x)\, :=\, \sup \Big\{v(x)\ :\  \underline{u}\, \leq \, v \, \leq \, \overline{u} \ \text{ and $v$ subsolution}\Big\}
\eeq
defines a Viscosity Solution of the Dirichlet problem.
}

\ms

For the proof we will need the simple notion of \emph{semicontinuous envelopes}. 

\ms

\noi \textbf{Definition 4 (Semicontinuous envelopes).} \emph{Let $u : \Om \sub \R^N \ri \R$. We set
\begin{align}
u^*(x)\, &:=\, \lim_{r \ri 0} \, \sup_{\mB_r(x)}u,\ms \label{5.3}\\
u_*(x)\, &:=\, \lim_{r \ri 0} \, \inf_{\mB_r(x)}u. \label{5.4}
\end{align}
We call $u^* : \Om \ri (-\infty,+\infty]$ the upper semicontinuous envelope of $u$ and $u_* : \Om \ri [-\infty,+\infty)$ the lower semicontinuous envelope of $u$.}

\ms

\noi $u^*$ can be characterised as \emph{the smallest function $v\in USC(\Om)$ for which $v\geq u$}. Symmetrically, $u_*$ can be characterised as \emph{the largest function $v\in LSC(\Om)$ for which $v\leq u$}.

\ms

\noi \textbf{Example 5.} The function
\[
u(x)\, :=\, 
\left\{
\begin{array}{l}
-1, \ \ x<0\\
\ \ 0,\ \ x=0\\
+1, \ \ x>0,
\end{array}
\right.
\]
has upper and lower semicontinuous envelopes as shown in the figure.
\[
\includegraphics[scale=0.2]{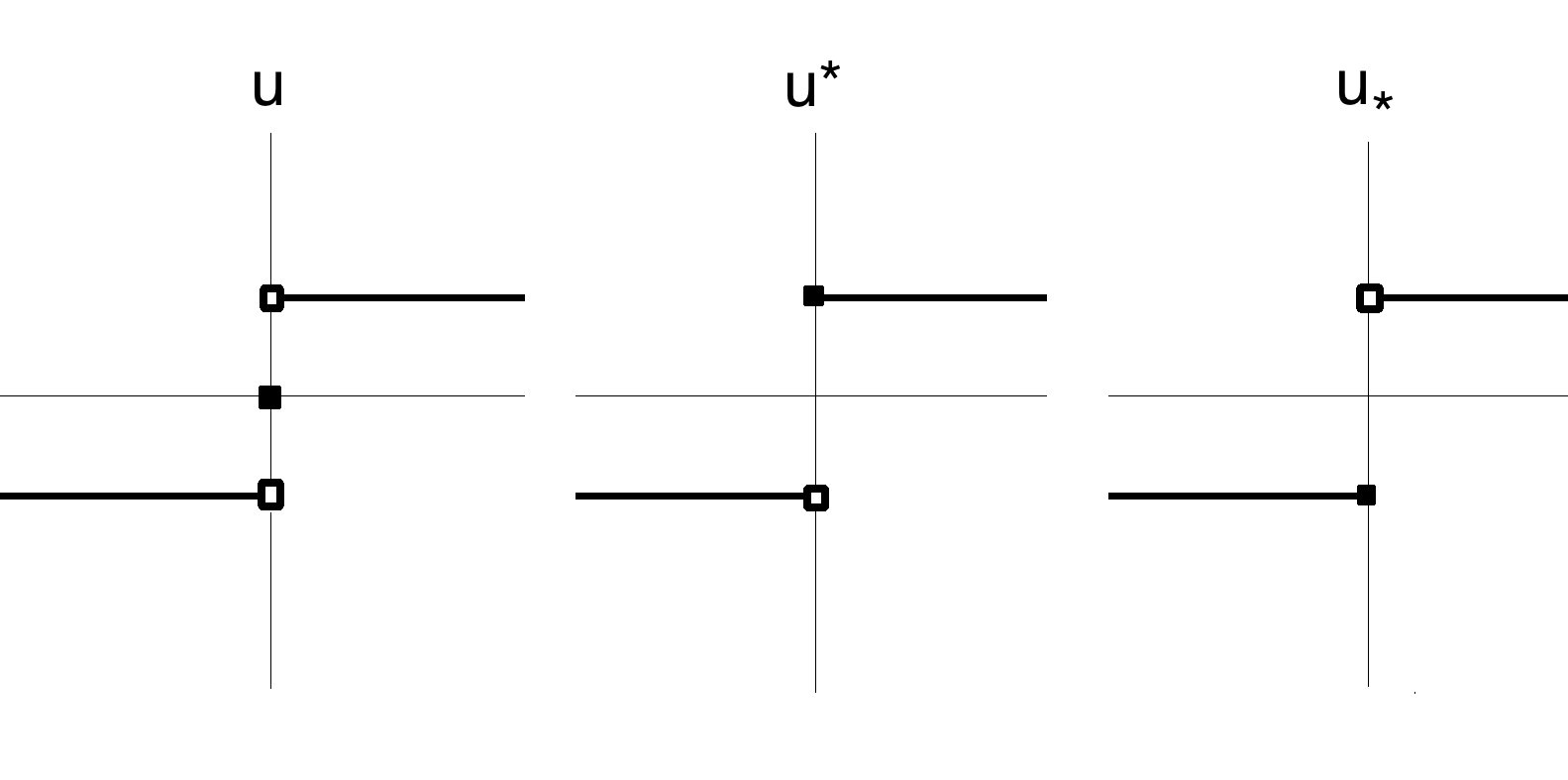}
\]

The first ingredient needed in the proof of Theorem 4 is the following minor variant of Theorem 12 of Chapter 3:

\ms

\noi \textbf{Theorem 6 (USC envelope of suprema of subsolutions is subsolution).} \emph{Let $\Om\sub \R^N$ and let also $F \in C^0\big(\Om \by \R \by \R^n \by \S(n) \big)$. Suppose that $\mU \sub USC(\Om)$ is a family of Viscosity Subsolutions, i.e.\ for each $u \in \mU$,
\[
F\big(\cdot,u,Du,D^2u\big)\,\geq\, 0,
\]
on $\Om$. We set
\[
U(x) \, :=\, \sup_{u \in \mU}\,u(x),\ \ x\in \Om. 
\]
If the upper semicontinuous envelope $U^*$ is finite on $\Om$ (i.e. $U^*(x)<\infty$ for all $x\in \Om$), then $U^*$ is  a Viscosity Subsolution on $\Om$ as well.
}

\ms

We refrain from giving the proof of Theorem 6, which is almost identical to that of Theorem 12 of Chapter 3. The next ingredient is the next lemma, which roughly says

{\center{
\fbox{\parbox[pos]{330pt}{
\ms

If a subsolution fails to be a supersolution, then it is not maximal.

}}

}}

\ms

\noi More precisely, we have

\ms

\noi \textbf{Lemma 7 (The ``Bump construction").} \emph{Let $F \in C^0\big(\Om \by \R \by \R^n \by \S(n) \big)$, where  $\Om\sub \R^N$. Suppose $u \in USC(\Om)$ is a  Viscosity Subsolution of
\[
F\big(\cdot,u,Du,D^2u\big)\,\geq\, 0,
\]
on $\Om$. Assume that there exists $\hat{x}\in \Om$ such that $u_*$ fails to be a Viscosity Supersolution at $\hat{x}$ \ms
}

( i.e. there is $(p,X)\in \J^{2,-}u_*(\hat{x})$ such that $F(\hat{x},u(\hat{x}),p,X)>0$ ). \ms

\noi \emph{Then, for each $r>0$ small enough, there exists a Viscosity Subsolution $U\in USC(\Om)$ strictly greater than $u$ near $\hat{x}$, i.e.
\begin{align}
U \, \geq \, u \, \text{ on }\, \mB_r(\hat{x}),\ \ U \, \equiv \, u \, \text{ on }\, \Om \set \mB_r(\hat{x}),  \nonumber\\
U \, > \, u \, \text{ on a nonempty subset of }\ \mB_r(\hat{x}). \ \, \nonumber
\end{align}
}

\ms

\noi \textbf{Proof.} We simplify things by assuming $\hat{x}=0\in \Om$. This can be easily achieved by a translation. Since $u_*$ is not a subsolution at $0$, there is $(p,X)\in \J^{2,-}u_*(0)$ such that
\[
F(0,u_*(0),p,X)\, > \, 0.
\]
We fix $\ga,\de,r>0$ and define the function
\[
w_{\ga,\de}(z)\, :=\, \de \, +\, u_*(0)\, +\,  p\cdot z\, +\, \frac{1}{2}X:z\ot z\, -\, \frac{\ga}{2}|z|^2.
\]
By continuity, $w_{\ga,\de}$ is classical subsolution on the ball $\mB_r(0)$, if $\ga,\de,r$ are small enough: indeed, we have
\begin{align}
F\big(\cdot,w_{\ga,\de},D w_{\ga,\de}, D^2 w_{\ga,\de}  \big)\, &=\, F\Big(\cdot,\de + u_*(0) +O(|r|), p+O(|r|), X-\ga I \Big) \nonumber\\
&\geq\, F(0,u_*(0),p,X)\, +\, o(1) \ \ \ \ \ \ \Big( \text{ as }\ga,\de,r \ri 0\ \Big)\nonumber\\
&\geq \, 0, \nonumber
\end{align}
on  $\mB_r(0)$. We fix such a triple $(\ga,\de,r)$. Since $(p,X)\in \J^{2,-}u_*(0)$, we have
\[
u(z)\, \geq\, u_*(z)\, \geq\, u_*(0)\, +\, p\cdot z\, +\, \frac{1}{2}X:z\ot z\, +\, o(|z|^2),
\]
as $z\ri 0$. By keeping $\ga$ constant and decreasing $r,\de$ further, we can achieve $u>w_{\ga,\de}$ outside the half ball $\mB_{r/2}(0)$: indeed, for $r/2<|z|<r$, we have
\begin{align}
w_{\ga,\de}(z)\,& =\, \de \, +\, u_*(0)\, +\,  p\cdot z\, +\, \frac{1}{2}X:z\ot z\, -\, \frac{\ga}{2}|z|^2 \nonumber\\
&\leq \, \de\, +\, u_*(z)\, +\, o(|z|^2) -\, \frac{\ga}{2}|z|^2  \ \ \ \ \ \Big( \text{ as }z \ri 0\ \Big) \nonumber\\
&\leq \,  u_*(z)\, +\, \de\, +\, \Big(o(1) -\, \frac{\ga}{4}\Big)r^2 \ \ \ \ \ \ \Big( \text{ as }r \ri 0\ \Big)\nonumber\\
&\leq \,  u_*(z) \nonumber\\
&\leq \,  u(z). \nonumber
\end{align}
We set 
\[
U\, :=\,
\left\{
\begin{array}{l}
\max\big\{u, w_{\ga,\de}\big\}, \ \ \ \text{ on }\mB_r(0),\\ 
u, \hspace{75pt} \text{ on }\Om \set \mB_r(0).
\end{array}
\right.
\] 
Since both $u$ and $w_{\ga,\de}$ are upper-semintinuous subsolutions, by invoking Theorem 6 we conclude that $U$ is a subsolution as well. We also have that
\begin{align}
\underset{z\ri 0}{\lim\sup}\, \big(U\, - \,u\big)(z)\, &=\, U(0)\, -\, \underset{z\ri 0}{\lim\inf}\, u(z)\nonumber\\
&=\, U(0)\, -\, u_*(0)  \nonumber\\
& = \, \de \, +\, u_*(0)\, -\, u_*(0) \nonumber\\
&>\, 0.  \nonumber
\end{align}
\[
\includegraphics[scale=0.22]{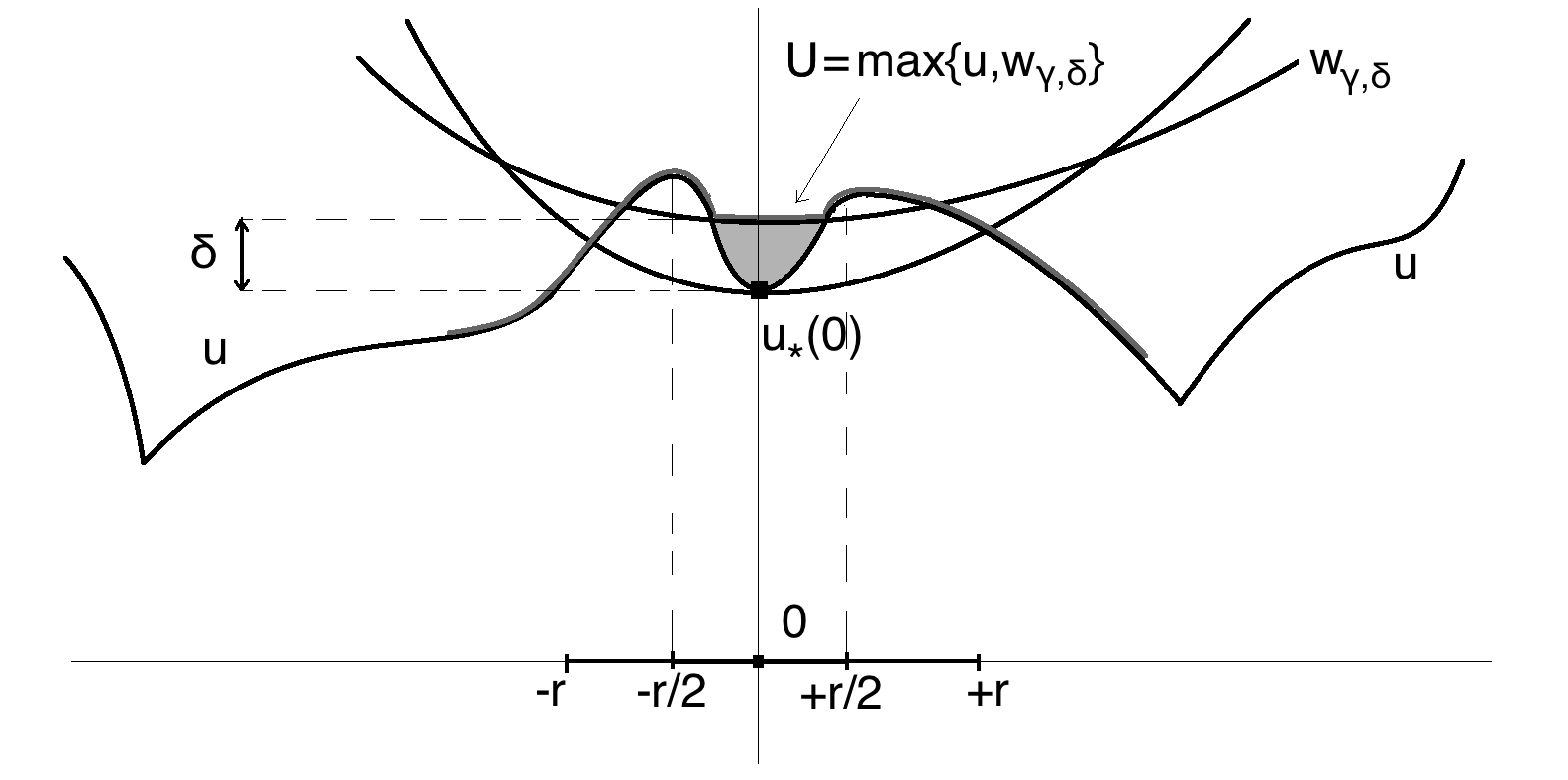}
\]
Hence, we conclude that $U$ is a subsolution which satisfies the claimed properties. The lemma ensues.                     \qed 

\ms

We may now prove Theorem 3.

\ms

\noi \textbf{Proof.} Let $V$ be given by \eqref{5.2} and note that $\underline{u} \leq V \leq \overline{u}$. By the definitions \eqref{5.3},  \eqref{5.4}, it follows that the envelopes respect inequalities. Hence,
\[
\underline{u}_* \, \leq \, V_*\, \leq \, V\, \leq \, V^*\, \leq \, \overline{u}^*,
\]
on $\Om$. Since $\underline{u}$ and $\overline{u}$ are sub-/super- solutions of the Dirichlet problem, we have 
\[
V_*\, =\, V\, = \, V^*\, =\, b 
\]
on $\p \Om$. By Theorem 7, we have that $V^*$ is a subsolution. Since (by assumption) the Comparison Principle holds and since $\underline{u}$ is in $USC(\Om)$, we have 
\[
\underline{u}\, = \, \underline{u}^* \, \leq \, V^* \, \leq \, \overline{u}. 
\]
By \eqref{5.2}, we see that $V$ is maximal and hence $V\geq V^*$. On the other hand, by definition of the envelope, we have $V\leq V^*$. Consequently, $V=V^* \in USC(\Om)$ and hence $V$ is a subsolution. We claim that $V$ is a supersolution as well. If the supersolution property fails at some point $\hat{x}\in \Om$, then by the ``bump construction" of Lemma 7 we can find $U\geq V$ with $U>V$ near $\hat{x}$ and $U\equiv V $ on $\overline{\Om} \set \mB_r(\hat{x})$. By the Comparison Principle, we have $U \leq \overline{u}$. This leads to a contradiction to the maximality of $V$, and hence we conclude that $V$ is a (continuous) Viscosity Solution of the Dirichlet problem.        \qed

\ms

\noi \textbf{Remark 8.} \textit{Theorem 3 provides the existence of solution for \eqref{5.1}, provided that the PDE satisfies the Comparison Principle and provided that a subsolution and a supersolution of  \eqref{5.1} have been constructed.}

\ms

Sufficient conditions guaranteeing the validity of the Comparison Principle will be given in the next chapter. However, the results in this chapter leave open the question of how to find a pair of sub-/super- solutions. The next result gives a clue of how to do that in an interesting case.

\ms

\noi \textbf{Application 9 (Construction of Sub-/Super- Solutions for degenerate elliptic PDE decoupled in $u$).} 

\noi \emph{Let $F \in C^0\big(\Om \by \R^n \by \mS(n)\big)$, $f\in C^0(\R)$ and $\Om \Subset \R^n$ with $\p \Om$ $C^2$ smooth. We assume that $X\mapsto F(x,p,X)$ is monotone and satisfies $F(x,0,0)\geq0$, and also that $f(0)=0$ and that $f$ is strictly increasing. Then, the Dirichlet problem
\beq \label{5.5}
\left\{
\begin{array}{l}
F\big(\cdot,Du,D^2u\big)\,=\,f(u), \ \text{ in }\Om, \ms\\
\hspace{67pt}u\, =\, 0,\ \ \, \text{ on }\p \Om,
\end{array}
\right.
\eeq
has a pair $(\underline{u},\overline{u})$ of sub-/super- solutions under a technical assumption given later.}

\ms

We begin by simplifying \eqref{5.5} to the equivalent problem
\beq \label{5.6}
\left\{
\begin{array}{l}
G\big(\cdot,Du,D^2u\big)\,=\,u, \ \text{ in }\Om, \ms\\
\hspace{67pt}u\, =\, 0,\ \ \, \text{ on }\p \Om.
\end{array}
\right.
\eeq
This is possible since $f$ is strictly increasing. Indeed, by setting $G:=f^{-1}\circ F$, we have that $X\mapsto G(x,p,X)$ is monotone as well and hence the PDE
\[
G\big(\cdot,Du,D^2u\big)\,-\,u\, =\, 0
\]
is degenerate elliptic. Since $F(x,0,0)\geq0$ and $f^{-1}$ is increasing and $f^{-1}(0)=0$, it readily follows that 
\[
G(x,0,0)\, =\, f^{-1}\big(F(x,0,0)\big)\, \geq\, 0. 
\]
Thus,
\[
\underline{u}\, \equiv\, 0
\]
is a subsolution of \eqref{5.6}. Now we work towards the construction of a supersolution. We fix $M,L>0$ and set
\[
u^+(x)\, :=\, M\left(1\, -\, e^{-L\, d(x)} \right).
\]
Here $d$ denotes the distance function from the boundary:
\[
d(x) \, \equiv\, \dist(x,\p \Om)\, :=\, \inf_{y\in \p \Om}\, |y-x|.
\]
For $\e>0$ small, consider the inner neighbourhood of the boundary of width $\e$:
\[
\Om_\e\, :=\, \big\{x\in \Om\ :\ d(x) \, <\, \e \big\}.
\]
Since $\p \Om$ is $C^2$ smooth, if $\e$ is small enough, it is a well known consequence of the implicit function theorem that $d \in C^2(\overline{\Om_\e})$. We have
\begin{align}
G\big(\cdot,Du^+,&D^2u^+\big)\,-\,u^+ \nonumber\\ 
&\leq\,  G\Big(\cdot,MLe^{-L\, d}Dd,MLe^{-L\, d}\big(D^2d-L Dd \ot Dd\big) \Big). \nonumber
\end{align}
For  $x\in \Om_{1/L}$, we have $d(x)\leq 1/L$ and hence
\[
\frac{1}{e}\, \leq \, e^{-L \, d}\, \leq \, 1
\]
on $\Om_{1/L}$. We choose $M$ large such that
\[
M\left(1-\frac{1}{e}\right)\, >\,1 \, +\, G(\cdot,0,0)
\]
on $\Om$. We \textbf{assume} that there is a $L>0$ large such that
\[
G\Big(\cdot,Mae^{-L\, d}Dd,Mae^{-L\, d}\big(D^2d-a Dd \ot Dd\big) \Big) \, \leq\, 0
\]
for all $a\in [M/e,M]$, on $\Om_{1/L}$. Then, we obtain that $u^+$ is a supersolution on $\Om_{1/L}$ which satisfies the boundary condition. We extend $u^+$ to a global $\overline{u}$ on $\Om$ by setting
\[
\overline{u}\, :=\,
\left\{
\begin{array}{l}
\min\big\{u^+,K\big\}, \ \ \ \text{ on }\Om_{1/L},\\ 
K, \hspace{68pt} \text{ on }\Om \set \Om_{1/L},
\end{array}
\right.
\] 
where $K>0$ is chosen such that
\[
G(\cdot,0,0)\, \leq\, K \, \leq\, M\left(1-\frac{1}{e}\right).
\]
By invoking Theorem 6, we conclude.                \qed

\ms

\ms

\ms

\noi \textbf{Remarks on Chapter 5.} The method we expounded on in this chapter has been introduced by Oscar Perron in \cite{P} for the Laplacian. Ishii was the first to adapt this method to viscosity solutions in \cite{I2}. There exist ad hoc variants of this method which utilise the ``variational structure" of certain PDE and can be alternatively used to prove existence. For instance, a very nice application related to the material of these notes which provides existence for the $\infty$-Laplacian can be found in Aronsson-Crandall-Juutinen \cite{ACJ}. See also Crandall \cite{C1}. Our exposition in this chapter follows Crandall-Ishii-Lions \cite{CIL} closely.

\chapter[Comparison and Uniqueness]{Comparison results and Uniqueness of Solution to the Dirichlet Problem}

In this chapter we consider the problem of uniqueness of Viscosity Solutions to the Dirichlet problem
\beq \label{6.1}
\left\{
\begin{array}{l}
F\big(\cdot,u,Du,D^2u\big)\,=\,0, \, \text{ in }\Om,\ms\\
\hspace{78pt}u\, =\, b, \, \text{ on }\p \Om,
\end{array}
\right.
\eeq
for fully nonlinear degenerate elliptic PDE defined by the continuous coefficient
\[
F\ :\ \Om \by \R \by\R^n \by \mS(n) \larrow \R
\]
where $\Om\Subset \R^n$ and $b\in C^0(\overline{\Om})$. Uniqueness is obtained by means of a \textbf{Comparison Principle} between sub- and super- solutions. \textit{One of the most outstanding features of the theory is that}

{\center{

\fbox{\parbox[pos]{310pt}{
\ms

\centerline{Viscosity Solutions are as unique as Classical Solutions.} 

\ms
}}

\ms\ms }}

\noi By this statement we mean that, typically, a PDE whose classical solutions are unique, under certain assumptions its viscosity solutions are also unique. We begin by motivating the ideas to follow.

\ms

\noi \textbf{Motivation of the Comparison Principle.} \textit{Suppose that $u \in C^0(\overline{\Om})\cap C^2(\Om)$ is a classical solution of 
\[
F\big(\cdot,u,Du,D^2u\big)\,\geq \,0
\]
on $\Om$, $v \in C^0(\overline{\Om})\cap C^2(\Om)$ is a classical solution of 
\[
F\big(\cdot,v,Dv,D^2v\big)\,\leq \,0
\]
on $\Om$, and
\[
u\, \leq\, v \ \ \text{ on }\p \Om.
\]
We seek to show that $u\leq v$ in $\Om$. This will not succeed without some assumptions on $F$. Suppose that $F$ is degenerate elliptic and the partial function $r\mapsto F(x,r,p,X)$ is strictly decreasing:
\begin{align}
X\, \leq \, Y \ \ &\Longrightarrow \ \ F(x,r,p,X)\, \leq\, F(x,r,p,Y) \nonumber\\
r\, < \, s\ \ \ &\Longrightarrow \ \ F(x,r,p,X)\, >\, F(x,s,p,X). \nonumber
\end{align}
Suppose for the sake of contradiction that $u>v$ somewhere in $\Om$. Then, since $u\leq v$ on $\p \Om$, there exists $\bar{x} \in \Om$ such that $u-v$ has a (positive) maximum there:
\[
(u\, -\, v)(\bar{x})\, \geq\, u\, -\, v, \ \text{ on }\Om.
\]
\[
\includegraphics[scale=0.2]{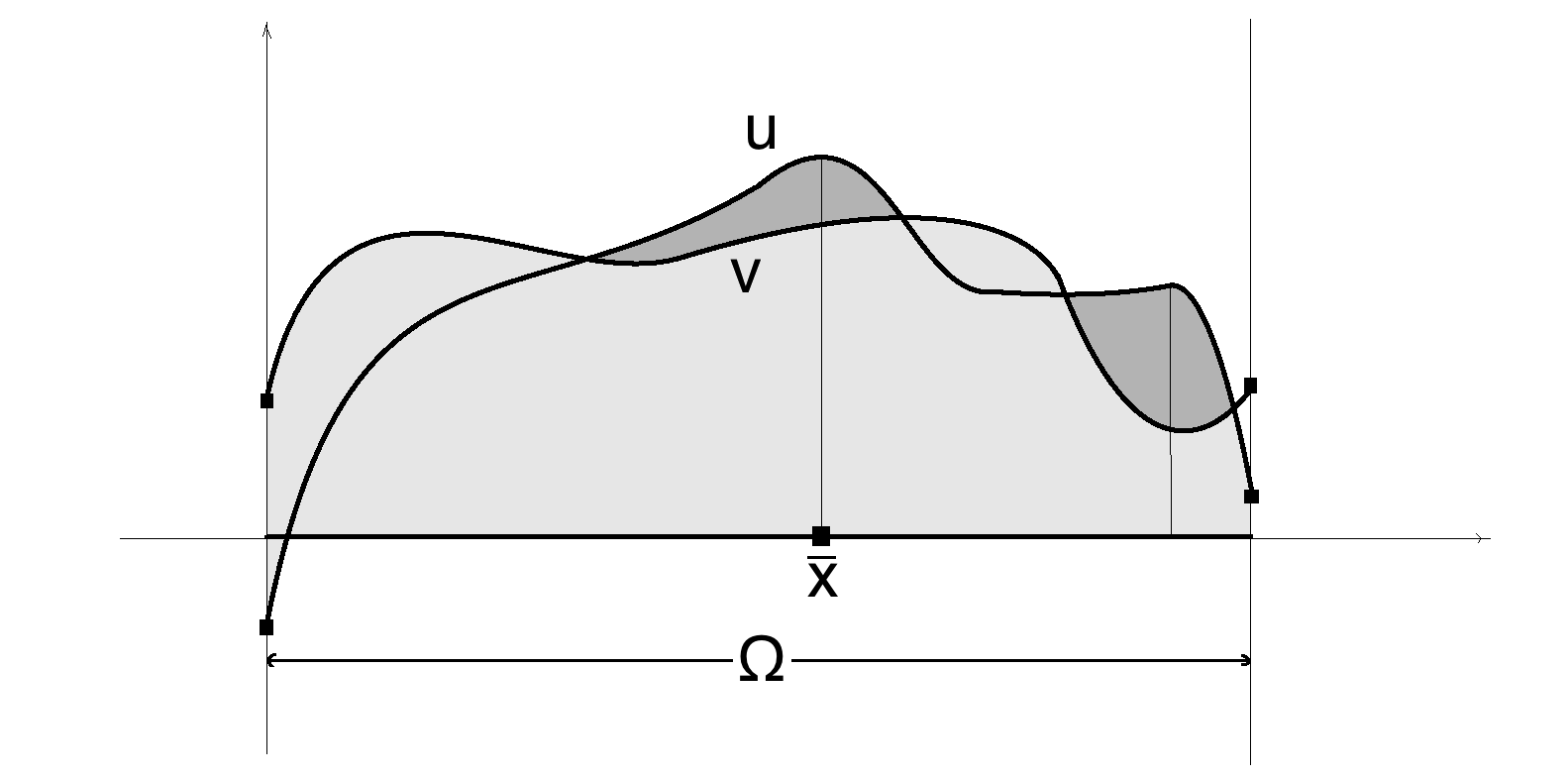}
\]
Hence, at $\bar{x}$ we have
\beq \label{6.2}
\left\{
\begin{array}{r}
D(u\, -\, v)(\bar{x})\, =\, 0, \ms\\
D^2(u\, -\, v)(\bar{x})\, \leq\, 0.
\end{array}
\right.
\eeq
By our assumptions on $F$, $u$, $v$ and in view of \eqref{6.2}, we have
\begin{align}  \label{6.3}
F\big(\bar{x},u(\bar{x}),Du(\bar{x}),D^2u(\bar{x})\big)\, \geq\, 0\, &\geq\, F\big(\bar{x},v(\bar{x}),Dv(\bar{x}),D^2v(\bar{x})\big) \nonumber\\
&=\, F\big(\bar{x},v(\bar{x}),Du(\bar{x}),D^2v(\bar{x})\big) \\
&\geq\, F\big(\bar{x},v(\bar{x}),Du(\bar{x}),D^2u(\bar{x})\big). \nonumber
\end{align}
Then, \eqref{6.3} gives
\[
F\big(\bar{x},u(\bar{x}),\bar{p},\bar{X}\big)\, \geq\, F\big(\bar{x},v(\bar{x}),\bar{p},\bar{X}\big)
\]
for $(\bar{p},\bar{X})=(Du(\bar{x}),D^2u(\bar{x}))$, which by anti-monotonicity of $r\mapsto F(x,r,p,X)$ contradicts that
\[
u(\bar{x})\, >\, v(\bar{x}).
\]
Hence, we conclude that
\[
u\, \leq \, v \ \text{ on }\overline{\Om}.
\]}
\textbf{The problem} with the above otherwise elegant argument is that 

{\center{

\fbox{\parbox[pos]{330pt}{
\ms

\centerline{THIS TECHNIQUE DOES NOT APPLY TO GENERAL $C^0$} 

\centerline{VISCOSITY SOLUTIONS, SINCE $\J^{2,\pm}u$, $\J^{2,\pm}v$ MAY WELL}

\centerline{BE \textbf{EMPTY} AT POINTS OF MAXIMA OF $u-v$ ! ! !}

\ms
}}

\ms\ms }}

\noi \textbf{The resolution.} The idea (which by now is pretty standard in uniqueness questions) is \textit{to double the number of variables and in the place of $u,v$, to consider instead maximisation of
\[
\overline{\Om} \by \overline{\Om} \ni \ (x,y)\,  \mapsto \, u(x)\, -\, v(y)\  \in \R
\]
on the product $\R^n \by \R^n$. Then we penalise the doubling of variables, in order to push the maxima on the diagonal $\{x=y\}$ of $\, \overline{\Om} \by \overline{\Om}$. The tool is to maximise the function
\[
W_\al(x,y)\,  := \, u(x)\, -\, v(y)\, -\, \frac{\al}{2}|x-y|^2
\]
and then let $\al \ri +\infty$.}

\ms

Let us now state the main results of this chapter. For simplicity, we restrict attention to equations with \textit{decoupled dependence in $x$}:
\[
F(u,Du,D^2u)\, =\, f
\]
and to \textit{continuous} sub-/super- solutions. The results however hold true for the general case of fully coupled dependence and semi-continuous sub-/super- solutions, but extra assumptions are required and the proof is a little more complicated. We have

\ms

\noi \textbf{Theorem 1 (Comparison Principle for Viscosity Solutions).}  \emph{Let $F\in C^0\big(\R \by \R^n \by \mS(n)\big)$, $f\in C^0(\overline{\Om})$, $\Om \Subset \R^n$ and assume that
\begin{align}
X\, \leq \, Y \ \ &\Longrightarrow \ \  F(r,p,X)\, \leq\, F(r,p,Y), \nonumber\\
r\, \leq \, s\ \ \ &\Longrightarrow \ \ F(r,p,X)\, \geq\, F(s,p,X) \, +\, \ga(s-r), \nonumber
\end{align}
for some $\ga>0$. Suppose that $u\in C^0(\overline{\Om})$ is a Viscosity Solution of 
\[
F(u,Du,D^2u)\, \geq\, f
\]
on $\Om$, and $v\in C^0(\overline{\Om})$ is a Viscosity Solution of 
\[
F(v,Dv,D^2v)\, \leq\, f
\]
on $\Om$. Then, we have
\[
u \, \leq\, v \ \text{ on }\p \Om\ \ \Longrightarrow\ \ u \, \leq\, v \ \text{ on } \overline{\Om}.
\]
}
Theorem 1 readily implies 

\ms

\noi \textbf{Theorem 2 (Uniqueness for the Dirichlet problem).}  \emph{Let $F\in C^0\big(\R \by \R^n \by \mS(n)\big)$, $f,b \in C^0(\overline{\Om})$, $\Om \Subset \R^n$ and assume that
\begin{align}
X\, \leq \, Y \ \ &\Longrightarrow \ \  F(r,p,X)\, \leq\, F(r,p,Y), \nonumber\\
r\, \leq \, s\ \ \ &\Longrightarrow \ \ F(r,p,X)\, \geq\, F(s,p,X) \, +\, \ga(s-r), \nonumber
\end{align}
for some $\ga>0$.  Then, the boundary value problem
\[
\left\{
\begin{array}{l}
F\big(u,Du,D^2u\big)\,=\,f, \, \text{ in }\Om,\ms\\
\hspace{69pt}u\, =\, b, \, \text{ on }\p \Om,
\end{array}
\right.
\]
admits at most one Viscosity Solution $u\in C^0(\overline{\Om})$.
}

\ms

\noi \textbf{Remark 3.} \textit{Sufficient conditions for the existence of solution to the Dirichlet problem are given in Chapters 3 and 5, either by means of stability and approximation, or by the Perron method. In particular, Perron's method implies existence of solution by assuming uniqueness and existence of a pair of sub-/super- solutions.}

\ms

Now we start working towards the proof of Theorem 1. The first step is the next device, which ``pushes the maxima to the diagonal", in the slightly more general form, which is needed later:

\ms

\noi \textbf{Lemma 4 (Penalisation).} \emph{Let $\mO \Subset \R^N$ for some $N\geq 1$, $\Phi,\Psi \in C^0(\overline{\mO})$, $\Psi \geq 0$. For $\al>0$, we assume that $\Phi -\al \Psi$ has a maximum at $z_\al$:
\[
\Phi \,-\, \al \Psi\, \leq \, (\Phi \,-\, \al \Psi)(z_\al),
\]
in $\overline{\mO}$. Then,  
\[
\lim_{\al \ri +\infty}\, \al \Psi(z_\al) \, =\, 0
\]
and if $z_\al \ri \hat{z}$ along a sequence as $\al \ri \infty$, we have $\Psi(\hat{z})=0$ and
\[
\max_{\{\Psi=0\}}\, \Phi\, =\, \Phi(\hat{z})\, =\,  \lim_{\al \ri +\infty}\, \max_{\overline{\mO}}\, \big\{\Phi \,-\,\al \Psi \big\}.
\]
}

\noi \textbf{Proof.} We set
\[
M_\al \, :=\, \max_{\overline{\mO}}\, \big\{\Phi \, -\, \al \Psi \big\}.
\]
The limit $\lim_{\al \ri \infty}M_\al$ exists and is finite since $\al \mapsto M_\al$ is decreasing and 
\[
+\infty\, >\, M_\al \, \geq \, (\Phi \,-\,\al \Psi)(z_0)\, >\, -\infty,
\]
for any $z_0 \in \mO$. Moreover,
\begin{align}
M_{\al/2}\, & =\, \max_{\overline{\mO}}\, \Big\{\Phi \, -\, \frac{\al}{2} \Psi \Big\} \nonumber\\
&\geq\, \Phi(z_\al)\, -\, \frac{\al}{2}\Psi(z_\al)  \nonumber\\
& =\, \Phi(z_\al)\, -\, \al\Psi(z_\al) \, +\, \frac{\al}{2}\Psi(z_\al) \nonumber\\
& =\, M_\al \,  +\, \frac{\al}{2}\Psi(z_\al) . \nonumber
\end{align}
Thus, as $\al \ri \infty$, we get
\[
\underset{\al \ri \infty}{\lim\sup}\, \big( \al\Psi(z_\al)\big) \, = \, 2\, \underset{\al \ri \infty}{\lim\sup}\, \big( M_{\al/2}\, -\, M_\al \big) \, = \, 0.
\]
Then, we have
\begin{align}
\max_{\{\Psi=0\}}\, \Phi \, &=\, \sup_{\{z\, :\, \Psi(z)=0\}}\, \big\{ \Phi(z)\, -\, \al \Psi(z) \big\}  \nonumber\\
 &\leq\, \max_{z\in \overline{\mO}}\, \left\{\Phi(z)\, -\, \al \Phi(z)\right\}  \nonumber\\
&=\,  M_\al \nonumber\\
&=\,  \Phi(z_\al)\, -\, \al\Psi(z_\al)  \nonumber\\
&\leq\, \Phi(z_\al), \nonumber
\end{align}
as $\al \ri \infty$. Suppose that $z_\al \ri \hat{z}$ along a sequence of $\al$'s. Then, by passing to the limit we get $\Psi(\hat{z})=0$ and
\[
\sup_{\{\Psi=0\}}\, \Phi \, =\,  \lim_{\al \ri +\infty}\, M_\al \, =\, \Phi(\hat{z}). 
\]
The lemma ensues.     \qed

\ms

\noi \textbf{Proof of Theorem 1.} Assume for the sake of contradiction that 
\beq \label{6.4}
\sup_{\Om}\, (u\, -\, v) \, \equiv\, \de \, >\, 0. 
\eeq
We fix $\al>0$ and define $W_\al \in C^0(\overline{\Om} \by \overline{\Om})$ by
\[
W_\al(x,y) \, :=\, u(x)\, -\, v(y)\, -\, \frac{\al}{2}|x-y|^2.
\]
Let $(x_\al,y_\al) \in \Om \by \Om$ be such that
\[
W_\al(x_\al,y_\al) \, =\, \sup_{\Om \by\Om}\, W_\al.
\]
\textbf{We first consider the case under the extra assumption that the sub-/super- solutions are in $C^2(\Om)$}. Then, \textbf{since $W_\al$ in twice continuously differentiable}, at $(x_\al,y_\al)$ we have
\begin{align}
DW_\al (x_\al,y_\al)\, &=\,0,\ \text{ in }\R^{2n},  \nonumber\\
 D^2W_\al (x_\al,y_\al)\, &\leq \,0,\ \text{ in }\mS(2n), \nonumber
\end{align}
which in turn says that
\[
\left[
\begin{array}{c}
Du(x_\al)\\
-Dv(y_\al)
\end{array}
\right]
\, =  \,  
\left[
\begin{array}{c}
\al(x_\al - y_\al)\\
\al(y_\al - x_\al)
\end{array}
\right]
\]
and
\[
\left[
\begin{array}{cc}
D^2u(x_\al) & 0\\
0 & -D^2v(y_\al)
\end{array}
\right]
\, \leq  \,  
\al
\left[
\begin{array}{cc}
I & -I\\
-I & I
\end{array}
\right]  .
\]
We set:
\[
\left[
\begin{array}{c}
p_\al\\
q_\al
\end{array}
\right]
\, :=  \,  
\left[
\begin{array}{c}
Du(x_\al)\\
Dv(y_\al)
\end{array}
\right]
\]
and
\[
\left[
\begin{array}{c}
X_\al\\
Y_\al
\end{array}
\right]
\, :=  \,  
\left[
\begin{array}{c}
D^2u(x_\al)\\
D^2v(y_\al)
\end{array}
\right] .
\]
Then, we have
\beq \label{6.5}
p_\al \, =\, \al(x_\al - y_\al)\, =\, q_\al
\eeq
and also by testing the matrix inequality against vectors of the form $(\xi,\xi)\in \R^{2n}$, we have
\begin{align}
\big(X_\al \, -\, Y_\al \big):\xi \ot \xi \, &=\,
\left[
\begin{array}{cc}
X_\al & 0\\
0 & -Y_\al
\end{array}
\right] 
:
\left[
\begin{array}{c}
\xi  \\
\xi
\end{array}
\right]
\ot
\left[
\begin{array}{c}
\xi\\
\xi
\end{array}
\right]
 \nonumber\\
& \leq  \,  
\al
\left[
\begin{array}{cc}
I & -I\\
-I & I
\end{array}
\right] 
:
 \left[
\begin{array}{c}
\xi\\
\xi
\end{array}
\right]
\ot
\left[
\begin{array}{c}
\xi\\
\xi
\end{array}
\right]
 \nonumber\\
&=\, 0.  \nonumber
\end{align}
Hence,
\beq \label{6.6}
X_\al \, \leq\,  Y_\al
\eeq
in $\mS(n)$. By using the fact that $u$ is a subsolution and that $v$ is a supersolution and that
\begin{align}
(p_\al,X_\al)\in \J^{2,+}u(x_\al), \label{6.6a}\\
 (q_\al,Y_\al)\in \J^{2,-}v(y_\al). \label{6.6b}
\end{align}
we have
\[
F\big(u(x_\al),p_\al,X_\al \big)\, -\, F\big(v(y_\al),q_\al, Y_\al \big) \, \geq\, f(x_\al)-f(y_\al).
\]
By using \eqref{6.5}, \eqref{6.6} and ellipticity, we have 
\begin{align}
F\big(u(x_\al),\al(x_\al-&y_\al),X_\al \big)\, -\, f(x_\al) \nonumber\\
 &\geq \, F\big(v(y_\al),\al(x_\al-y_\al),X_\al\big)\, -\, f(y_\al) .\nonumber
\end{align}
By using our assumption, we have
\begin{align} \label{6.7}
f(x_\al)-f(y_\al)\, & \leq \, F\big(u(x_\al),\al(x_\al-y_\al),X_\al \big)  \nonumber\\ 
&\ \ \ -\, F\big(v(y_\al),\al(x_\al-y_\al),X_\al \big)  \\
&\leq\, -\ga \big( u(x_\al)\, -\, v(y_\al)\big).  \nonumber
\end{align}
We now apply Lemma 4 to
\begin{align}
\mO\, :=\, \Om &\by \Om \, \sub\, \R^{n+n}\,=\,\R^N, \nonumber\\
z\,=\, & (x,y)\in \Om \by \Om, \nonumber\\
 \ \Phi(z)\,& :=\, u(x)\,-\,v(y), \nonumber\\
\Psi(z)\, &:=\, \frac{1}{2}|x-y|^2. \nonumber
\end{align}
to obtain that along a sequence of $\al$'s  we have 
\[
(x_\al,y_\al)\larrow (\hat{x},\hat{x}),\ \ \text{ as $\al \ri \infty$}
\]
and
\beq \label{6.8}
\sup_{\Om}\, \{u\,-\, v\}  \, =\, (u\, -\, v)(\hat{x})  \, =\,  \lim_{\al \ri +\infty}\, W_\al(x_\al, y_\al). 
\eeq
Then, by using \eqref{6.4} and \eqref{6.8}, our passage to the sequential limit as $\al \ri \infty$ in \eqref{6.7} gives the contradiction
\[
0\,=\, f(\hat{x})- f(\hat{x})\, \leq \, -\ga\big(u(\hat{x})-v(\hat{x})\big)\, =\, -\ga \de \, <\, 0. 
\]
Hence,
\[
u\, \leq \, v\ \text{ on }\overline{\Om}. 
\]
Theorem 1 will be fully proved once we \textbf{remove the assumption that $W_\al$ is smooth}, which was used in order to derive relations \eqref{6.5}, \eqref{6.6} and in turn to construct the Jets of $u$, $v$ in \eqref{6.6a},  \eqref{6.6b}. This gap is bridged right afterwards.   \qed

\ms

In order to complete the proof in the general case, we need

{\center{

\fbox{\parbox[pos]{330pt}{
\ms

\centerline{ \textbf{a strong tool which produces Jets of sums of functions}. } 

\ms
}}

\ms\ms }}

In order to state it, we need the next slight extension of the notion of jets:

\ms

\noi \textbf{Definition 5 (Jet Closures).} \emph{Let $u\in C^0(\Om)$ and $x\in \Om \sub \R^n$. The Jet closures $\overline{\J}^{2,\pm}u(x)$ of ${\J}^{2,\pm}u(x)$ are defined by
\begin{align}
\overline{\J}^{2,\pm}u(x)\, :=\, \Big\{ &(p,X)\in\, \R^n \by \mS(n)\ \Big|\ \exists\ (x_m,p_m,X_m)\ri (x,p,X)  \nonumber\\
&\ \text{ as }m\ri \infty,\ \text{ such that }\ (p_m,X_m)\in \J^{2,\pm}u(x_m)\Big\}.  \nonumber
\end{align}
}
The idea of Jet closures is that we ``attach jets of nearby points". Obviously,
\[
{\J}^{2,\pm}u(x) \, \sub  \, \overline{\J}^{2,\pm}u(x).
\]
It is a simple continuity fact that the definition of Viscosity Solutions remains the same if one replaces ${\J}^{2,\pm}u(x)$ by $\overline{\J}^{2,\pm}u(x)$. We record this fact in the next

\ms

\noi \textbf{Lemma 6 (An equivalent definition of Viscosity Solutions).} \emph{Let $u\in C^0(\Om)$, $x\in \Om \sub \R^n$ and $F\in C^0(\Om \by \R\by \R^n \by \mS(n))$. Then, $u$ is Viscosity Solution of 
\[
F\big(\cdot,u,Du,D^2u\big)\, =\, 0
\]
on $\Om$, if and only if 
\begin{align}
(p,X)\in\overline{\J}^{2,+}u(x)\ \ \Longrightarrow\ \ F(x,u(x),p,X)\, \geq \, 0,  \nonumber\\
(p,X)\in\overline{\J}^{2,-}u(x)\ \ \Longrightarrow\ \ F(x,u(x),p,X)\, \leq \, 0, \nonumber
\end{align}
for all $x\in \Om$.}

\ms

\noi \textbf{Proof.} The one direction is obvious. For the converse, let $(p,X)\in \overline{\J}^{2,+}u(x)$. Then, there exists a sequence  such that
\[
(x_m,p_m,X_m)\larrow (x,p,X)
\]
as $m\ri \infty$, while $ (p_m,X_m)\in \J^{2,+}u(x_m)$. Since $u$ is a Viscosity Subsolution, we have
\[
F\big(x_m,u(x_m),p_m,X_m \big)\, \geq \, 0.
\]
By letting $m\ri \infty$, we get $F\big(x,u(x),p,X \big) \geq 0$. The supersolution property is analogous and hence we are done.   \qed

\ms

We now have the

\ms

\noi \textbf{Theorem 7 (Theorem on Sums).} \emph{Let $u_i\in C^0(\overline{\Om_i})$, $\Om_i \Subset \R^{n_i}$, $1\leq i \leq k$. For $N:=n_1 +\cdots+n_k$ and $x=(x_1,\cdots,x_k)\in \R^N$, we set
\[
u(x)\, :=\, u_1(x_1)\, +\, ... \,+\, u_k(x_k),\ \ u\in C^0(\overline{\Om_1} \by ... \by \overline{\Om_k}).
\]
Suppose that 
\[
\big(p_1,...,p_k,X \big) \, \in \, \J^{2,+}u(x)
\]
for some $x \in \Om_1 \by ... \by \Om_k$.} 

\emph{Then, for each $\e>0$, there exist matrices $ X_i \in \mS(n_i)$ such that, for each $i=1,...,k$ we have 
\[
(p_i,X_i)\, \in \, \overline{\J}^{2,+}u_i(x_i)
\]
and
\beq \label{6.11}
-\left(\frac{1}{\e}\, +\, \|X\| \right)\, I\, \leq  \,  
\left[
\begin{array}{ccc}
X_1 & \cdots & 0\\
\vdots & \ddots & \vdots\\
0 & \cdots & X_k
\end{array}
\right]\,
\leq\, X\,+ \,\e X^2,
\eeq
in $\mS(N)$.}

\ms

Note carefully that the jets constructed are actually \textit{in the jet closures, rather than in the jets themselves}. We conclude this chapter with the proof of Theorem 7. We caution the reader that

{\center{

\fbox{\parbox[pos]{280pt}{
\ms

The proof of the Theorem on Sums is \textbf{more advanced} 

than the rest of the material appearing in these notes 

\ \ \ \ \ \ \ \ \ and should be \textbf{omitted} at first reading.

\ms
}}

\ms\ms }}

\noi \textbf{Proof of Theorem 7.} We begin by observing that we may in addition assume that
\begin{align}
x\, =\, (0,...,0)&\, \in\, \Om_1\by ... \by \Om_k, \nonumber
\\
u(0)\, &=\, 0,  \nonumber \\
(p_1, ...,p_k)\, &=\, (0,...,0) \,\in\, \R^{n_1 +... +n_k},  \nonumber\\
u(z)\, \leq\, \,\frac{1}{2} & X:z\ot z,\ \text{ as }z\ri 0 \ \ \ \text{ \big(no $o(|z|^2)$ terms in $\J^{2,+}$\big)}.  \nonumber
\end{align}
Indeed, it suffices to consider the result for the functions
\[
\tilde{u}_i(z_i)\,:=\, u_i(z_i+x_i) \, -\, u_i(x_i)\, -\, p_i \cdot z_i \,+\, \frac{\si}{2}|z_i|^2
\]
for $\si>$ small, and set
\[
\tilde{u}(z)\, :=\, \tilde{u}_1(z_1)\, +\, ... \, +\, \tilde{u}_k(z_k).
\]
Indeed, the normalisation technique is quite obvious. Hence, under these normalisation, it suffices to prove that given
\[
(0,...,0,X)\in \J^{2,+}u(0),
\]
there exist $k$ matrices $X_i \in \mS(n_i)$ such that
\[
(0,X_i) \in \overline{\J}^{2,+}u_i(0)
\]
and also  \eqref{6.11} holds. Suppose now that Theorem 7 holds in the case where the super-quadratic terms vanish. Then, for the general case \emph{with} super-quadratic terms, we have that if $|z|$ is small, then
\[
o(|z|^2)\, \leq \, \frac{\si}{2}|z|^2
\]
near the origin. Hence, if $(0,...,0,X)\in \J^{2,+}u(0)$, then
\[
u(z)\, \leq\, \,\frac{1}{2} (X+\si I):z\ot z
\]
and evidently $(0,...,0,X+\si I)\in \J^{2,+}u(0)$. We then apply the result to conclude that there are $X_1(\si),...,X_k(\si)$ such that $(0,X_i(\si)) \in \overline{J}^{2,+}u_i(0)$. By using the bounds \eqref{6.11} we can let $\si \ri 0$ and conclude in the general case.

The next step is the next matrix inequality:

\ms

\noi \textbf{Claim 8.} \emph{For any $X\in \mS(N)$, $a,b\in \R^N$ and $\e>0$, we have
\[
X:a \ot a \, \leq\, \big(X+\e X^2 \big):b\ot b\ +\ \left(\frac{1}{\e}\, +\, \|X\| \right)|b-a|^2
\]
where $\| \cdot\|$ is the following norm on $\mS(N)$:
\[
\|X\| \, :=\, \max |\si(X)|\, =\, \max_{|\xi| =1}|X: \xi \ot \xi |.
\]
}

\noi \textbf{Proof.} By the Cauchy-Schwartz inequality, we have
\begin{align}
X:a\ot a\, &=\, X:b \ot b \ +\ X:(a-b)\ot(a-b) \ +\ 2 X:(a-b)\ot b  \nonumber\\
&\leq\, X:b \ot b\ +\ \|X\|\,|b-a|^2\ +\ 2\left(\frac{a-b}{\sqrt{\e}}\right) \cdot \big(\sqrt{\e}Xb \big)  \nonumber\\
&\leq\, X:b \ot b\ +\ \|X\|\,|b-a|^2\ +\ \frac{|a-b|^2}{\e}\ +\ \e|Xb|^2.  \nonumber
\end{align}
By using the identity
\[
|Xb|^2\, =\, X_{ik}b_k\, X_{il}b_l\, =\, X_{ik} X_{il}\, b_k\,b_l\,=\, X^\top\! X :b \ot b
\]
and that $X^\top=X$, we get
\begin{align}
X:a\ot a\, &\leq\, \big(X \, +\, \e X^\top\! X\big):b \ot b\ +\ \left(\frac{1}{\e}\, +\, \|X\| \right)|b-a|^2 \nonumber\\
&=\, \big(X \, +\, \e X^2\big):b \ot b\ +\ \left(\frac{1}{\e}\, +\, \|X\| \right)|b-a|^2. \nonumber
\end{align}
The claim follows.   \qed

\ms

We now mollify things, by considering the sup-convolution approximations introduced in Chapter 4. We set
\[
\la\, :=\, \left(\dfrac{1}{\e}\, +\, \|X\| \right)^{-1}
\]
and consider the sup-convolutions $u_1^\la$, ..., $u_k^\la$ and $u^\la$ of $u_1$, ..., $u_k$ and $u$ respectively:
\begin{align} \label{6.12}
u_i^\la(z_i)\, &:=\, \sup_{y_i \in \Om_i}\left\{u_i(y_i)\, -\, \frac{|z_i-y_i|^2}{2\la} \right\},\ \ z_i \in \Om_i,\\  \label{6.13}
u^\la(z)\, &:=\, \sup_{y \in \Om_1 \by ... \by \Om_k}\left\{u(y)\, -\, \frac{|z-y|^2}{2\la} \right\}.
\end{align}
By orthogonality, we have
\[
|z_1-y_1|^2\, +...\, +\, |z_k-y_k|^2\, =\, |z-y|^2
\]
and hence by \eqref{6.12} and  \eqref{6.13}, we obtain
\[
u_1^\la(z_1)\, +\, ...\,+\, u^\la_k(z_k)\, =\, u^\la(z).
\]
By using Claim 8, eqref{6.13} and the definition of $\la$, we deduce
\beq \label{6.14}
u^\la(z)\, \leq\, \frac{1}{2}\big(X\, +\, \e X^2\big):z\ot z,
\eeq
on $\Om$. 

We now state and prove the following important general fact bout semiconvex functions, which we use right afterwards.

\ms

\noi \textbf{Lemma 9 (Jensen's Lemma).} \emph{Let $f\in C^0(\Om)$ be semiconvex on $\Om \sub \R^N$ with semiconvexity constant $\la>0$. Suppose that $f$ has a strict local maximum at some $x\in \Om$. For $\de,\rho>0$ small, we set
\[
K_{\de,\rho}\, :=\, \Big\{y\in \mB_\rho (x)\, \Big| \, \exists\ p \in \mB_\de(0)\, :\, z\mapsto f(z)+p\cdot(z-y)\, \text{ has max at }y \Big\}.
\]
Then, the Lebesgue measure of $K_{\de,\rho}$ is positive and we have the estimate
\[
\big|K_{\de,\rho}\big|\, \geq\, \al(n)(\la\de)^N
\]
where $\al(N)$ is the volume of the unit ball of $\R^N$.}

\ms

\[
\includegraphics[scale=0.2]{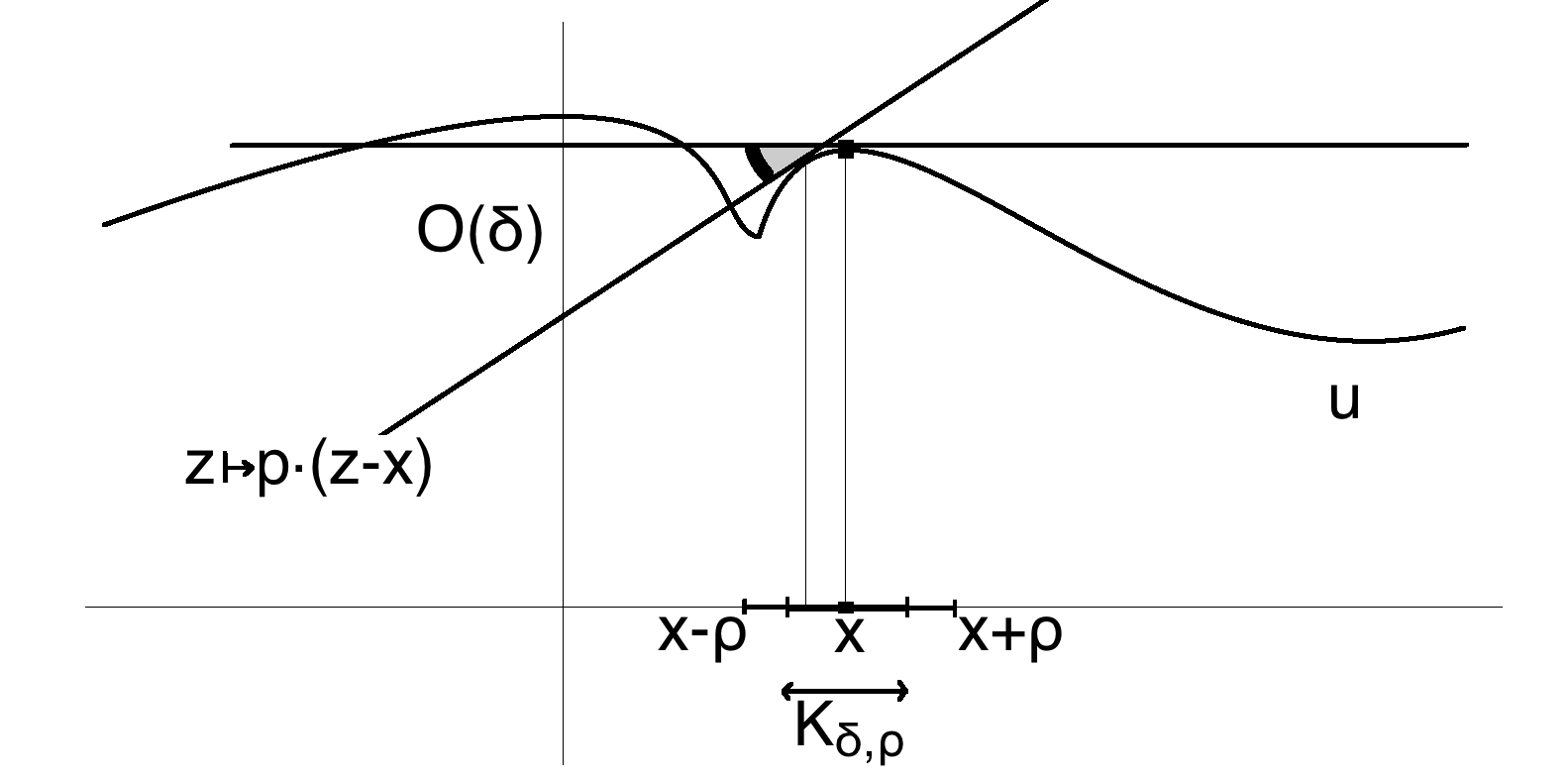}
\]

\noi \textbf{Proof.} We begin by observing that $K_{\de,\rho}\neq \emptyset$. Indeed, this is clear for $\de$ small by touching from above $f$ by the affine function $z\mapsto a+ p\cdot z$. Assume first that $f\in C^2(\Om)$. Then, if $|p|<\de$, at maxima of $f$ at $y$ we have $p=Df(y)$ and hence $y \in K_{\de,\rho}$. Thus,
\[
\mB_\de (0)\, \sub\, Df(K_{\de,\rho})
\]
and if $y\in K_{\de,\rho}$, then
\[
-\frac{1}{\la}I\, \leq \, D^2f(y)\, \leq\, 0.
\]
The left inequality is by semiconvexity and the right by the fact that $f$ has a max at $y$. Hence,
\[
 -\frac{1}{\la^N}\, \leq\, \det(D^2f) \, \leq\, 0.
\]
By the following well-known inequality
\[
\int_A\big|\det(D^2f) \big|\, \geq \, |Df(A)|
\]
which holds for any measurable set $A\sub \Om$ and is a consequence of the celebrated \emph{Co-Area Formula}, we have
\begin{align}
\al(N)\de^N\, &=\, |\mB_\de(0)|  \nonumber\\
&\leq\, \big| Df(K_{\de,\rho}) \big|  \nonumber\\
&\leq\, \int_A\big|\det(D^2f) \big|  \nonumber\\
&\leq \, |K_{\de,\rho}| \frac{1}{\la^N} . \nonumber
\end{align}
As a consequence, the estimate follows in that case. For the general case, we mollify $f$ in the standard way (see the beginning of Chapter 4) to $f^\e = f *\eta^\e$ and let $K^\e_{\de,\rho}$ denote the respective sets corresponding to $f^\e$. Then, we have
\[
K_{\de,\rho}\, \supseteq \, \bigcap_{i=1}^\infty \bigcup_{m=i}^\infty K^{1/m}_{\de,\rho}
\]
and hence
\[
\al(N)(\de \la)^N\, \leq\, \underset{i\ri \infty}{\lim} \left|\bigcup_{m=i}^\infty K^{1/m}_{\de,\rho} \right|\, \leq\, | K_{\de,\rho} |.
\]
The lemma ensues.  \qed

\ms

Now we use Lemma 9 in order to prove the next

\ms

\noi \textbf{Lemma 10 (Existence of ``approximate" jet).} \emph{Let $f\in C^0(\Om)$ be semiconvex on $\Om \sub \R^N$ with semiconvexity constant $\la>0$. If 
\[
f(z+x) \, \leq \, f(x) \, +\, \frac{1}{2}A:z\ot z
\]
for some $A\in \mS(N)$ near $x \in \Om$, then there exists $X\in \mS(N)$ such that
\[
(0,X) \, \in \, \overline{\J}^{2,+}f(x) \bigcap  \overline{\J}^{2,-}f(x)
\]
and also
\[
-\frac{1}{\la}I\, \leq\, X\, \leq\, A.
\]}

\noi \textbf{Proof.} Since $f$ is semiconvex, by Alexandroff's theorem is twice differentiable a.e.\ on $\Om$. By Jensen's Lemma, if $\de>0$ is small, the set of points $y_\de$ for which  $|p_\de|<\de$ and
\[
(p_\de,A)\, \in \, \J^{2,+}f(y_\de)
\]
has positive measure. Hence, there exist points $y_\de$ at which $f$ is twice differentiable and consequently, there we have
\[
|Df(y_\de)|\,  <\, \de \ , \ \ \ \  -\frac{1}{\la}I\, \leq\, D^2 f(y_\de)\, \leq\, A.
\]
By letting $\de \ri 0$, by compactness we get that there is $X\in \mS(N)$ such that along a sequence
\[
Df(y_\de) \ri 0 \ , \ \ \ \ D^2 f(y_\de)\ri X.
\]
In view of the definition of Definition 5, the lemma follows.  \qed

\ms

We may now conclude the proof of Theorem 7 by noting that the sup-convolutions $u_i^\la$ of each $u_i$ is $\la$-semiconvex. Hence, by applying Lemma 10 to  $u_i^\la$ for $x_i=0$ and
\[
A\, :=\, X\, +\e X^2,
\]
we get that there exists $X_i \in \mS(n_i)$ such that $(0,X_i) \in \overline{\J}^{2,\pm}u_i^\la(0)$. Then, the ``magic property" of sup-convolution in Theorem 7 of Chapter 4 implies that
\[
(0,X_i) \in \overline{\J}^{2,\pm}u_i(0).
\]
The desired matrix inequality follows by \eqref{6.14} and by the fact that $u^\al$ is semiconvex with constant $\la=1/(1\e+\|X\|)$. Hence, Theorem 7 ensues.    \qed
 
\ms

\ms

\ms

\noi \textbf{Remarks on Chapter 6.} This material of this chapter is taken from Crandall-Ishii-Lions \cite{CIL}, but the exposition is more thorough and we analyse technical points which in the original article are omitted. In \cite{CIL} the reader can also find a more general version of our Theorems 1 and 2, without decoupling the dependence in $x$, at the expense of some technical complexity (and of an extra assumption on the nonlinearity $F$). In these notes we do not consider the issue of uniqueness of Hamilton-Jacobi PDE. For this and further material (like for example uniqueness in unbounded domains), we redirect the reader to the article of Crandall \cite{C2} and to the notes of Koike \cite{Ko}. Viscosity solution apply to \textbf{monotone systems} as well, namely systems for which each equation of the system depends only on the derivatives of the respective component of the solution: $F_\al(\cdot,u_1,...,u_N,Du_\al,D^2u_\al)=0$, $\al=1,...,N$. For such systems the viscosity notions apply directly via ``component-wise" arguments, and so does the Perron method and the uniqueness theory. For details see Ishii-Koike \cite{IK}.

\chapter*{Applications}

\chapter[Minimisers and Euler-Lagrange PDE]{Minimisers of Convex Functionals and Existence of Viscosity Solutions to the Euler-Lagrange PDE}

In this chapter we consider applications of the theory of Viscosity Solutions to Calculus of Variations and to the Euler-Lagrange PDE. Some basic graduate-level familiarity with weak derivatives and functionals will be assumed, although all the basic notions will be recalled. The results herein will be needed in the next chapter as well, when we will consider the existence of solution to the $\infty$-Laplacian. 

\ms

\noi \textbf{Some basic ideas from Calculus of Variations.} Let $F \in C^2(\R^n)$ be a convex function and $\Om \sub \R^n$. We consider the \textbf{functional}
\beq \label{7.1}
E\  :\ W^{1,1}_{\text{loc}}(\Om) \by \B(\Om)  \larrow [0,\infty]
\eeq
defined by
\begin{align}
E(u,\B)\, &:=\, \int_\B F(Du),\  \text{ when\  }\ F(Du) \in L^1(\B), \ \B \in \B(\Om), \nonumber \ms\\
 &:=\,+\infty \hspace{40pt} \text{ otherwise.} \nonumber
\end{align}
Here $\B(\Om)$ denotes the family of Borel subsets of $\Om$ and $W^{1,1}_{\text{loc}}(\Om)$ denotes the Sobolev space of locally integrable weakly differentiable functions on $\Om$ with  locally integrable weak derivative, that is
\[
W^{1,1}_{\text{loc}}(\Om)\, :=\, \Big\{u\in L^1_{\text{loc}}(\Om)\ \big| \ \exists \ Du=(D_1,u,...,D_nu) \ : \ D_iu \in L^1_{\text{loc}}(\Om) \Big\}.
\]
We recall that for $r\geq 1$, the space of locally $r$-integrable functions is
\[
L^r_{\text{loc}}(\Om)\, := \,\left\{u \, : \, \Om\ri\R \ \text{(Lebesgue) measurable } \Big|\ \forall\, \Om' \Subset \Om,\, \int_{\Om'}|u|^r<\infty \right\},
\]
and the measurable vector function  $Du \ :\ \Om \sub \R^n \larrow \R^n$
is called the \textbf{weak derivative} of $u$, when for all $\psi \in C^1_c(\Om)$, we have
\[
\int_\Om u D\psi \, =\, -\int_\Om \psi Du.
\]
We note that

{\center{

\fbox{\parbox[pos]{260pt}{
\ms

\centerline{In effect, weak differentiation means that} 
\ms

\centerline{\!``integration by parts"  is taken as definition} 

\ms

\centerline{of a non-classical derivative !}
\ms
}}

\ms\ms }}

\noi The Sobolev space $W^{1,r}_{\text{loc}}(\Om)$ for $r>1$ can be defined analogously, by replacing $L^1$ with $L^r$. A ``functional" is nothing but a function defined on an ``infinite-dimensional space". \textbf{One of the most central concepts of Calculus of Variations is that one can find a solution $u$ of the PDE}
\beq \label{7.2}
\Div \big(F_A(Du) \big)\, =\, 0
\eeq
(where $F_{A}$ stands for the derivative of $F$ on $\R^n$) \textbf{by finding a $u$ which minimises the functional \eqref{7.1}.} The PDE \eqref{7.2} is the celebrated \textbf{Euler-Lagrange equation}. The heuristic idea is

{\center{

\fbox{\parbox[pos]{200pt}{
\ms

\centerline{AT MINIMUM POINTS OF $E$,} 
\ms

\centerline{WE HAVE ``$E'=0$" !} 

\ms
}}

\ms\ms }}

\ms

The finite-dimensional analogue of this idea is depicted in the next figure.

\includegraphics[scale=0.18]{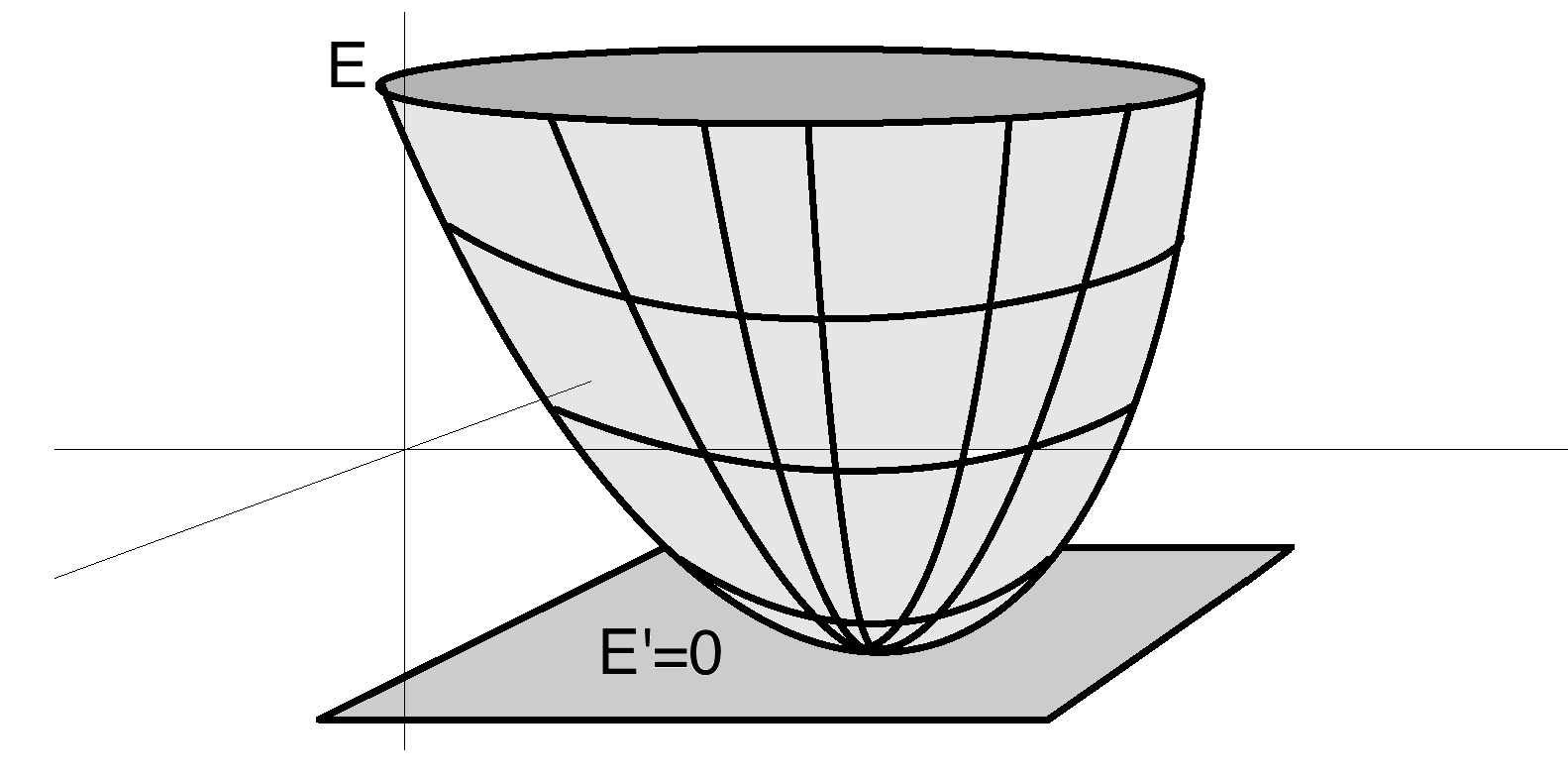}

\noi \textbf{Formal Argument.} Indeed, in order to see that, suppose
\[
E(u,\Om)\, \leq \, E(u+\psi,\Om),\  \text{ for }\psi \in C^1_c(\Om).
\]
Fix such a $\psi$ and for $\e \in \R$, consider the function
\[
f(\e)\, :=\, E(u+\e \psi,\Om'),\ \  f\ :\ \R \ri \R. 
\]
Then, the previous inequality says that
\[
f(0)\, \leq \,f(\e), \ \text{ for }\e\in \R.
\]
Consequently, $f$ has a minimum at zero, which gives $f'(0) = 0$. As a result,
\begin{align}
f'(0)\, &=\, \frac{d}{d\e}\Big|_{\e=0} E(u+\e \psi,\Om) \nonumber\\
&=\, \int_{\Om'}\frac{d}{d\e}\Big|_{\e=0} F(Du+\e D\psi) \nonumber \\
&=\, \int_{\Om} F_A(Du) \cdot D\psi. \nonumber
\end{align}
Thus, 
\beq \label{7.3}
\int_{\Om} F_A(Du) \cdot D\psi\, =\, 0
\eeq
and since this holds for all $\psi$, we obtain \eqref{7.2}.

\ms

\noi \textbf{The problem.}

{\center{

\fbox{\parbox[pos]{260pt}{
\ms

The ``gap" in the previous reasoning is that \textit{we do not know whether $f$ is differentiable at zero, not even that $f$ is finite and not $\pm \infty$ near zero.} Moreover, \emph{\eqref{7.3} does not imply \eqref{7.2} unless $u \in C^2(\Om)$.} 

\ms
}}

\ms\ms }}

\noi Although under standard coercivity and finiteness assumptions\footnote{The simplest case is to assume (on top of the convexity of $F$) that $F(A)\geq C|A|^q-\frac{1}{C}$ for some $C>0$ and $q>1$. Then, for any $b\in W^{1,q}(\Om)$ such that $E(b,\Om)<\infty$, $E$ attain its infimum over $W^{1,q}_b(\Om)$.} it can be achieved that $E$ has a minimum point for which $|E(u,\Om)|<\infty$ and also we can \textbf{define weak solutions of \eqref{7.2} as those $u$ satisfying \eqref{7.3} for all $\psi$}, one of the most fundamental problems in Calculus of Variations is that

{\center{

\fbox{\parbox[pos]{260pt}{
\ms

\textit{there is still a gap between the assumptions guaranteeing existence of minima and the assumptions}

\ \  \textit{guaranteeing the existence of weak solutions !}

\ms
}}

\ms\ms }}

\ms

\noi The obstruction is that \eqref{7.1} may not be directionally (Gateaux) differentiable at the minimiser $u$. Moreover, mimima of \eqref{7.1} exist even when $F$ is not even in $C^1(\R^n)$, is which case \textit{$F_A(Du)$ may not even be Lebesgue measurable!} Hence, the classical technique of variations recalled above does not apply and there is no way to infer that a function $u$ which minimises the functional is a weak solution of \eqref{7.2}, which of course means that $u$ would satisfy
\beq   \label{7.4}
\int_\Om F_A(Du) \cdot D\psi =0,\ \ \ \text{ for all } \psi \in C^1_c(\Om).
\eeq

\ms

\noi \textbf{Back to Viscosity Solutions.} 

\ms

\noi In this chapter we show that  \textit{continuous local minimisers of the functional \eqref{7.1} in the space $W^{1,1}_{\text{loc}}(\Om)$ are Viscosity Solutions of the Euler-Lagrange PDE \eqref{7.2} corresponding to \eqref{7.1}, but expanded (with the divergence distributed)}:
\beq  \label{7.5}
F_{AA}(Du):D^2u \, =\, 0.
\eeq
Of course, $F_{AA}$ above denotes the Hessian matrix of $F$.  By \textbf{Local Minimisers of $E$ in $W^{1,1}_{\text{loc}}(\Om)$}, we mean functions $u \in W^{1,1}_{\text{loc}}(\Om)$ for which 
\beq \label{7.6}
E(u,\Om')\, \leq \, E(u+\psi,\Om'),\ \ \text{ for }\psi \in W^{1,1}_0(\Om'), \ \, \Om' \Subset \Om.
\eeq
As an application, we extend the classical theorem of Calculus of Variations regarding existence of solution to the Dirichlet problem
\beq   \label{7.7}
\left\{
\begin{array}{r}
F_{AA}(Du):D^2u \, =\, 0, \ \ \text{ in }\Om,\ \, \ms\\
 u\, =\, b, \ \text{ on }\p \Om,
\end{array}
\right.
\eeq
when $\Om \Subset \R^n$ and $b \in W^{1,1}_{\text{loc}}(\Om) \cap C^0(\overline{\Om})$ has finite energy on $\Om$ (that is $E(b,\Om)<+\infty$). These assumptions are much weaker than those guaranteeing the existence of weak solutions.

\ms

\noi \textbf{Theorem 1 (Minimisers as Viscosity Solutions of the E.-L. PDE).}  

\noi \emph{Fix $\Om\sub \R^n$ open and let $F\in C^2(\R^n)$ be convex. Then, continuous local minimisers of the functional \eqref{7.1} in $W^{1,1}_{\text{loc}}(\Om)$ are Viscosity Solutions of the expanded Euler-Lagrange PDE \eqref{7.7} on $\Om$.
}

\ms

\noi \textbf{Remark 2.}

{\center{

\fbox{\parbox[pos]{330pt}{
\ms

ALTHOUGH THE EULER-LAGRANGE PDE IS IN DIVERGENCE FORM AND $F$ IS SMOOTH CONVEX, WEAK SOLUTIONS IN $W^{1,r}(\Om)$ MAY NOT EXIST, BECAUSE UNLESS
\[
F_A(p)\, \leq\, C|p|^{r-1}
\]
THE FUNCTIONAL $E$ MAY NOT BE DIRECTIONALLY (GATEAUX) DIFFERENTIABLE AT THE MINIMISER ! ! !
\ms
}}

\ms\ms }}

\ms

\noi Indeed, for minimisers $u \in W^{1,r}(\Om)$ $r\geq 1$, assuming that such a bound on $F$ holds, H\"older's inequality implies 
\begin{align}
\left| \int_\Om F_A(Du) \cdot D\psi \right|\, &\leq \,  \int_\Om \big|F_A(Du)\big| \, |D\psi|  \nonumber\\
&\leq \, C\int_\Om |Du|^{r-1} \, |D\psi|  \nonumber\\
&\leq \, C \left(\int_\Om \big(|Du|^{r-1}\big)^{\frac{r}{r-1}}\right)^{\frac{r-1}{r}} \left( \int_\Om  |D\psi|^r \right)^r  \nonumber\\
&= \, C \left(\int_\Om |Du|^{r}\right)^{\frac{r-1}{r}} \left( \int_\Om  |D\psi|^r \right)^r \nonumber
\end{align} 
and hence $E$ is directionally differentiable at the minimiser $u$ since $F_A(Du) \cdot D\psi \in L^1(\Om)$.

\ms

\noi \textbf{Proof of Theorem 1.} Assume for the sake of contradiction that there is a $u \in(C^0\cap W^{1,1}_{\text{loc}})(\Om)$ which satisifes \eqref{7.6}, but $u$ is not a Viscosity Subsolution of \eqref{7.5}. Then, there exists $x\in \Om$, a smooth $\psi \in C^2(\R^n)$ and an $r>0$ such that
\[
u-\psi \,<\, 0\, = \,(u-\psi)(x)
\]
on $\mB_r(x)\set \{x\}$, while for some $c>0$ we have
\[
F_{AA}(D\psi(x)):D^2\psi(x) \, \leq \,-2c\, <\, 0.
\]
Since the function $F(D\psi)$ is in $C^2\big(\mB_r(x)\big)$, we have
\beq
F_{AA}(D\psi):D^2\psi\, =\, \Div\big(F_A(D\psi)\big)
\eeq
on $\mB_r(x)$. By restricting $r$ even further, we get
\beq  \label{7.9}
-\Div\big(F_A(D\psi)\big)\, \geq c,
\eeq
on $\mB_r(x)$. By strictness of the maximum of $u-\psi$,  there is a $k>0$ small such that, by sliding $\psi$ downwards to some $\psi-k$, we have
\[
\Om^+\, :=\, \Big\{u-\psi+k\, >\, 0 \Big\} \, \sub \mB_r(x)
\]
and also 
\[
u\, =\, \psi-k \, \text{ on }\, \p \Om^+. 
\]

\[
\includegraphics[scale=0.26]{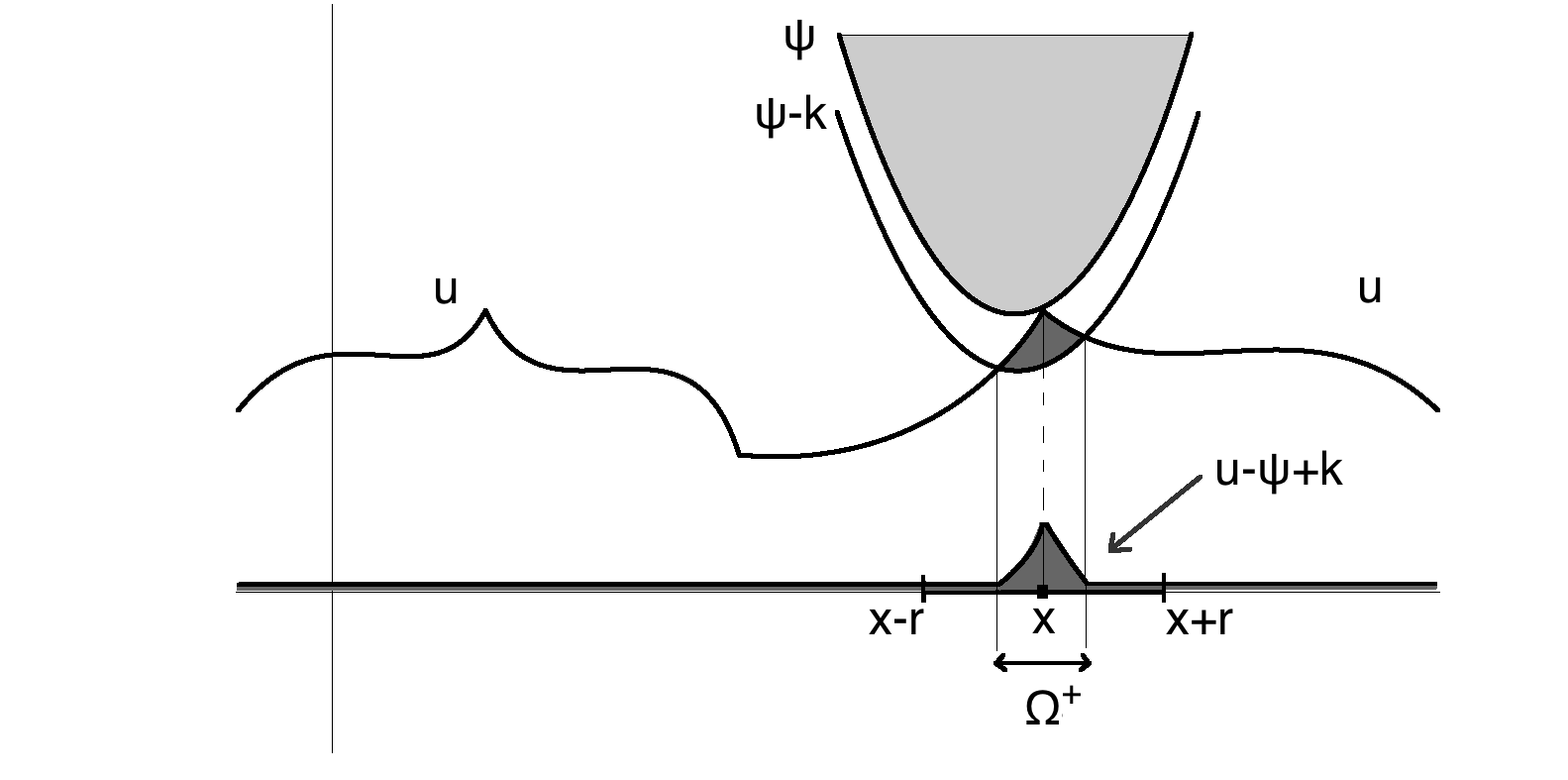}
\]
By multiplying \eqref{7.9} by the function 
\[
u-\psi+k \ \in \, W^{1,1}_0(\Om^+)
\]
and integrating by parts, we obtain
\beq \label{7.9}
\int_{\Om^+}F_A(D\psi) \cdot (Du-D\psi) \, \geq\, c\int_{\Om^+}\big|u-\psi+k \big|.
\eeq
Since $F$ is convex on $\R^n$, we have the elementary inequality
\beq \label{7.10}
F_A(a) \cdot (b-a) \, \leq\, F(b)-F(a).
\eeq
Indeed, let us derive \eqref{7.10} for the sake of completeness. By convexity of $F$, for any $\la \in (0,1]$ and $a,b \in \R^n$, we have
\[
F\big(\la b +(1-\la)a \big)\, \leq \, \la F(b) \, +\, (1-\la)F(a).
\]
Hence,
\[
\frac{1}{\la}\Big( F\big(a+\la(b-a) \big) \, -\, F(a) \Big)\, \leq \, F(b) \, -\, F(a)
\]
and by letting $\la \ri 0^+$, \eqref{7.10} follows. Going back to \eqref{7.9}, we see that by utilising \eqref{7.10}, we have
\begin{align} 
c\int_{\Om^+}|u-\psi+k| \, &\leq\, \int_{\Om^+}F(Du) \,- \, \int_{\Om^+}F(D(\psi-k)) \nonumber\\
     &=\, E\big(u,\Om^+\big) \, -\, E\big(\psi-k,\Om^+\big).  \nonumber
\end{align}
Hence, we obtain that 
\[
E\big(u,\Om^+\big) \,-\, E\big(\psi-k,\Om^+\big)\, \leq \, 0,
\]
while $u=\psi-k$ on $\p \Om^+$. Since the set $\Om^+$ is open, we deduce that
\[
\Om^+\, =\, \emptyset, 
\]
which is a contradiction. Hence, $u$ is a Viscosity Solution of 
\[
F_{AA}(Du):D^2u \,\geq\, 0
\]
on $\Om$. The supersolution property follows in the similar way, and so does the theorem.
\qed

\ms

The above result implies the following existence theorem:

\ms

\noi \textbf{Theorem 3 (Existence for the Dirichlet Problem).} \emph{Assume that the convex function $F \in C^2(\R^n)$ satisfies
\beq   \label{7.12}
 F(A)\, \geq \, C_1 |A|^s -C_2,
\eeq
for some $C_1,C_2>0$ and some $s>n$. Let also $\Om \sub \R^n$ be bounded and choose
\[
b \in  W^{1,1}_{\text{loc}}(\Om)\cap C^0(\overline{\Om})
\]
which has finite energy on $\Om$ i.e.\ $E(b,\Om)<+\infty$.\ms
}

\emph{Then, the Dirichlet Problem \eqref{7.7} has a Viscosity Solution, which is (globally) minimising for $E$ in $W^{1,s}_b(\Om)$.
}

\ms

\noi \textbf{Proof.} The argument is a simple implementation of the \emph{direct method of Calculus of Variations}, which we include it for the sake of completeness. Since $E(b,\Om)<+\infty$, if $|\cdot|$ denotes the Lebesgue measure, it follows that 
\begin{align}
C_1 \int_\Om |Db|^s \, &\leq \, C_2|\Om|\, + \, \int_\Om F(Db)  \nonumber \\
& \leq\, C_2|\Om|\, +\, E(b,\Om) \nonumber
\end{align}
which gives
\[
\|Db \|_{L^s(\Om)}\, <\, \infty.
\]
This implies that the infimum of $E$ in the affine space $W^{1,s}_b(\Om)$ is finite:
\[
-\infty\, <\, -C_2|\Om|\, \leq\, \inf_{W^{1,s}_b(\Om)}E\, \leq E(b,\Om)\, <\, \infty.
\]
Let $(u^m)_1^\infty$ be a minimising sequence. This means that
\[
E(u^m,\Om) \larrow  \inf_{W^{1,s}_b(\Om)}E
\]
as $m\ri \infty$. Let us recall the fundamental Poincar\'e inequality 
\[
 \big\|u_m-b  \big\|_{L^s(\Om)}\, \leq \, C \big\|Du_m-Db  \big\|_{L^s(\Om)},
\]
which applies to $u_m-b \in W^{1,r}_0(\Om)$. This implies
\begin{align}
\|u_m\|_{L^s(\Om)}\, &\leq  \,  \big\|u_m-b  \big\|_{L^s(\Om)}   \, +\, \|b \|_{L^s(\Om)}  \nonumber\\
&\leq \, C  \big\|Du_m-Db  \big\|_{L^s(\Om)}   \, +\, \|b \|_{L^s(\Om)}  \nonumber\\
&\leq \, C  \|Du_m \|_{L^s(\Om)}   \, +\, \|b \|_{W^{1,s}(\Om)}.  \nonumber
\end{align}
By Poincar\'e inequality, we have that $b\in W^{1,s}(\Om)$. Hence, we have the uniform bound 
\[
\|u^m\|_{W^{1,s}(\Om)}\, \leq \, C. 
\]
By \textit{weak compactness}, there exists a subsequence along which we have $u^m \rightharpoonup u$ as $m\ri \infty$ weakly in $W^{1,s}(\Om)$. In view of our assumption of $F$ and a well-know result in Calculus of Variations, the functional \eqref{7.1} is \textit{weakly lower-semicontinuous  in $W^{1,s}(\Om)$}, that is, if
\[
u^m \lharpoonup u \ \ \text{ weakly in $W^{1,s}(\Om)$ \  as }m \ri \infty,
\]
then
\[
E(u,\Om)\, \leq \, \liminf_{m\ri \infty}\, E(u^m,\Om).
\]
Since
\[
\liminf_{m\ri \infty}\, E(u^m,\Om)\, =\, \inf_{W^{1,s}_b(\Om)}E\,  <\, +\infty,
\]
we deduce that $u$ is minimiser of $E$ in $W_b^{1,s}(\Om)$. This means that
\[
E(u,\Om)\, \leq \, E(u+\psi, \Om)
\]
for all $\psi \in W^{1,s}_0(\Om)$. Since $s>n$, By Morrey's estimate, we have
\[
\big|u(x)-u(y) \big|\, \leq \, C |x-y|^{1-\frac{n}{s}}
\]
for $x,y\in \Om$ and hence (the ``precise representative" of) $u$ is in $C^0(\Om)$. By Theorem 1, $u$ is a Viscosity Solution of the PDE and also $u=b$ on $\p \Om$. Hence $u$ solves the Dirichlet problem and the theorem follows.     \qed

\ms

\ms

\ms

\noi \textbf{Remarks on Chapter 7.} For background material on Sobolev spaces, weak derivatives, the Calculus of Variations and the Euler-Lagrange equation, the reader may consult the textbooks of Evans \cite{E4} and Dacorogna \cite{D3}. The idea of the proof of Theorem 1 is taken from Barron-Jensen \cite{BJ}. The equivalence of weak solutions and viscosity solutions is an important and largely unexplored topic in PDE theory. For the model case of the $p$-Laplacian the reader may consult Juutinen-Lindqvist-Manfredi \cite{JLM}.

\chapter[Existence for the $\infty$-Laplacian]{Existence of Viscosity Solutions to the Dirichlet Problem for the $\infty$-Laplacian} 

In this chapter we prove existence of a Viscosity Solution $u\in C^0(\Om)$ to the Dirichlet problem for the $\infty$-Laplacian
\beq \label{8.1}
\left\{
\begin{array}{l}
\De_\infty u \, =\, 0,\ \ \text{ in }\Om, \ms\\
\hspace{20pt} u\, =\, b,\ \text{ on }\p \Om,
\end{array}
\right.
\eeq
where $\Om \Subset \R^n$ has Lipschitz boundary and $b\in W^{1,\infty}(\Om)$\footnote{It is well-known that under these assumptions on $\Om$ we have $b\in C^0(\overline{\Om})$. For simplicity, we shall require the boundary condition to be satisfied in the classical pointwise sense. However, regularity of $\p \Om$ is \emph{not} actually needed, once we interpret the boundary condition in the sense $u-b \in W^{1,\infty}_0(\Om)$}. Throughout this chapter, as in the previous one, we assume a basic degree of familiarity with the basic notions of measure theory and of Sobolev spaces. We recall from Chapter 1 that the $\infty$-Laplacian is the analogue of the ``Euler-Lagrange" PDE for the supremal functional
\[
E_\infty(u,\Om)\, :=\, \|Du\|_{L^\infty(\Om)}
\]
in $W^{1,\infty}(\Om)$. The $\infty$-Laplace operator is given by
\[
\De_\infty u \, :=\, Du \ot Du :D^2u\, =\, D_i u\, D_j u\, D^2_{ij}u
\]
and is a quasilinear non-divergence structure differential operator, which is degenerate elliptic. Indeed, by setting
\[
F_\infty(p,X)\, :=\, X:p\ot p
\]
the $\infty$-Laplacian PDE
\beq  \label{8.2}
\De_\infty u\, =\,0
\eeq
can be written as
\[
F_\infty(Du,D^2u)\, =\, 0
\]
and we have
\[
X\leq Y\text{ in }\mS(n)\ \ \Longrightarrow\ \ F_\infty(p,X)\, \leq F_\infty(p,Y).
\]
The PDE \eqref{8.2} is in non-divergence form and weak solutions can not be defined. In fact, as we observed in Chapter 1, \eqref{8.2} does not have solution is the classical senses. We prove existence of a Viscosity Solution to \eqref{8.1} by utilising the \textit{stability results of Viscosity Solutions} of Chapter 3. 

As in Chapter 5, by \emph{Solution to the Dirichlet problem \eqref{8.1}} we mean a function $u_\infty \in C^0(\overline{\Om})$ which satisfies $u_\infty=b$ on $\p \Om$ and for any $x\in \Om$, we have
\beq   \label{8.3}
\inf_{(p,X)\in \J^{2,+}u_\infty(x)}\, X:p \ot p\, \geq \, 0,
\eeq
and also
\beq   \label{8.4}
\sup_{(p,X)\in \J^{2,-}u_\infty(x)}\, X:p \ot p\, \leq \, 0. 
\eeq
Hence, the viscosity inequalities \eqref{8.3},  \eqref{8.4} give a rigorous meaning to the statement \emph{$u_\infty$ is an $\infty$-Harmonic function on $\Om$.}

The simplest way to prove existence of solution to \eqref{8.1} is \textit{by approximation}. It seems that the natural choice of approximating sequences of PDE and respective solutions is that of \textit{$p$-Harmonic functions} as $p\ri \infty$:
\beq    \label{8.5}
\De_p u \, :=\, \Div\big(|Du|^{p-2}Du \big)\, =\, 0.
\eeq

\ms

\noi \textbf{Formal derivation of the $\De_\infty$ from $\De_p$ as $p\ri \infty$.} \textit{Before proceeding to the rigorous arguments, let us derive the $\infty$-Laplacian from the $p$-Laplacian as $p \ri \infty$. By distributing derivatives in \eqref{8.5}, we have
\[
D_i\big(|Du|^{p-2}D_iu \big)\ =\ |Du|^{p-2}D^2_{ii}u\ +\ (p-2)|Du|^{p-4}D_i u\, D_j u\, D^2_{ij}u.
\]
Hence,
\[
|Du|^{p-2}\De u\ +\  (p-2)|Du|^{p-4}Du \ot Du :D^2u\ = \ 0.
\]
We normalise by multiplying by 
\[
(p-2)|Du|^{p-4}, 
\]
to obtain
\[
Du \ot Du :D^2u\ +\ \frac{|Du|^{2}}{p-2}\De u \ = \ 0.
\]
Then, as $p\ri \infty$ we formally obtain
\[
Du \ot Du :D^2u \ = \ 0.
\]}

\ms

The idea is to prove existence by making the previous formal derivation \emph{rigorous}. An extra gift of this method is that the solution we obtain is not just continuous, but also Lipschitz continuous. Accordingly, we have

\ms

\noi \textbf{Theorem 1 (Existence of $\infty$-harmonic functions with prescribed boundary values).}

\noi \emph{For any $\Om \Subset \R^n$ with Lipschitz boundary and $b\in W^{1,\infty}(\Om)$, the Dirichlet problem
\[
\left\{
\begin{array}{l}
\De_\infty u \, =\, 0,\ \ \text{ in }\Om, \ms\\
\hspace{20pt} u\, =\, b,\ \text{ on }\p \Om,
\end{array}
\right.
\]
has a Viscosity solution $u_\infty \in W^{1,\infty}(\Om)$ (which by definition satisfies \eqref{8.3}, \eqref{8.4} in $\Om$).}

\ms

For the proof we need two lemmas. Roughly, these lemmas say

{\center{

\fbox{\parbox[pos]{240pt}{
\ms

\ \ The Dirichlet problem for the $m$-Laplacian 
\[
\left\{
\begin{array}{l}
\De_m u \, =\, 0,\ \ \text{ in }\Om, \ms\\
\hspace{20pt} u\, =\, b,\ \text{ on }\p \Om,
\end{array}
\right.
\]
\ \ has a \textbf{Viscosity Solution} $u_m \in W^{1,m}_b(\Om)$,

\ms
}}

\ms\ms }}

and

{\center{

\fbox{\parbox[pos]{250pt}{
\ms

\ \ The family of $m$-Harmonic functions $\{u_m\}_{m=1}^\infty$ \ms

\hspace{55pt} is \textbf{precompact} in $C^0(\overline{\Om})$.

\ms
}}

\ms\ms }}

\noi The rest of the proof is a simple application of the stability of Viscosity Solutions under locally uniform limits. A subtle point is that the \emph{$m$-Laplacian is divergence PDE and actually the Euler-lagrange PDE of the $m$-Dirichlet functional
\beq \label{8.6}
E_m(u,\Om)\, :=\, \int_{\Om}|Du|^m.
\eeq
Thus, the solutions that can be constructed by minimising the functional \eqref{8.6} are \textbf{weak solutions}.} Happily, Theorem 3 of the previous chapter allows to deduce that actually

{\center{

\fbox{\parbox[pos]{260pt}{
\ms

\textit{ The minimisers of the $m$-Dirichlet functional are} 

\ms

\textit{ Viscosity Solutions of the $m$-Laplacian expanded!}

\ms
}}

\ms\ms }}

\noi For, we have

\ms

\noi \textbf{Lemma 2 (Existence of $m$-Harmonic functions in the Viscosity Sense).}

\noi \emph{Fix $\Om \Subset \R^n$, $b\in W^{1,\infty}(\Om)$ and $m> n\geq 2$. Consider the Dirichet problem for the $m$-Laplacian expanded
\[
\left\{
\begin{array}{r}
|Du|^{m-2}\De u\ +\  (m-2)|Du|^{m-4}Du \ot Du :D^2u\, = \, 0,\ \ \text{ in }\Om, \  \ms\\
u\, =\, b,\ \text{ on }\p \Om. 
\end{array}
\right.
\]
Then, the problem has a Viscosity Solution $u_m\in W^{1,m}_b(\Om)$, which minimises the $m$-Dirichlet functional \eqref{8.6} in $W^{1,m}(\Om)$.}

\ms

\noi \textbf{Proof.} Consider the minimisation problem
\[
E_m(u_m,\Om) \ =\ \inf_{v \in W^{1,\infty}(\Om)} \, E_m(v,\Om) 
\]
where
\[
E_m(u,\Om)\, =\, \int_{\Om}f(Du)
\]
and
\[
f(p)\, := \, |p|^m.
\]
Since $m>n$, by Theorem 3 of the previous chapter, $f \in C^2(\R^n)$ is convex and satisfies the assumptions of the statement for $s\equiv m$. Hence, there exists a viscosity solution of the PDE
\[
f_{AA}(Du) : D^2u\, =\, 0
\]
with $u=b$ on $\p \Om$. By differentiating $f$ we have 
\[
f_{A_i}(p)\, =\, m|p|^{m-2}p_i, 
\]
and hence
\[
f_{A_i A_j}(p)\ =\ m|p|^{m-2}\de_{ij}\ +\ m(m-2)|p|^{m-4}p_i \, p_j .
\]
Hence, $u_m \in W^{1,m}_b(\Om)$ is a Viscosity solution to the Dirichlet problem and the Lemma follows.      \qed

\ms

An intermediate technical step we need is packed in the next

\ms

\noi \textbf{Lemma 3 (Normalisation).} \emph{Let $\Om \sub \R^n$, $m>n \geq 2$ and $u\in C^0(\Om)$. Then, the PDE
\beq \label{8.7}
|Du|^{m-2}\De u\ +\  (m-2)|Du|^{m-4}Du \ot Du :D^2u\, = \, 0
\eeq
and the PDE
\beq  \label{8.8}
Du \ot Du :D^2u\ +\ \frac{|Du|^{2}}{m-2}\De u \, = \, 0
\eeq
are equivalent in the Viscosity sense. Namely, $u$ is a Viscosity solution of \eqref{8.7} if and only if it is a Viscosity solution of \eqref{8.8}.}

\ms

\noi \textbf{Proof.} Let $(p,X)\in \J^{2,+}u(x)$ for some $x\in \Om$. Then, if $u$ solves  \eqref{8.7}, we have
\[
|p|^{m-2}X:I\ +\  (m-2)|p|^{m-2}\frac{p}{|p|} \ot \frac{p}{|p|}  :X\, \geq \, 0.
\]
If $p\neq 0$, the latter inequality is equivalent to
\[
p \ot p :X \ +\ \frac{|p|^{2}}{m-2}X:I \, \geq \, 0.
\]
On the other hand, if $p=0$ then both the above inequalities vanish since $m>n\geq 2$. The supersolution property follows similarly. Hence,  the lemma ensues.                  \qed

\ms

Next, we have

\ms

\noi \textbf{Lemma 4 (Compactness).} \emph{Let $\{u_m\}_n^\infty$ be the family of viscosity solutions to the Dirichlet problems
\[
\left\{
\begin{array}{r}
Du_m \ot Du_m :D^2u_m\ +\ \dfrac{|Du_m|^{2}}{m-2}\De u_m \, = \, 0,\ \ \text{ in }\Om, \ \,  \ms\\
u_m\, =\, b,\ \text{ on }\p \Om.
\end{array}
\right.
\]
There is a function $u_\infty \in W^{1,\infty}_b(\Om)$ and a sequence of $m's$, such that, along this sequence,
\[
u_m \larrow u_\infty
\]
in $C^0(\Om)$, as $m\ri \infty$.}

\ms

\noi \textbf{Proof.} Since each $u_m$ minimises the $m$-Dirichlet functional, we have
\[
E_m(u_m,\Om)\, \leq\, E_m(b,\Om)\, < \, \infty.
\]
Hence, by H\"older inequality,
\begin{align}
\|Du_m\|_{L^m(\Om)}\, & =\, E_m(u_m,\Om)^m  \nonumber\\
&\leq \, E_m(b,\Om)^m  \nonumber\\
&=\, \|Db\|_{L^m(\Om)}.  \nonumber
\end{align}
Fix a $k>n$. For any $m\geq k$, we have
\begin{align}
\|Du_m\|_{L^k(\Om)}\, &\leq \, \|Du_m\|_{L^m(\Om)} |\Om|^{1/k - 1/m} \nonumber\\
&\leq\, |\Om|^{1/k} \|Db\|_{L^\infty(\Om)}. \nonumber
\end{align}
Consequently,
\beq \label{8.10}
\|Du_m\|_{L^k(\Om)}\,  \leq \, |\Om|^{1/k} \|Db\|_{L^\infty(\Om)}.
\eeq
By Poincar\'e inequality, we have
\begin{align}
\|u_m\|_{L^k(\Om)}\, &\leq  \,  \big\|u_m-b  \big\|_{L^k(\Om)}   \, +\, \|b \|_{L^m(\Om)}  \nonumber\\
&\leq \, C  \big\|Du_m-Db  \big\|_{L^k(\Om)}   \, +\, \|b \|_{L^k(\Om)}  \nonumber\\
&\leq \, C  \|Du_m \|_{L^k(\Om)}   \, +\, \|b \|_{W^{1,k}(\Om)} \nonumber
\end{align}
and hence
\beq \label{8.11}
\|u_m\|_{L^k(\Om)}\, \leq  \, C  \Big(\|Du_m \|_{L^k(\Om)}   \, +\, \|b \|_{W^{1,\infty}(\Om)} \Big).
\eeq
By \eqref{8.10},  \eqref{8.11} we have that for each $k>n$ fixed,  we have the bounds
\[
\sup_{m>k}\, \|u_m\|_{W^{1,k}(\Om)}\, \leq  \, C.
\]
By weak compactness, there exists a $u_\infty \in W^{1,k}(\Om)$ such that, along a sequence,
\[
u_m \lharpoonup u_\infty, \, \text{ weakly in }W^{1,k}(\Om) \, \text{ as }m\ri \infty.
\]
This holds for all $k>n$. By a Cantor diagonal argument, we can extract a common subsequence of $m$'s such that, 
\[
u_m \lharpoonup u_\infty, \, \text{ weakly in }W^{1,k}(\Om) \, \text{ for all $k>n$, as }m\ri \infty.
\]
By Morrey's estimate, we also have that along this sequence
\[
u_m \larrow u_\infty \, \text{ in }C^0(\overline{\Om}) \, \text{ as }m\ri \infty.
\]
By letting $m\ri \infty$ in \eqref{8.10}, the sequential weak lower-semicontinuity of the $L^k$-norm, implies
\begin{align}
\|Du_\infty\|_{L^k(\Om)}\,  &\leq \,\underset{m\ri \infty}{\lim\inf}\, \|Du_m\|_{L^k(\Om)} \nonumber\\
&  \leq \, |\Om|^{1/k} \|Db\|_{L^\infty(\Om)}. \nonumber
\end{align}
By letting now $k\ri \infty$, since the $L^k$ norm convergence pointwise to the $L^\infty$ norm, we conclude that
\begin{align}
\|Du_\infty\|_{L^\infty(\Om)}\,  &\leq \,\underset{m\ri \infty}{\lim}\, \|Du_\infty\|_{L^k(\Om)} \nonumber\\
&\leq \,\underset{m\ri \infty}{\lim}\,  |\Om|^{1/k} \|Db\|_{L^\infty(\Om)}\nonumber\\
&  \leq \, \|Db\|_{L^\infty(\Om)}. \nonumber
\end{align}
Consequently, $u_\infty \in W^{1,\infty}(\Om)$ and the lemma ensues.     \qed

\ms

\noi \textbf{Proof of Theorem 1.} We set
\[
F_m(p,X)\, :=\, X:p\ot p\ +\ \frac{|p|^2}{m-2}X:I
\]
and recall that
\[
F_\infty(p,X)\, =\, X:p\ot p.
\]
We have that the expanded $m$-Laplacian can be written as
\[
F_m(Du_m,D^2u_m)\, = \, 0
\]
and the $\infty$-Laplacian as
\[
F_\infty(Du_\infty,D^2u_\infty)\, = \, 0.
\]
By Lemmas 2 and 3, we have that the Dirichlet problem
\[
\left\{
\begin{array}{r}
F_m(Du_m,D^2u_m)\, = \, 0,\ \ \text{ in }\Om, \ \,  \ms\\
u_m\, =\, b,\ \text{ on }\p \Om,
\end{array}
\right.
\]
has a Viscosity Solution. By Lemma 4, there is a $u_\infty \in C^0(\overline{\Om})$, such that, along a sequence of $m$'s we have
\[
u_m \larrow u_\infty
\]
uniformly on $\Om$ as $m\ri \infty$, while $u_\infty \in W^{1,\infty}_b(\Om)$. Since
\[
F_m \larrow F_\infty
\]
locally uniformly on $\R^n \by \mS(n)$ as $m\ri \infty$, it follows that $u_\infty$ solves
\[
\left\{
\begin{array}{r}
F_\infty(Du_\infty,D^2u_\infty)\, = \, 0,\ \ \text{ in }\Om, \ \, \ms\\
u_\infty\, =\, b,\ \text{ on }\p \Om,
\end{array}
\right.
\]
in the viscosity sense. The theorem follows. \qed

\ms

\ms

\ms

\noi \textbf{Remarks on Chapter 8.} The $\infty$-Laplacian has been discovered by Aronsson in \cite{A3} and was studied until the mid 1980s by (exclusively) the same author, when he observed the existence of ``singular solutions" in \cite{A6, A7}. Viscosity solutions (which had already been extended to the 2nd order case) found a very nice application to Calculus of Variations in $L^\infty$. Bhattacharya-DiBenedetto-Manfredi were the first to prove existence for the Dirichlet problem for the $\infty$-Laplacian in \cite{BDM}. Uniqueness was a big step forward to the theory and was first made by Jensen in \cite{J2}. For a more recent proof of this fact, which is based on properties of the $\infty$-Laplacian discovered later see Armstrong-Smart \cite{AS}. For further material on Calculus of Variations in $L^\infty$, the reader may consult Barron \cite{B}, Crandall \cite{C1} (see also Juutinen \cite{Ju}). In particular, lower semicontinuity properties of supremal functionals have been studied in Barron-Jensen-Wang \cite{BJW2}. One of the main open problems (to date!) is the question of $C^1$ regularity of $\infty$-harmonic maps, which is not know in higher than $2$ dimensions; see Savin \cite{S} and Evans-Savin \cite{ES} and Evans-Smart \cite{ESm}. A slightly more general PDE than the $\infty$-Laplacian is the so-called \emph{Aronsson equation}, which reads $A_\infty u:=H_P(Du)\ot H_P(Du):D^2u=0$ and corresponds to the functional $\|H(Du)\|_{L^\infty(\Om)}$ for a Hamiltonian $H\in C^2(\R^n)$. This PDE has a corresponding interesting existence-uniqueness theory, see  Barron-Jensen-Wang \cite{BJW1}, Armstrong-Crandall-Julin-Smart \cite{ACJS} and also Barron-Evans-Jensen \cite{BEJ}. The corresponding regularity problem is also open, see Wang-Yu \cite{WY} and also \cite{K}. The vector-valued theory for the $\infty$-Laplacian is very promising but still at its birth, see \cite{K1} and also Sheffield-Smart \cite{SS}.

\chapter[Miscellaneous and extensions]{Miscellaneous topics and some extensions of the theory}

In this chapter we collect, mostly without proofs, some very important facts related to the main text and to the applications we presented so far which we were not able to analyse to a wide extent due to space and time considerations.

\subsection*{Fundamental Solutions of the $\infty$-Laplacian}

Let $u\in C^0(\Om)$, $\Om\sub \R^n$. We have already seen that the $\infty$-Laplacian 
\[
\De_\infty u \, :=\, Du \ot Du :D^2u\, =\, 0
\]
is an all-important PDE, and the archetypal example of ``Euler-Lagrange" analogue, associated to the model supremal functional
\[
E_\infty (u,\Om)\, :=\, \|Du\|_{L^\infty(\Om)},
\]
placed in $W^{1,\infty}(\Om)$. In this section we record a set of equivalent characterisations of Viscosity Solutions to $\De_\infty u =0$, including a precise meaning to the exact $L^\infty$ variational principle governing the $\infty$-Laplacian. A \textbf{very important} point is that, roughly\footnote{This statement is meant in the sense that the subsolutions of $\De_p$ can be characterised via comparison with a class of ``fundamental" solutions only for $p=2,\infty$. See \cite{CZ1, CZ2}.},

{\center{

\fbox{\parbox[pos]{330pt}{
\ms

\textit{The Laplacian $\De=\De_2$ and the $\infty$-Laplacian $\De_\infty$ are the only operators of the $p$-Laplace class that have fundamental solutions!}

\ms
}}

\ms\ms }}

\noi A vague comment would be that

{\center{

\fbox{\parbox[pos]{300pt}{
\ms

\textit{Harmonic functions can be characterised via convolution and duality,} 

\ms
\centerline{while}  
\ms 

\textit{$\infty$-Harmonic functions can be characterised via sup-convolution and nonlinear duality.}

\ms
}}

\ms\ms }}

\noi A further slightly humorous comment is that 

{\center{

\fbox{\parbox[pos]{300pt}{
\ms

the $p$-Laplacian $\De_p$ ``does not exist", since it is a quasilinear 

\ \ \ \ combination of the Laplacian and the $\infty$-Laplacian:

\[
\De_p u\, =\, |Du|^{p-2}\De u\, +\, (p-2)|Du|^{p-4}\De_\infty. 
\]

\ms
}}

\ms\ms }}

\noi We now record the following result

\ms

\noi \textbf{Theorem 1 (Characterisations of $\De_\infty$).} Let $u\in C^0(\Om)$ and $\Om\sub \R^n$. 
The following statements are equivalent:

\ms

\noi (a) \textbf{(Viscosity Solutions of $\De_\infty$)}  \emph{$u$ is Viscosity Solution of the $\infty$-Laplacian, that is, for all $x\in \Om$,
\begin{align}
\inf_{(p,X)\in \J^{2,+}u_\infty(x)}\, X:p \ot p\, &\geq \, 0, \nonumber\\
\sup_{(p,X)\in \J^{2,-}u_\infty(x)}\, X:p \ot p\, &\leq \, 0,\nonumber 
\end{align}}
\ms
 
\noi (b) \textbf{(Absolute Minimisers of $E_\infty$)}  \emph{$u$ is locally Lipschitz continuous on $\Om$ (i.e.\ $u\in W_{\text{loc}}^{1,\infty}(\Om)$) and for all $\Om'\Subset \Om$ and all $\psi \in W_0^{1,\infty}(\Om')$, we have
\[
E_\infty (u,\Om')\, \leq\, E_\infty(u+\psi,\Om').
\]}

\ms
 
\noi (c) \textbf{(Absolutely Minimising Lipschitz extensions)}   \emph{$u$ is locally Lipschitz continuous on $\Om$ and for all $\Om'\Subset \Om$, we have
\[
\Lip(u,\Om')\, =\, \Lip(u,\p \Om')
\]
where ``Lip" is the Lipschitz constant functional:
\[
\Lip(u,A)\, :=\, \sup\left\{\frac{|u(x)-u(y)|}{|x-y|}\ \, : \ \, x,y\in A,\ x\neq y \right\}, \ \ A\sub \R^n.
\]}
\ms
 
\noi (d) \textbf{(Comparison with cones - ``Fundamental Solutions")} \emph{$u$ satisifes the comparison principle with respect to all cone functions 
\[
C_{x,L}(z)\, :=\, a\, +\, L|z-x|
\]
(where $a,L \in \R$) on subdomain for which the vertex $x$ is not included. Namely, for any cone $C_{x,L}$ and any $\Om' \sub \Om \set \{x\}$, we have
\begin{align}
\sup_{\Om'}\, \big(u- C_{x,L}\big) \, \leq \, \max_{\p \Om'} \,\big(u- C_{x,L}\big), \nonumber\\
\inf_{\Om'}\, \big(u- C_{x,L}\big) \, \leq \, \min_{\p \Om'} \,\big(u- C_{x,L}\big).\nonumber 
\end{align}}
\[
\includegraphics[scale=0.22]{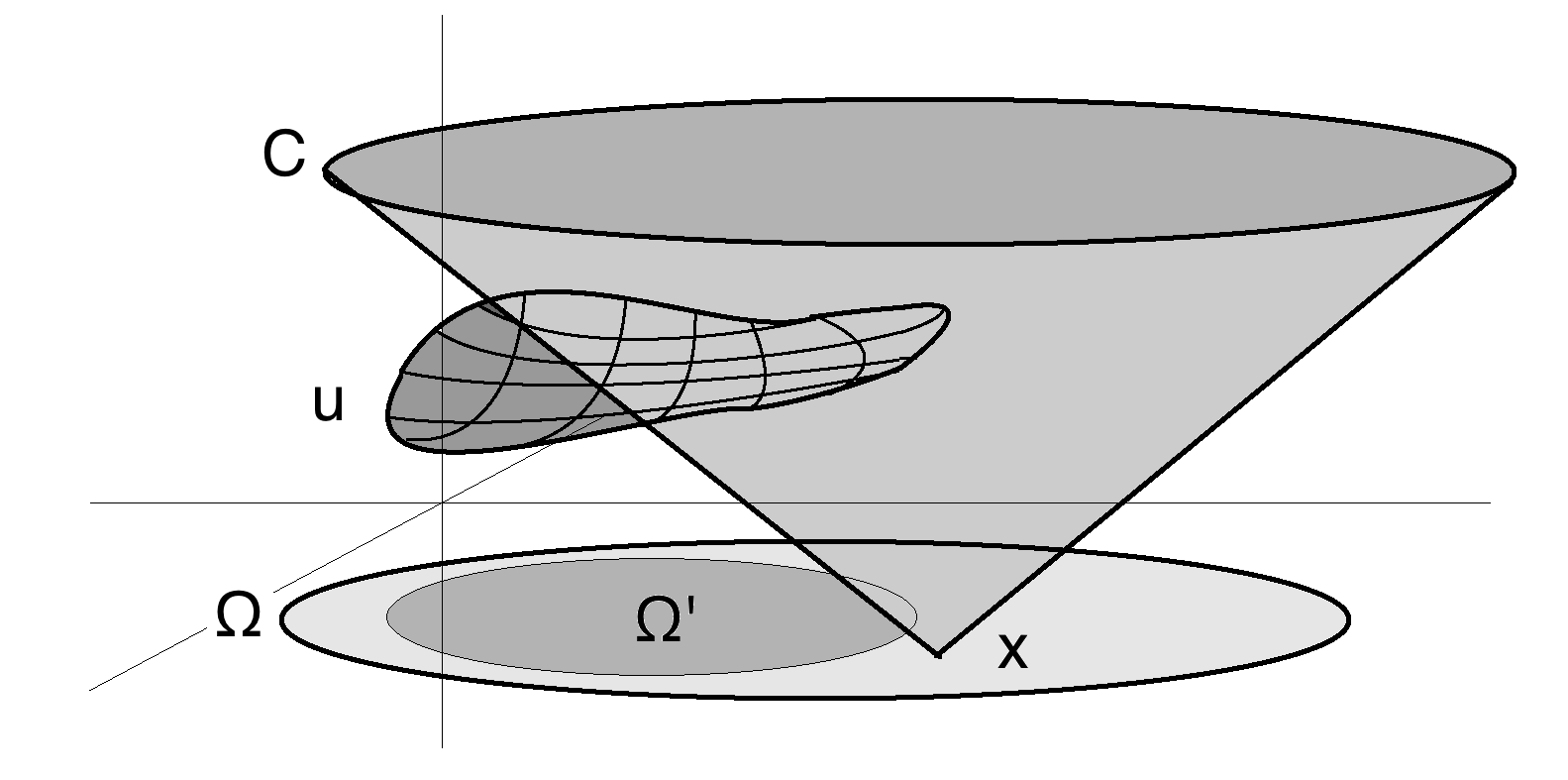}
\]

\subsection*{The $\infty$-Laplacian and Tug-of-War Differential Games}

A recent trend in PDE theory is the quest to find game-theoretic interpretations of differential equations. Except for the intrinsic interest in finding hidden connections between these two fields, there is also the underlying hope that new tools will be discovered in order to solve PDE. 

Such an approach as the one just described has been introduced by Peres-Schramm-Sheffield-Wilson. Given an open domain $\Om \sub \R^n$ and a fixed $\e>0$, we introduce a random two-person zero-sum game as follows: A token is placed at an arbitrary point $x_0 \in \Om$. Two players, Player I and Player II, flip a fair coin. The player who wins moves the token from $x_0$ to $x_1=x_0 +\e e$, where $|e|=1$. If Player I wins, he chooses the direction $e$ such that the expected payoff of the game is maximised, where if Player II wins, he chooses the direction $e$ such that the expected payoff of the game is minimised. The process is repeated and the token is moved from $x_1$ to $x_2=x_1 +\e e$ by the player who wins the 2nd coin flip. The game ends when for some $x_*=x_k$ we reach the boundary of $\Om$ (and exit $\Om$). The expected payoff then is $b(x_*)$, where $b: (\p \Om)^\e \ri \R$ is a function defined on a strip of width $\e$ around the boundary $\p \Om$. By employing the Dynamic Programming Principle, it follows that the payoff $u^\e$ satisfies the differential identity
\[
u^\e(x)\, =\, \frac{1}{2}\left( \max_{\p \mB_\e(x)}u^\e \, +\, \min_{\p \mB_\e(x)}u^\e\right)
\]
when starting at $x=x_0\in \Om$. The idea is that Player I selects the direction on the sphere which maximises, while Player II selects the direction which minimises, and both events may happen with equal probability $1/2$ since the coin is fair. It is a remarkable observation that if $u^\e \larrow u$ as $\e \ri 0$, then $u$ is $\infty$-Harmonic and satisfies $u=b$ on $\p \Om$:
\[
\left\{
\begin{array}{r}
\De_\infty u  = 0, \ \ \text{ in }\Om, \\
u =b, \text{ on }\p \Om .
\end{array}
\right.
\]
We give some \textbf{formal} heuristics of this derivation. If $L : \R^n \larrow \R$ is a linear non-constant function, then
\begin{align}
\max_{\p \mB_\e(x)}\, L\, =\, L\left(x +\e\frac{DL(x)}{|DL(x)|} \right), \nonumber\\
\min_{\p \mB_\e(x)}\, L\, =\, L\left(x -\e\frac{DL(x)}{|DL(x)|} \right).  \nonumber
\end{align}
Hence, for $\e>0$ small, if $Du(x)\neq 0$,
\begin{align}
\max_{\p \mB_\e(x)}\, u\, \cong\, u \left(x +\e\frac{Du (x)}{|Du (x)|} \right), \nonumber\\
\min_{\p \mB_\e(x)}\, u \, \cong \, u \left(x -\e\frac{Du (x)}{|Du (x)|} \right).  \nonumber
\end{align}
Consequently, by Taylor expansion,
\begin{align}
\max_{\p \mB_\e(x)}\, u \, \cong\, u (x)  + \e|Du (x)| + \frac{\e^2}{2}\frac{Du (x)}{|Du (x)|}\ot \frac{Du (x)}{|Du(x)|}:D^2u (x) +o(\e^2), \nonumber\\
\min_{\p \mB_\e(x)}\, u  \, \cong \, u (x)  - \e|Du(x)| + \frac{\e^2}{2}\frac{Du (x)}{|Du (x)|}\ot \frac{Du  (x)}{|Du (x)|}:D^2u (x) +o(\e^2).  \nonumber
\end{align}
Hence, we reach to the identity
\[
\frac{Du(x) \ot Du(x)}{|Du(x)|^2} :D^2u(x)\, =\, o(1)\, +\, \frac{1}{\e^2}\left(\max_{\p \mB_\e(x)}\, u\, +\, \min_{\p \mB_\e(x)}\, u\, -\, 2u(x)\right)
\]
which, at least formally, gives the connection between Tug-of-War games and the $\infty$-Laplacian. The above derivation can be made rigorous in the context of Viscosity Solutions.

\subsection*{Discontinuous coefficients, Discontinuous solutions}

In these notes, we have defined continuous Viscosity Solutions of PDE defined by continuous nonlinearities. However, neither the solutions nor the ``coefficients" is actually needed to be continuous. It turns out that all the results hold true if the ``solution function"  and the ``coefficient function" are replaced by its semi-continuous envelopes. Accordingly, we have the next

\ms

\noi \textbf{Definition 1 (Discontinuous Viscosity Solutions of discontinuous PDE).} \emph{Consider the PDE
\[
F\big(\cdot,u,Du,D^2u\big)\, =\, 0,
\]
where
\[
F\ :\ \Om\by\R \by\R^n \by\mS(n) \larrow \R
\]
is degenerate elliptic, namely satisfies
\[
X\leq Y\ \text{ in }\, \mS(n)\ \ \Longrightarrow \ \ F(x,r,p,X)\leq F(x,r,p,Y).
\]}
(a) \emph{ The function $u : \Om \sub \R^n \ri \R$ is a Viscosity Subsolution of the PDE (or Viscosity Solution of $F\big(\cdot,u,Du,D^2u\big)\geq 0$) on $\Om$, when 
\beq \label{2.8}
(p,X)\in \J^{2,+}u^*(x)\ \ \Longrightarrow \ \ F^*\big(x,u^*(x),p,X \big)\geq 0,
\eeq
for all $x\in \Om$. Here, $u^*$ and $F^*$ denote the upper-semicontinuous envelopes of $u$ and $F$ respectively.}

(b) \emph{ The function $u : \Om \sub \R^n \ri \R$  is a Viscosity Supersolution of the PDE (or Viscosity Solution of $F\big(\cdot,u,Du,D^2u\big)\leq 0$) on $\Om$, when 
\beq  \label{2.9}
(p,X)\in \J^{2,-}u_*(x)\ \ \Longrightarrow \ \ F_*\big(x,u_*(x),p,X\big)\leq 0,
\eeq
for all $x\in \Om$. Here, $u_*$ and $F_*$ denote the lower-semicontinuous envelopes of $u$ and $F$ respectively.}

(c) \emph{ The function $u : \Om \sub \R^n \ri \R$  is a Viscosity Solution of the PDE when it is both a Viscosity Subsolution and a Viscosity Supersolution.}

\subsection*{Barles-Perthame Relaxed Limits (1-sided Uniform Convergence) and generalised 1-sided stability}

In Chapter 3 we saw that Viscosity Sub-/Super- Solutions pass to limits under $\pm \Gamma$-convergence, a very weak notion of 1-sided convergence. A particular case of this kind of convergence is described by the so-called ``Barles-Perthame relaxed limits".

\ms

\noi \textbf{Definition 1 (Barles-Perthame relaxed limits).} \emph{Let $\{u_m\}_1^\infty : \Om \sub \R^n \larrow \R$ be a sequence of functions. We set
\begin{align} \label{9.1}
\overline{u}(x)\, &:=\, \underset{m\ri \infty}{\lim{\sup}^*}\, u_m(x) \\
 &\equiv\, \underset{j\ri \infty}{\lim}\, \sup_{m> j}\, \sup_{|y-x|<1/j} u_m(y), \nonumber
\end{align}
and 
\begin{align} \label{9.2}
\underline{u}(x)\, &:=\, \underset{m\ri \infty}{\lim{\inf}_*}\, u_m(x) \\
 &\equiv\, \underset{j\ri \infty}{\lim}\, \inf_{m> j}\, \inf_{|y-x|<1/j} u_m(y). \nonumber
\end{align}}

\ms

\noi \textbf{Remark 2.} \textit{It is a simple exercise that  the function $\overline{u}$, $\underline{u}$ are upper and lower semi-continuous respectively:
\[
\overline{u} \, \in \, USC(\Om),\ \ \  \underline{u} \, \in \, LSC(\Om).
\]}

A basic property of these limits is, as the reader can easily check, that they imply $\pm \Gamma$-convergence. This gives

\ms

\noi \textbf{Lemma 3 (Stability of Jets under relaxed limits).} \emph{Let $\Om\sub \R^n$ and suppose $\{u_m\}_1^\infty \sub USC(\Om)$. Let $\overline{u}$ be given by \eqref{9.1}. If $x\in \Om$ and $(p,X)\in \J^{2,+}\overline{u}(x)$, then there exist $\{x_m\}_1^\infty \sub \Om$ and 
\[
(p_m,X_m)\in \J^{2,+}u_m(x_m)
\]
such that 
\[
(x_m,p_m,X_m) \larrow (x,p,X)
\]
as $m\ri \infty$.
}

\ms

A symmetric statement holds for $\underline{u}$ and subjets $\J^{2,-}\underline{u}$. An even more surprising statement is that if each $u_m$ is a subsolution of a PDE $F_m= 0$, then $\overline{u}$ is a subsolution of $\overline{F}$:

\ms

\noi \textbf{Theorem 4 (Relaxed limits as subsolutions).} \emph{Let $\Om\sub \R^n$ and suppose $\{u_m\}_1^\infty \sub USC(\Om)$. Consider also functions (not necessarily continuous!)
\[
F_m\ :\ \Om \by \R \by \R^n \by \mS(n) \larrow \R
\]
which are proper, that is, they satisfy the monotonicity assumptions
\begin{align}
X\leq Y\ \text{ in }\mS(n)\ \ &\Longrightarrow \ \ F(x,r,p,X)\leq F(x,r,p,Y), \nonumber\\
r\leq s\ \ \text{ in }\R\ \ \ \ \ \, &\Longrightarrow \ \ F(x,r,p,X)\geq F(x,s,p,X),  \nonumber
\end{align}
for $X,Y \in \mS(n)$, $r,s \in \R$ and $(x,p)\in \Om \by \R^n$. Let $\overline{u}$ be given by \eqref{9.1} and suppose $\overline{u}(x)<\infty$ for all $x\in \Om$. If each $u_m$ is a Viscosity Solution of
\[
F_m\big (\cdot, u_m, Du_m,D^2u_m \big)\ \geq \, 0
\]
on $\Om$, then $\overline{u}$ solves
\[
\overline{F} \big(\cdot, \overline{u}, D\overline{u},D^2\overline{u} \big)\ \geq \, 0
\]
on $\Om$, where
\[
\overline{F}(x,r,p,X)\, :=\, \underset{m\ri \infty}{\lim{\sup}^*}\, F_m(x,r,p,X).
\]
}

\ms

\noi The proof is near-trivial and hence omitted (see Chapter 3).

\ms

One basic property of the relaxed limits is that the are ``1-sided locally uniform convergence:"

\ms

\noi \textbf{Lemma 5 (Relaxed limits as 1-sided $C^0$ convergence).} \emph{Let $\{u_m\}_1^\infty : \Om \sub \R^n \larrow \R$ be a sequence of functions. If
\[
\overline{u} \ \equiv  \ \underline{u}
\]
on $\Om$, where $\overline{u}$, $\underline{u}$ are given by \eqref{9.1},  \eqref{9.2}, then 
\[
u_m \larrow u \ \text{ in }C^0(\Om),
\]
as $m \ri \infty$.
}

\ms

\noi \textbf{Proof.} We denote the common value of $\overline{u}$ and $\underline{u}$ by $u$. Then, since  $\overline{u}$ is upper semicontinuous and $\underline{u}$ is lower semi-continuous, we have that $u$ is continuous. Suppose now for the sake of contradiction that uniform convergence fails on a compact set $K \sub \Om$. Then, for some $\e_0>0$, there would be a subsequence $\{m_k\}_1^\infty$ such that 
\[
\big|u_{m_k}(x_k)\, -\, u(x_k) \big|\, \geq\, \e_0
\]
for some $\{x_k\}_1^\infty \sub K$. By compactness, we may assume that $x_k \ri x$ (along perhaps a further subsequence) as $k\ri \infty$. Assume first that
\[
u_{m_k}(x_k)\, -\, u(x_k) \, \geq\, \e_0
\]
Since $u \in C^0(\Om)$, we have 
\[
\e_0\, \leq\, \underset{k\ri \infty}{\lim{\sup}}\, \big( u_{m_k}(x_k)\, -\, u(x_k) \big) \, \leq\, \overline{u}(x)\, -\, u(x),
\]
which is a contradiction. The case of the opposite inequality follows similarly. Hence, the lemma ensues.  \qed

\ms

\ms

\subsection*{Boundary Jets and Jets relative to non-open sets}

In these notes we have defined Jets $\J^{2,\pm}u(x)$ at all \textbf{interior} points $x$ in an \textbf{open} set $\Om$ and for functions $u : \Om \sub \R^n \ri \R$. However, these pointwise generalised derivatives can be defined for functions
\[
u\ : \ \mO \sub \R^n \ri \R
\]
where $\mO$ is an \textbf{arbitrary} subset of $\R^n$. Notwithstanding, some restrictions are needed in order the results we presented herein to be true for general sets $\mO$ and this restriction is \textbf{local compactness}. The most interesting case of locally compact $\mO$ is
\[
\mO\ =\ \overline{\Om},\ \ \ \Om \sub \R^n \, \text{ open}.
\]
The latter choice is of major importance in the consideration of \textit{nonlinear boundary conditions}.

\ms

\noi \textbf{Definition 1 (Jets relative to general sets).} \emph{Let $u : \mO \sub \R^n \ri \R$ and $x\in \mO$. we define the sets
\begin{align} 
\J^{2,+}_\mO u(x)\, :=\, \Big\{ (p,X) \in \R^n \by \mS(n)& \ \Big|  \ \text{ as }\mO \ni z\ri0, \text{ we have}: \nonumber\\
u(z+x)\, \leq\, u(x) \, &+\, p\cdot z\, +\, \frac{1}{2}X:z\ot z \, +\, o(|z|^2) \Big\},\nonumber
\end{align}
and
\begin{align} 
\J^{2,-}_\mO u(x)\, :=\, \Big\{ (p,X) \in \R^n \by \mS(n)& \ \Big|  \ \text{ as }\mO \ni z\ri0, \text{ we have}: \nonumber\\
u(z+x)\, \geq\, u(x) \, &+\, p\cdot z\, +\, \frac{1}{2}X:z\ot z \, +\, o(|z|^2) \Big\}. \nonumber
\end{align}
We call $\J^{2,+}_\mO u(x)$ the 2nd order Super-Jet of $u$ relative to $\mO$  at $x$ and $\J^{2,-}_\mO  u(x)$ the 2nd order Sub-Jet of $u$  relative to $\mO$ at $x$.}

\ms

It should be quite evident that if $x\in \Om$ and $\Om$ is open, then
\[
\J^{2,\pm}_\Om u(x)\ =\ \J^{2,\pm} u(x).
\]

\noi \textbf{Remark 2.} \textit{Let $\mO= \overline{\Om}$ and $u \in C^2(\overline{\Om})$. The latter means that $u \in C^2(\Om)$ and $Du$, $D^2u$ extend as continuous maps to an open neighbourhood of $\overline{\Om}$. If $x \in \p \Om$, then the characterisation
\[
\J^{2,+}_{\overline{\Om}} u(x)\, =\, \Big\{ \big(Du(x),D^2u(x)+A \big)\ :\ A\geq 0 \Big\}
\]
\textbf{fails}. The reader is invited to verify this claim by attempting to see why the proof fails if $x \in \p \Om$.
}

\ms

Accordingly, we quote without proof the next fact:

\ms

\noi \textbf{Lemma 3 (Jets of smooth functions at regular boundary points.)} \emph{Let $\Om \sub \R^n$ be open, such that $\p \Om$ is a smooth compact hypersurface with exterior normal vector field $\nu$. Let also $u\in C^2(\overline{\Om})$.}

\emph{For an $x\in \p \Om$ fixed, let $T_x \p\Om$ denote the tangent hyperplane at $x$, $\Pi_x : T_x\p \Om \ri \R^n$ the projection on the tangent space and $I\!\!I$ the 2nd fundamental form of $\p \Om$ at $x$ (oriented with respect to $\nu$ and extended on $\R^n$ by zero). Then, we have
\[
(p,X) \in \J^{2,+}_{\overline{\Om}} u(x)
\]
if and only if either
\[
(p,X) \, =\, \big(Du(x),D^2u(x)+A \big)
\]
for some $A\geq 0$, or
\[
p\, =\, Du(x)\, -\, \la \nu, \ \ \ \Pi \big(D^2u(x) -X\big) \Pi \, \leq\,  -\, \la I\!\!I
\]
for some $\la>0$.
 }

\ms

\subsection*{Nonlinear Boundary Conditions}

In these notes we have considered only the Dirichlet problem for fully nonlinear PDE
\[
\left\{
\begin{array}{l}
F\big(\cdot,u,Du,D^2u\big)\,=\,0, \, \text{ in }\Om,\ms\\
\hspace{78pt}u\, =\, b, \, \text{ on }\p \Om,
\end{array}
\right.
\]
in which case the boundary condition is linear and of zeroth order. However, most of our analysis carries over to more general boundary conditions, even fully nonlinear:
\beq \label{9.5}
\left\{
\begin{array}{l}
F\big(\cdot,u,Du,D^2u\big)\,=\,0, \, \text{ in }\Om,\ms\\
B\big(\cdot,u,Du,D^2u\big)\,=\,0, \, \text{ on }\p \Om,
\end{array}
\right.
\eeq
where
\begin{align}
F\ &: \Om \by \R \by \R^n \by \mS(n) \larrow \R, \nonumber\\
B\ &: \p \Om \by \R \by \R^n \by \mS(n) \larrow \R.  \nonumber
\end{align}
However, this can not be done without extra assumption and there are some tricky points. As an illustration, assume that $\p \Om$ is smooth and that $B$ is of first order, that is
\[
B(\cdot,u,Du)\, =\, 0, \ \ \text{ on}\p \Om. 
\]
Suppose \textbf{we attempt} to define \textit{boundary conditions in the viscosity sense} as
\begin{align}
(p,X)\in \J^{2,+}_{\overline{\Om}}u(x) \ \ \Longrightarrow \ \ B(x,u(x),p)\, \geq \, 0, \nonumber\\
(p,X)\in \J^{2,-}_{\overline{\Om}}u(x) \ \ \Longrightarrow \ \ B(x,u(x),p)\, \leq \, 0, \nonumber
\end{align}
when $ x\in \p \Om$. By Lemma 6 above, if $u \in C^2(\overline{\Om})$ then we have
\[
p\, =\, Du(x)\, -\, \la \nu(x)
\]
and hence the above ``definition" is \textit{not compatible with classical satisfaction} of the boundary condition, unless
\[
t \, \mapsto\, B\big(x,r,p- t \nu(x)\big)
\]
is non-decreasing for $t\geq 0$, where $\nu : \p \Om \larrow \R^n$ is the outwards pointing unit normal vector field. Moreover, the degenerate choices 
\begin{align}
F(x,r,p,X)\, &:= \, r-h(x),\ \ |h|\, \geq\, 1,\nonumber\\
B(x,r,p,X)\, &:= \, 0, \nonumber
\end{align}
in \eqref{9.5} imply that \textit{the boundary condition $u\equiv 0$ on $\p \Om$ can not be satisfied classically by the unique solution $u\equiv h \neq 0$ in $\Om$}!  We remedy these problem by introducing the next

\ms

\noi \textbf{Definition 1 (Viscosity formulation of nonlinear Boundary Value Problems).} \emph{Consider the boundary value problem
\[
\left\{
\begin{array}{l}
F\big(\cdot,u,Du,D^2u\big)\,=\,0, \, \text{ in }\Om,\ms\\
B\big(\cdot,u,Du,D^2u\big)\,=\,0, \, \text{ on }\p \Om,
\end{array}
\right.
\]
where $\Om \sub \R^n$ is open and
\begin{align}
F\ &: \overline{\Om} \by \R \by \R^n \by \mS(n) \larrow \R, \nonumber\\
B\ &: \p \Om \by \R \by \R^n \by \mS(n) \larrow \R.  \nonumber
\end{align}
Then, $u\in C^0(\overline{\Om})$ is a viscosity solution to the boundary value problem, when for all $x\in \Om$, we have
\begin{align}
(p,X)\in \J^{2,+}u(x) \ \ \Longrightarrow \ \ F \big(x,u(x),p,X \big)\, \geq \, 0, \nonumber\\
(p,X)\in \J^{2,-}u(x) \ \ \Longrightarrow \ \ F \big(x,u(x),p,X \big)\, \leq \, 0, \nonumber
\end{align}
and for all $x\in \p \Om$, we have
\begin{align}
(p,X)\in \J^{2,+}_{\overline{\Om}}u(x) \ \ &\Longrightarrow \ \ \max\Big\{F \big(x,u(x),p,X \big),B \big(x,u(x),p,X \big) \Big\}\, \geq \, 0, \nonumber\\
(p,X)\in \J^{2,-}_{\overline{\Om}}u(x) \ \ &\Longrightarrow \ \  \min \Big\{F \big(x,u(x),p,X \big),B \big(x,u(x),p,X \big) \Big\} \, \leq \, 0. \nonumber
\end{align}
}

\noi A particular interesting case of boundary operator is the one giving the Neumann boundary condition
\[
B \big(x,u(x),p,X \big)\, :=\, p\cdot \nu(x)\, +\, f(x,r)
\]
where $\nu$ is the outer normal unit vector field of $\p \Om$. One can obtain existence and uniqueness results similar to those obtained for the Dirichlet problem under certain assumptions on $f$, $F$ and $\p \Om$. 

\ms

\noi \textbf{Remark 2 (Generalised normal vectors).} \textit{We conclude this section by giving a ``weak formulation" of the set of (perhaps non-unique or non-existent) exterior normal vectors of  $\, \overline{\Om} \sub \R^n$ at $\hat{x} \in \p \Om$ is
\[
N(\hat{x})\, :=\, \Big\{\nu \in \R^n\ \Big| \ \nu \cdot (x-\hat{x})\, \leq\, o(|x-\hat{x}|),\ \text{ as }\Om \ni x\ri \hat{x} \Big\}.
\]
$N(\hat{x})$ is the normal cone at $\hat{x}$. The figure below show the normal cones at a singular and a regular point of the boundary.
\[
\includegraphics[scale=0.22]{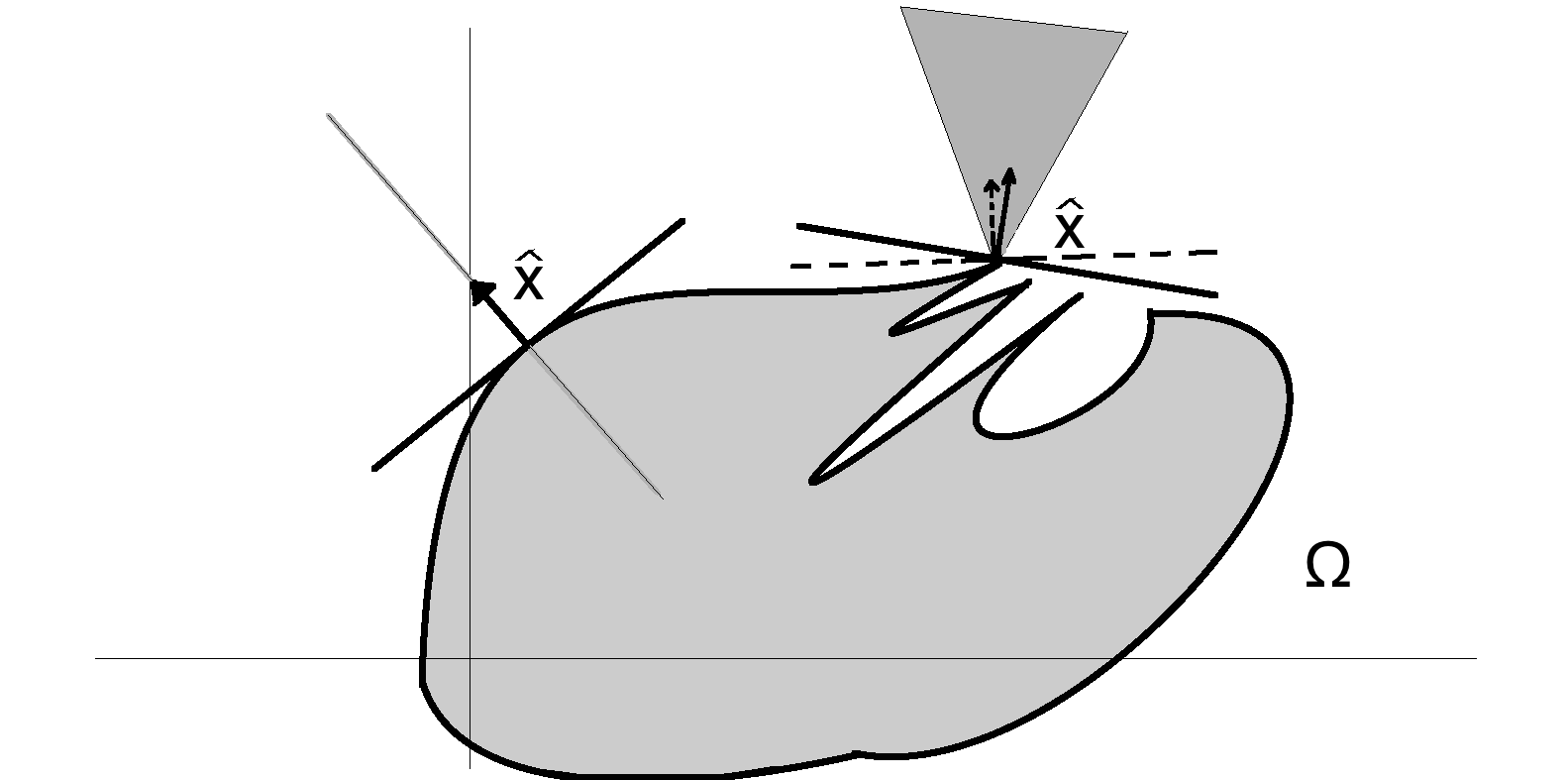}
\]
}

\subsection*{Comparison Principle for Viscosity Solutions without decoupling in the $x$-variable.} In Chapter 6 we proved the comparison principle between sub-/super- solutions for the case of PDE of the form
\[
F(u,Du,D^2u)\, =\, f.
\]
This was just for simplicity in the exposition, since the result is true in the general case of PDE of the form
\[
F(\cdot,u,Du,D^2u)\, =\, 0.
\]
For, we record without proof the following extension the comparison result, which can be easily deduced by slight variations of the arguments of Chapter 6:

\ms

\noi \textbf{Theorem 1 (Comparison Principle for Viscosity Solutions).}  \emph{Let $F\in C^0\big(\Om \by \R \by \R^n \by \mS(n)\big)$, $\Om \Subset \R^n$ and assume that
\begin{align}
X\, \leq \, Y \ \ &\Longrightarrow \ \  F(x,r,p,X)\, \leq\, F(r,p,Y), \nonumber\\
r\, \leq \, s\ \ \ &\Longrightarrow \ \ F(,x,r,p,X)\, \geq\, F(x,s,p,X) \, +\, \ga(s-r), \nonumber
\end{align}
for some $\ga>0$. Assume further that there is an $\om : [0,\infty]\ri [0,\infty]$ with $\om(0^+)=0$, such that, if $\al>0$ and
\[
-\, 3\al\, 
\left[
\begin{array}{cc}
I & 0\\
0 & I
\end{array}
\right] 
                            \, \leq\,  
\left[
\begin{array}{cc}
X& 0\\
0 & -Y
\end{array}
\right]   
                          \, \leq\, 
3\al\,
\left[
\begin{array}{cc}
I & -I\\
-I & I
\end{array}
\right] 
\]
then
\[
F\big(x,r,\al(x-y),X \big) \, -\, F\big(y,r,\al(x-y),Y \big)\, \leq\, \om\big(\al|x-y|^2+|x-y| \big).
\]
Suppose that $u\in C^0(\overline{\Om})$ is a Viscosity Solution of 
\[
F(\cdot,u,Du,D^2u)\, \geq\,0
\]
on $\Om$, and $v\in C^0(\overline{\Om})$ is a Viscosity Solution of 
\[
F(\cdot, v,Dv,D^2v)\, \leq\, 0
\]
on $\Om$. Then, we have
\[
u \, \leq\, v \ \text{ on }\p \Om\ \ \Longrightarrow\ \ u \, \leq\, v \ \text{ on } \overline{\Om}.
\]
}

\ms

\ms

\ms

\noi \textbf{Remarks to Chapter 9.} The material of this Chapter is intended to be an ``eye-opener" for the interested reader, in order to stimulate the interest for the subject by indicating that there exist extensions to ``every possible direction" which might apply to numerous cases and for PDE not treated herein. By no means though the list of extensions is complete. The equivalences on the $\infty$-Laplacian are basically motivated by Crandall \cite{C1} (see also Aronsson-Crandall-Juutinen \cite{ACJ}). Our exposition on tug-of-war games and the $\infty$-Laplacian follows closely the article of Evans \cite{E5}. The first who discovered the connection of $\De_\infty$ to differential games were Perez-Schramm-Sheffield-Wilson \cite{PSSW}. For further material on the $\infty$-Laplacian and related problems we also refer to Juutinen-Lindqvist-Manfredi \cite{JLM2, JLM3}. The material of the last five subsections is taken from the handbook of Crandall-Ishii-Lions \cite{CIL}, in which numerous motivations and references are given.

\newpage

\end{document}